\documentclass[12pt,a4paper]{article}   
  
\usepackage{tensor}

\usepackage{tikz}
\usetikzlibrary{matrix,arrows,decorations.pathmorphing,shapes.geometric}

\usepackage{tikz-cd}

\usepackage{amsthm}
\usepackage{amsmath}
\usepackage{amsfonts} 
\usepackage{amssymb} 
\usepackage{stmaryrd}
\usepackage{enumerate} 
\usepackage{amscd}
\usepackage{graphicx,scalerel}
\usepackage{mathtools}
\usepackage{bbm}
\usepackage{makeidx}
\usepackage{wasysym}

\usepackage[colorlinks=true,citecolor=red,linkcolor=blue]{hyperref} 
\usepackage[capitalise,noabbrev]{cleveref}

\newcounter{Def}[section]

\theoremstyle {definition}
\newtheorem{Proposition}[Def]{Proposition}

\newtheorem{Example}[Def]{Example}

\newtheorem{Remark}[Def]{Remark}
\newtheorem{Lemma}[Def]{Lemma}
\newtheorem{Theorem}[Def] {Theorem}
\newtheorem{Cor}[Def] {Corollary}

\newenvironment{Definition}{ \refstepcounter{Def}\mbox{}\\ 
\noindent\sl\textbf{Definition
\arabic{section}.\arabic{Def}} }
{\vspace{0.3cm}}

\newcommand\Cite[2] {\cite[#1]{#2}}

\def\AA            {{\!A}}

\def\act           {\vartriangleright}
\def\actle         {\negtriangleright}
\def\actlE         {\,{\actle}\,}

\def\actler        {\negtriangleleft}
\def\actleR        {\,{\actler}\,}
\def\actre         {\vartriangleright}
\def\actrE         {\,{\actre}\,}
\def\actrer        {\vartriangleleft}
\def\actreR        {\,{\actrer}\,}
\def\actrp         {\vartriangleright'}
\def\actrP         {\,{\actrp}\,}
\def\alph          {\alpha}  
\def\Alg           {\mathcal Alg}

\def\be            {\begin{equation}}
\def\bearl         {\begin{array}{l}}
\def\bearll        {\begin{array}{ll}}
\def\bimod         {\text-{\rm bimod}}

\def\C             {{\ensuremath\calc}}
\newcommand\cc[1] {\overline{#1}}

\def\calc          {{\mathcal C}}
\def\calcopp       {{\mathcal C}\opp}
\def\calctop       {{\mathcal C}\top}

\def\cald          {{\mathcal D}}

\def\calm          {{\mathcal M}}
\def\caln          {{\mathcal N}}
\def\calo          {{\mathcal O}}

\def\Cat           {{\mathcal Cat}}
\def\cir           {\,{\circ}\,}
\def\coev          {\mathrm{coev}}
\def\coevp         {\widetilde{\mathrm{coev}}}
\def\coevl         {\coev^\ihl}
\def\coevr         {\coev^\ihr}
\def\coHom         {\mathrm{coHom}}

\def\Colon         {:\quad}
\def\comoD           {{\rm comod}\text-}

\def\Cong          {\,{\cong}\,}

\def\deltal        {\delta^\ihl}
\def\deltar        {\delta^\ihr}
\def\distl         {\deltal}  
\def\distr         {\deltar}  

\def\ee            {\end{equation}}
\def\eear          {\end{array}}
\def\End           {\mathrm{End}}
\def\Enumerate     {\def\leftmargini{1.34em}~\\[-1.42em]\begin{enumerate}}
\def\eps           {\varepsilon}
\def\eq            {\,{=}\,}

\newcommand\equ[1] {\stackrel{\eqref{#1}}=}
\def\ev            {\mathrm{ev}}

\def\evl           {\ev^\ihl}
\def\evr           {\ev^\ihr}

\def\ftransport    {flipping transport} 
\def\Ftransport    {Flipping transport} 
\def\Fun           {\mathrm{Fun}}
\def\Hom           {\mathrm{Hom}}

\def\HomC          {\mathrm{Hom}_\calc}
\def\HomM          {\mathrm{Hom}_\calm}

\def\icH           {\underline{\coHom}}
\def\icoHom        {\icH}
\def\icoHomr       {\icH}  
\def\icoHoml       {\icH^{\ihl}}
\def\icoev         {\underline{\mathrm{coev}}}

\def\id            {\mathrm{id}}
\def\Id            {\mathrm{id}}
\def\idA           {\id_\AA}

\def\iDelta        {\underline{\Delta}{}}

\def\iev           {\underline{\mathrm{ev}}}

\def\iH            {\underline{\Hom}}
\def\ihl           {{\mathrm l}}
\def\iHom          {\iH}
\def\iHomr         {\iH}  
\def\iHoml         {\iH^{\ihl}}

\def\ihr           {{\mathrm r}}
\def\imu           {\underline{\mu}{}}

\def\iN            {\,{\in}\,}

\def\Ind           {\mathrm{Ind}}
\def\inv           {^{-1}}
\def\kapp          {\kappa}
\def\kappo         {\overline\kapp}
\def\ko            {{\ensuremath{\Bbbk}}}
\def\GVMod         {\mathrm{GV\text-Mod}}
\def\GVModrela     {\mathrm{GV\text-Mod}^{\actre,\text{lax}}}
\def\GVModreol     {\mathrm{GV\text-Mod}^{\actre,\text{oplax}}}
\def\GVModlela     {\mathrm{GV\text-Mod}^{\actle,\text{lax}}}
\def\GVModleol     {\mathrm{GV\text-Mod}^{\actle,\text{oplax}}}
\def\L             {\mathrm{L}}

\def\Lax           {\mathrm{Lax}}
\def\la            {^{\rm l.a.}}

\def\letore        {\actle \mapsto \actre}

\def\moD           {{\rm mod}\text-}

\def\modC          {{\rm mod}_\calc\text-}
\def\mopp          {^{\tpr\rm opp}}

\def\NAl           {\mho_\AA^\ihl}
\def\NAr           {\mho_\AA^\ihr}
\newcommand\nxl[1] {\\[#1mm]}
\newcommand\Nxl[1] {\\[-1.3em]\\[#1mm]}
\def\ohr           {\reflectbox{$\rho$}}

\def\oloa          {\tpl^{A}}
\def\oloA          {\,{\oloa}\,}
\def\one           {\TI}
\def\opp           {^{\rm opp}}
\def\ota           {\otimes_{\!A}}
\def\otA           {\,{\ota}\,}
\def\oti           {\,{\otimes}\,}
\def\otik          {\,{\otimes_\ko}\,}

\def\pik           {\varpi}     
\def\piv           {\pi}        
\def\Psil          {\Psi}

\def\Psili         {\Psi^{-1}}
\def\Psir          {\Psi'}

\def\Psiri         {\Psi'{}^{-1}_{}}
\def\ra            {^{\rm r.a.}}
\def\ractle        {\negtriangleleft}
\def\ractlE        {\,{\ractle}\,}
\def\ractre        {\vartriangleleft}
\def\ractrE        {\,{\ractre}\,}
\newcommand\rarr[1]{\xrightarrow{~#1~}}
\newcommand\Rarr[1]{\,{\xrightarrow{\,#1\,}}\,}

\def\retole        {\actre \mapsto \actle}

\def\R             {\mathrm{R}}

\def\tev           {\mathrm{ev}}
\def\tevl          {\tev^{\ihl}}
\def\tevr          {\tev^{\ihr}}
\def\TI            {1}
\def\Times         {\,{\times}\,}
\def\To            {\,{\to}\,}
\def\tocong        {\rarr\cong}
\def\Tocong        {\Rarr\cong}
\def\top           {^{\otimes\rm opp}}

\def\tpl           {\negtimes}
\def\tpL           {\,{\tpl}\,}
\def\tplp          {\tpl'}

\def\tpr           {\otimes}
\def\tpR           {\,{\tpr}\,}

\def\tr           {\mathrm{tr}}

\def\vect          {\mathrm{vect}}

\newcommand\void[1]{}

\def\walml         {\widehat\calm^\tpl}
\def\walmr         {\widehat\calm^\tpr}

\def\captimes{\mathbin{\ooalign{$\cap$\cr%
             \hfil\raise0.42ex\hbox{$\scriptscriptstyle\times$}\hfil\cr}}}
\def\cuptimes{\mathbin{\ooalign{$\cup$\cr%
             \hfil\raise0.42ex\hbox{$\scriptscriptstyle\times$}\hfil\cr}}}

\makeatletter
\newcommand*{\relrelbarsep}{.386ex}
\newcommand*{\relrelbar}{%
  \mathrel{%
      \mathpalette\@relrelbar\relrelbarsep }}
\newcommand*{\@relrelbar}[2]{%
  \raise#2\hbox to 0pt{$\m@th#1\relbar$\hss}%
    \lower#2\hbox{$\m@th#1\relbar$}}
\providecommand*{\rightrightarrowsfill@}{%
      \arrowfill@\relrelbar\relrelbar\rightrightarrows }
\providecommand*{\leftleftarrowsfill@}{%
       \arrowfill@\leftleftarrows\relrelbar\relrelbar }
\providecommand*{\xrightrightarrows}[2][]{%
          \ext@arrow 0359\rightrightarrowsfill@{#1}{#2}}
\providecommand*{\xleftleftarrows}[2][]{%
    \ext@arrow 3095\leftleftarrowsfill@{#1}{#2}}
\makeatother

\newcommand{\scirc}{\textup{\CIRCLE}}

\newcommand{\negtimes}{{\textcolor{gray}{\scirc}\mkern-14.03mu{\textcolor{white}{\times}}}}
\newcommand{\negtriangleright}{{\textcolor{gray}{\blacktriangleright}}}
\newcommand{\negtriangleleft}{{\textcolor{gray}{\blacktriangleleft}}}

\begin{document}

\numberwithin{equation}{section}


\begin{flushright}  
   {\sf ZMP-HH/24-11}\\
   {\sf Hamburger$\;$Beitr\"age$\;$zur$\;$Mathematik$\;$Nr.$\;$966}\\[2mm] May 2024
\end{flushright}

\vskip 2.0em     

\begin{center}
	{\bf \Large Grothendieck-Verdier module categories,
        \\[5pt]
        Frobenius algebras and relative Serre functors
	}

\vskip 2.6em

{\large
J\"urgen Fuchs\,$^{a},\,$
Gregor Schaumann\,$^{b},\,$
Christoph Schweigert\,$^{c},\,$
Simon Wood\,$^{d}$
}

\vskip 15mm

 \it$^a$
 Teoretisk fysik, \ Karlstads Universitet \\
 Universitetsgatan 21, \  S\,--\,651\,88\, Karlstad
 \\[9pt]
 \it$^b$
 Mathematische Physik, \ Institut f\"ur Mathematik, \ Universit\"at W\"urzburg \\
 Emil-Fi\-scher-Stra\ss e 31, \ D\,--\,97\,074\, W\"urzburg
 \\[9pt]
 \it$^c$
 Fachbereich Mathematik, \ Universit\"at Hamburg\\
 Bereich Algebra und Zahlentheorie\\
 Bundesstra\ss e 55, \ D\,--\,20\,146\, Hamburg
 \\[9pt]
 \it$^d$
 School of Mathematics, \ Cardiff University\\
 Abacws, \ Senghennydd Road, \ Cardiff, CF24 4AG

\end{center}

\vskip 4.7em

\noindent{\sc Abstract}\\[3pt]
We develop the theory of module categories over a Grothendieck-Verdier category
$\calc$, i.e.\ a  monoidal category with a dualizing object and hence a duality
structure more general than rigidity. Such a category comes with two monoidal 
structures $\tpr$ and $\tpl$ which are related by non-invertible
morphisms and which we treat on an equal footing. Quite generally, non-invertible
structure morphisms play a dominant role in this theory.
In any Grothendieck-Verdier module category $\calm$ we find two important
subcategories $\widehat{\calm}^{\tpl}$ and $\widehat{\calm}^{\tpr}$. The internal 
Hom $A_m \,{:=}\, \iHom(m,m)$ of an object in $\widehat{\calm}^{\tpr}$ that is a
$\calc$-generator is an algebra such that $\moD A_m $ is equivalent to $\calm$ as 
a module category. We also introduce a partially defined relative Serre functor
which furnishes an equivalence $S\colon \widehat{\calm}^{\tpr}
\Rarr~ \widehat{\calm}^{\tpl}$. Any isomorphism $m\Rarr\cong S(m)$ endows 
the algebra $A_m$ with the structure of a Grothendieck-Verdier Frobenius algebra.

	    \newpage
	    \tableofcontents
	    \newpage


\section{Introduction}

Monoidal categories arise in a wide variety of contexts, for instance as cobordism
categories and as representation categories of diverse
algebraic structures. Accordingly, monoidal 
categories have found numerous applications, ranging from representation theory to
quantum topology and logic. Rigidity is a natural duality property for an object in
a monoidal category. A rigid monoidal category is a monoidal category in which every
object has rigid duals; such categories are pervasive in quantum topology.

To motivate the results of this paper, let us consider the following subclass of 
rigid categories: Let $\ko$ be an algebraically closed field.
A \emph{finite tensor category} over $\ko$ is a $\ko$-linear abelian rigid monoidal
category that is finite, i.e.\ is equivalent, as an abelian category, to the category
of finite-dimensional modules over a finite-dimensional $\ko$-algebra. The class of
finite tensor categories includes e.g.\ categories of finite-dimensional modules
over finite-dimensional Hopf algebras; for these, rigidity is rather directly
inherited from basic properties of finite-di\-men\-sional vector spaces. 
Finite tensor categories are well understood \cite{EGno}. In particular there is a
satisfactory theory of exact module categories over finite tensor categories. This
theory is deeply linked with rich algebraic notions like Drinfeld centers and 
relative Serre functors. In fact, the bicategory of exact module categories is as
indispensable for the understanding of a finite tensor category as the category
of modules is for the understanding of a ring.

On the other hand, for a general monoidal category rigidity is a rather restrictive
property. For instance, it forces the tensor product of a finite tensor category to
be an exact functor. It should  therefore not come as a surprise that many concrete
monoidal categories arising in representation theory, algebraic geometry or linear
logic are not rigid. However, a lot of interesting tensor categories still exhibit
the more general duality structure of a star-autonomous or \emph{Grothendieck-Verdier}
duality. Roughly speaking, a Grothendieck-Verdier category is a monoidal category
$(\calc,\otimes)$ together with the additional structure of a \emph{dualizing object}
$K\iN\calc$ (see Definition \ref{def:thedef}), which entails a contravariant duality
functor $G\colon \calc\Rarr~\calcopp$ that is an anti-equivalence. 
Examples of Grothendieck-Verdier categories include suitable representation 
categories of Hopf algebroids \cite{alleR} and of a large class of vertex algebras
\cite{alsW}, as well as categories of finite-dimensional bimodules over
finite-dimensional algebras, see e.g.\ the discussion in \cite{fssw}.

Apart from the fact that they naturally arise as representation categories for large
classes of algebraic structures, Grothendieck-Verdier categories also have several 
other appealing features. In particular, while their tensor product is not necessarily
exact, the structure is still sufficiently restrictive so as to guarantee the 
existence of internal Homs. Owing to the fact that the duality functor $G$ of a
Grothendieck-Verdier category $(\calc,\tpr,K)$ is not necessarily monoidal, $\calc$
naturally comes with a second tensor product $\tpl$ that is left exact.
This also explains why such categories have been considered in categorical approaches
to logic, which require two monoidal operations ``and'' and ``or'' that are related 
by a negation $G$. Each of the two tensor products comes with associator isomorphisms
that obey a respective pentagon identity. In addition they are linked by 
``mixed associators'', called \emph{distributors}, which, unlike ordinary associators,
are not necessarily isomorphisms, but still obey pentagon identities. Structure
morphisms that are not isomorphisms will be ubiquitous in our discussion.

In this paper we develop a theory of module categories over
Grothendieck-Verdier categories. We are motivated by the desire to reach a deeper
structural understanding of Grothendieck-Verdier categories, as well as by potential 
applications of module categories in the construction of two-dimensional conformal 
field theories. (The latter application actually needs at least pivotal module
categories over pivotal Grothendieck-Verdier categories.)

\medskip

Let us summarize some of our results that we regard as most significant. We consider
left module categories $(\calm,\actre)$ over a Grothendieck-Verdier category 
$(\calc,\tpr,K)$ for which the functors   
  \be
  c \actrE {-} \Colon \calm \rarr~ \calm \qquad \text{and} \qquad
  {-} \actrE m \Colon \calc \rarr~ \calm
  \ee
have a right adjoint for every $c\iN\calc$ and $m\iN\calm$. (The action of $\calc$ on
itself, the regular module, gives an example of such a module category.)
If $\calc$ is even rigid, then this reduces to the class of exact module categories 
considered in \cite{EGno}. We show in Proposition \ref{proposition:actle} that the
$\tpr$-module category $(\calm,\actrE)$ also comes with a natural structure of a
module category $(\calm,\actlE)$ over the monoidal category $(\calc,\tpl)$ with
monoidal structure given by the left exact tensor product on $\calc$.
In fact, 
we take it as a guiding principle of our analysis that the two monoidal structures 
of a Grothendieck-Verdier-ca\-tegory should be treated on an equal footing. In
particular, constructions that are defined using one of the two monoidal structures,
such as module categories and module functors, should exhibit corresponding
features with respect to the other monoidal structure as well.

Next we study module functors. Given two module categories $\calm$ and $\caln$ over
the same Gro\-then\-dieck-Verdier category $\calc$, these are functors 
$F\colon \calm \Rarr~ \caln$ with additional constraint data. In our setting it is 
not reasonable to require that these constraints are isomorphisms. This leads us to
consider four classes of module functors, being lax and oplax for either of the two
actions $\actle$ and $\actre$, respectively, constituting four a priori different
bicategories. In Theorem \ref{Thm:bicat-GV-mod} we show that two of these 
bicategories can be naturally identified. Our insights about module functors allow us
to introduce, in Definition \ref{Definition:module-distrib}, \emph{module distributors}
which relate the two actions $\actlE$ and $\actrE$. In Section \ref{sec:adjoints}, 
we use profunctors and adjoints to exhibit various relations between different 
classes of module functors; they are summarized in the table \eqref{eq:table}.
Our theory allows 
us to exhibit a plentitude of six pentagon diagrams for Grothendieck-Verdier module 
categories, and 32 pentagon diagrams for Grothendieck-Verdier bimodule categories. 

It is well-known that Grothendieck-Verdier categories admit evaluation and 
coevaluation morphisms. Our constructions allow for a concise proof that these obey
snake relations that involve the distributors. We finally reveal in Proposition
\ref{Proposition:weak-str-internal-hom} in which way internal Hom (and the coHom)
functors possess natural structures of module functors. 

In Section \ref{sec:Frob} we turn our attention to the study of algebras and 
coalgebras in Grothen\-dieck-Ver\-dier categories. These are straightforward to define, 
and it is immediate that internal (co)Homs of Grothendieck-Verdier module categories 
provide examples. We introduce three possible definitions of the notion of a Frobenius
algebra and show in Theorem \ref{thm:4.14} that they are equivalent. Afterwards we
turn to the following issue: While it is obvious that for any algebra $A$ in a 
Grothendieck-Verdier category $\calc$ the category $\moD A$ of right $A$-modules is a
$\calc$-module category, it turns out to be a much more subtle question whether a 
given $\calc$-module category $\calm$ can be represented as a category of modules 
over some algebra in $\calc$. And indeed generically this is not possible; we call
those module categories which \emph{can} be represented in this way 
\emph{algebraic}.

When investigating this question we are led to the Definition 
\ref{definition:Admissible} of two interesting subcategories of $\calm$: one of
these is the full subcategory $\widehat\calm^\tpr$ of $\tpr$-\emph{admissible}
objects, i.e.\ of those objects $m\iN\calm$ for which the functor
$\iHom(m,-) \colon \calm\Rarr~\calc$ is a strong $\actre$-module functor and 
has a right adjoint. Our principle to treat the two monoidal structures $\tpr$ and
$\tpl$ on the same footing leads us to also consider the corresponding structures for
the $\actlE$-action. We thus obtain another subcategory $\widehat\calm^\tpl$ of 
objects $m\iN\calm$ for which the functor $\icoHom(m,-) \colon \calm\Rarr~\calc$
is a strong $\actle$-module functor and has a left adjoint.
  
The categories $\widehat\calm^{\tpr,\tpl}$ are not abelian, but they are still nicely
compatible with the monoidal structure on $\calc$. As we show in Proposition
\ref{proposition:admissible-module}, by restriction of the action of $\calc$ on 
$\calm$ the category $\widehat{\calm}^{\tpr}$ is a left 
$\widehat{\calc}^{\tpr}$-module category, while $\widehat{\calm}^{\tpl}$ is a left 
$\widehat{\calc}^{\tpl}$-module category. In particular, the subcategories
$\widehat{\calc}^{\tpr}$ and $\widehat{\calc}^{\tpl}$ of a Grothendieck-Verdier
category $\calc$ are monoidal subcategories.
With these subcategories at hand, we are in a position to characterize, in 
Proposition \ref{Prop:i(co)Hom-equivalence} and Theorem \ref{Thm:i(co)Hom-equivalence}, 
algebraic module categories: objects
$m\iN\widehat\calm^\tpr$ which are $\calc$-generators are precisely those for which
the internal End $A_m \eq \iHom(m,m)$ is an algebra such that $\moD A_m$ is 
equivalent to $\calm$. A dual statement holds for $\calc$-cogenerators
$m\iN\widehat\calm^\tpl$, the coalgebra $C_m \eq \icoHom(m,m)$ and $C_m$-comodules.

We consider the existence of the subcategories $\widehat{\calm}^{\tpr}$ and 
$\widehat{\calm}^{\tpl}$ as one of the main insights of this paper. 
They constitute a profound structure with important consequences. In particular,
they enter into the definition of \emph{relative Serre functors}
which we introduce in Section \ref{sec:iHom+repfun}: two functors
  \be
  S \Colon \widehat{\calm}^{\tpr} \rarr{~} \calm \qquad\text{and}\qquad
  \widetilde{S} \Colon \widehat{\calm}^{\tpl} \rarr{~} \calm \,.
  \ee
such that $\iHomr(n,Sm) \Cong G (\iHomr(m,n))$ and
$\iHomr(\widetilde{S}m,n) \Cong G (\iHomr(n,m))$. Theorem \ref{thm:Serre=equiv} then
asserts that $S$ and $\widetilde{S}$ provide an equivalence of categories between
$\widehat{\calm}^{\tpr}$ and $\widehat{\calm}^{\tpl}$. This allows us to generalize
the insights of \cite{shimi20} to show that every choice of isomorphism 
$p \colon m \Rarr~ Sm$ in $\calm$ induces the structure of a Frobenius algebra
on the algebra $\iHom(m,m)$.

In the present paper we discuss pivotal Grothendieck-Verdier categories only in the
following contexts: In Section \ref{sec:symmFrob} we characterize \emph{symmetric}
Frobenius algebras in pivotal Grothendieck-Verdier categories, and at the end of
Section \ref{sec:Serre+Frob} we explain how rigid duals and Grothendieck-Verdier
duals are related in pivotal categories. In view of potential applications to modular
functors and conformal field theories, pivotal Grothendieck-Verdier categories
deserve a further study.


\section{Preliminaries}

We assume that the reader is familiar with the concept of rigid duality in monoidal
categories. In the present section we review the notion of Grothendieck-Verdier
duality, contrast it to rigid duality, and recall concepts and results from the
theory of module categories over monoidal categories.

In our discussion we consider simultaneously two types of categories -- linear
and not necessarily linear ones. More explicitly, if no extra
assumptions are made, then categories and functors will be arbitrary. If instead
we deal with a linear category, then we assume that $\ko$ is an arbitrary but
fixed field and that the category is $\ko$-linear without any further finiteness 
assumptions: a functor between linear categories is assumed to be linear, 
natural transformations are linear maps, and for a
module category over a linear monoidal category the action functor is assumed to be 
bilinear. Following \Cite{Defs.\,1.8.5\,\&\,1.8.6}{EGno}, a \emph{finite category}
is for us a linear category that  is equivalent as an abelian category to the 
category of finite-dimensional modules over a finite-dimensional \ko-algebra.


\begin{Definition} \label{def:thedef}
\Enumerate
 \item
An object $K$ in a monoidal category $(\calc,\tpr,\TI,\alph,l,r)$ is said to be a
\emph{dualizing object} if, for every $y\iN \calc$, the functor
$x \,{\mapsto}\, \Hom(x \tpR y,K)$ is representable by some object $Gy \iN \calc$
and the so defined contravariant functor $G\colon \calc \Rarr~ \calc$ is an 
anti-equivalence, so that there are isomorphisms
  \be
  \pik_{x,y}^{} \Colon \Hom(x\tpR y,K) \rarr\cong \Hom(x,Gy)
  \label{eq:GV1}
  \ee
natural in $x,y \iN \C$. $G$ is called the \emph{duality functor with respect to $K$}.
 \item
A \em{Grothendieck-Verdier category} -- or \em{GV-category} for short -- is a monoidal
category together with a choice of a dualizing object $K\iN\calc$.
\end{enumerate}
\end{Definition}

The assignment of the functor $G$ on a morphism $f\colon x \Rarr~ y$ is obtained by 
noting that the pullback induces the natural transformation in the left
column of the diagram 
  \be
  \begin{tikzcd}[row sep=2.2em]
  \Hom(- \tpR y,K) \ar{r}{\pik_{-,y}} \ar{d}[swap]{(\id \tpr f)^{*}}
  & \Hom(-, Gy) \ar{d}{Gf_{*}}
  \\
  \Hom(- \tpr x, K) \ar{r}{\pik_{-,x}} & \Hom(-, Gx)
  \end{tikzcd}
  \label{eq:def:Gonmorphisms}
  \ee
Transporting this natural transformation to the one in the right column of
\eqref{eq:def:Gonmorphisms} defines, by the Yoneda lemma, the morphism $Gf$.

In general, $G(y)\tpR G(x)$ is not isomorphic to $G(x\tpR y)$. However, the covariant
functor $G^2$ comes with a canonical monoidal structure \cite[Prop.\,5.2]{boDr}.
The choice of $K$ is \emph{structure}. The following result
\cite[Prop.\,2.3]{boDr} clarifies the freedom given by this structure:
  
\begin{Proposition}\label{prop:invdual}
Let $(\calc ,K)$ be a Grothendieck-Verdier category with duality functor
$G \,{\equiv}\, G_K$.
 \\[-1.4em]
\Enumerate
 \item
The functor $G$ is an anti-equivalence between the full subcategory of  
invertible objects $U$ in $\calc$ and the full subcategory of dualizing objects.
Indeed, $G$ satisfies $G(U) \eq K\tpR U^{-1}$ for any invertible object $U \iN \C$.
 \\
Analogous statements hold for the functor $G\inv$, with $G^{-1}(U) \eq U^{-1}\tpR K$.
 \item
If an object $U\iN\calc$ is invertible, then $G^2(U)$ is invertible as well,
and one has a canonical isomorphism $K \tpR U^{-1} \Rarr\cong (G^2U)^{-1}\tpR K $.
\end{enumerate}
\end{Proposition}

The duality functor $G$ allows for the following definition, which is central to 
our work. 

\begin{Definition}
The functor $\tpl \colon \calc\Times \calc \Rarr~ \calc$ is defined by the mapping 
  \be
  x\tpL y:= G^{-1} (Gy\tpR Gx) 
  \label{eq:def:seccirc}
  \ee
on objects, and analogously on morphisms.
\end{Definition}

This functor provides a second monoidal structure on $\calc$:

\begin{Proposition}
Let  $(\calc ,\tpr , \TI ,\alph, l^\tpr, r^\tpr, K)$ be a GV-category.
Then the functor $\tpl \colon \calc\Times \calc\To \calc$ defined by
\eqref{eq:def:seccirc} together with unit constraints 
  \be
  l_x^\tpl := G^{-1}((r_{Gx}^\tpr)^{-1}) \qquad\text{and}\qquad
  r_x^\tpl := G^{-1}((l_{Gx}^\tpr)^{-1})
  \ee
and associator
  \be
  \alph^\tpl_{x, y, z} := G^{-1} ( \alph^{-1}_{Gz,Gy,Gx}) 
  \ee
endows $\calc$ with the structure of a monoidal category 
$(\calc,\tpl,K,\alph^\tpl,l^\tpl,r^\tpl)$ with monoidal unit $K$.
\end{Proposition}

One may also define a functor $\tplp\colon \C\Times\C \Rarr~ \C$ by
$x\tplp y \,{:=}\, G (G^{-1}y\tpR G^{-1}x)$. This is not independent, however:
The $\tpr$-monoidal structure on the functor \(G^2\) implies a canonical 
identification of \(\tpl\) and \(\tplp\) \Cite{Sect.\,4.1}{boDr}.
The functors $\tpr$ and $\tpl$ are isomorphic if and only if $G$ is monoidal, and
thus in particular if \C\ is rigid.
Also note that the functor $G^{2}$ is monoidal both for $\tpr$ and for $\tpl$.

If we change the choice of dualizing object from $K$ to 
$\widetilde{K} \,{:=}\, g \tpR K$ for an invertible object $g \iN \calc$, then
the monoidal functor $G^{2}$ gets replaced by $\widetilde{G}^{2}$ defined as 
  \be
  \widetilde{G}^{2}(y) := g \tpr G^{2}(y) \tpr g^{-1}
  \label{eq:widetildeG**2}
  \ee
for $y\iN\calc$. 
To see this, just note that the functor $G^{2}$ is determined by the isomorphism
$\Hom(x \tpr y, K) \Cong \calc(G^{2}y \tpR x, K)$ for all $x,y \iN \calc$
and that we have
  \be
  \begin{aligned}
  \Hom(x \tpR y, \widetilde{K}) &= \Hom(x \tpR y, g \tpR K)
  \cong \Hom(g^{-1} \tpR x \tpR y, K) 
  \nxl1
  & \cong \Hom(G^{2}y \tpR g^{-1} \tpR x, K)
  \nxl1
  & \cong \Hom(g \tpR G^{2}y \tpR g^{-1} \tpR x, g \tpR K)
  = \Hom(\widetilde{G}^{2}y \tpR x, \widetilde{K})
  \end{aligned}
  \ee
for $x,y \iN \calc$.

      \bigskip

For $a,b,c\iN\calc$ we have    
  \be      
  \Hom_\calc(a\tpR b, c) \cong \Hom_\calc(a, c\tpL Gb)
  \label{eq:195}
  \ee
and
  \be
  \Hom_\calc(a\tpR b, c) \cong \Hom_\calc(b, G^{-1}\!a\tpL c) \,.
  \label{eq:196}
  \ee
Informally, $G$ behaves like a right duality and $G^{-1}$ like a left duality 
if we take the $\tpr$-tensor product in the domain and the $\tpl$-tensor product in 
the codomain of a morphism.

In general the two monoidal structures are different. However, there are mixed 
associators, called \emph{distributors}, which relate them. These are coherent 
morphisms $(a \tpL b) \tpR c \Rarr~ a \tpL (b \tpR c)$ and
$a \tpR (b \tpL c) \Rarr~ (a \tpR b) \tpL c$, for $a,b,c \iN \calc$; they are not 
necessarily isomorphisms. In Section \ref{sec:app:dist,dual} we will examine the 
distributors in detail and in particular prove the coherence diagrams satisfied by them.

\begin{Example} \label{Example:comm-alg}
Let $A$ be a finite-dimensional \ko-algebra. The category $A\bimod$ of 
finite-di\-men\-si\-onal $A$-bimodules is a GV-category with the usual tensor product 
over $A$ as the right exact tensor product and with the dualizing functor given by 
$G(m) \eq m^{*}$ for $m \iN A\bimod$, with $m^{*}$ denoting the linear dual \cite{fssw}.
If $A$ is in addition commutative, then the category $\moD A$ of finite-dimensional
right $A$-modules is a full subcategory of $A\bimod$ by regarding a right module as
a bimodule in the natural way. Moreover, the GV-structure on $A\bimod$ then induces
the structure of a GV-category on $\moD A$.
 \\[2pt]
Other examples of GV-categories include representation categories of vertex operator
algebras \cite{alsW} and of Hopf algebroids \cite{alleR}.
\end{Example}

\begin{Remark} \label{rem:K=injective}
In an abelian GV-category $\calc$, $x \iN \calc$ is projective if and only if $G(x)$ 
is injective: $x$ is projective if and only if $\Hom(x,-)$ is exact, which is the 
case iff $\Hom(G(-),G(x))$ is exact, which (since $G$ is an equivalence) is the case 
iff $\Hom(-,G(x))$ is exact, which is the case if and only if $G(x)$ is injective.
\end{Remark}


Crucial for our work is the connection between module categories and categories
of modules.

\begin{Definition} \label{def:mod}
Let $\calc$ be a monoidal category. A \emph{left module category over $\calc$}
is a category $\calm$ together with a functor 
  \be
  \actre \Colon \calc \Times\calm \rarr~ \calm \,,
  \ee
called the \emph{action functor}, together with a natural family of
isomorphisms $c \actrE (d \actrE m) \Rarr\cong (c \tpr d) \actrE m$
for all $c,d \iN \calc$ and all $m \iN \calm$ and isomorphisms
$\TI \actrE m \Rarr\cong m$ for all $m \iN \calm$, which obey the appropriate
pentagon and triangle constraints. 
\end{Definition}

Analogously one considers right module categories over $\calc$. For $\cald$ another
monoidal category, a $(\calc,\cald)$-bimodule category is a left $\calc$- and a right
$\cald$-module category with an additional family of coherent isomorphisms that 
interchange the ordering of the two actions.

\begin{Definition} \label{def:alg,coalg}
A $($unital, associative$)$ \emph{algebra} in a monoidal category
$(\calc,\tpr)$ is an algebra object in $\calc$, i.e.\ a triple
$(A,\mu,\eta)$ with $A\iN\calc$, $\mu\iN\Hom(A\tpR A,A)$ and $\eta\iN\Hom(\one,A)$
satisfying $($including the associator 
$\alph_{\AA,A,A} \colon (A \tpR A) \tpR A \Rarr~ A \tpR (A \tpR A$ and unitors
$l_\AA\colon \TI\tpR A \Rarr~ A$ and $r_\AA\colon A\tpR\TI \Rarr~ A)$ 
  \be
  \mu \circ (\mu \tpR \idA) = \mu \circ (\idA \tpR \mu) \circ \alph_{\AA,A,A}
  \label{eq:GVass}
  \ee
as morphisms in $\Hom((A\tpR A) \tpR A,A)$ and
  \be
  \mu \circ (\eta \tpR \idA) \circ l_\AA\inv = \idA
  = \mu \circ (\idA\tpR \eta) \circ r_\AA\inv
  \label{eq:GVuni}
  \ee
as morphisms in $\End(A)$.  
 \\[2pt]
Dually, a $($counital, coassociative$)$ \emph{coalgebra} in
$(\calc,\tpr,\one,\alph,l,r)$ is a triple $(C,\Delta,\eps)$ with
$C\iN\calc$, $\Delta\iN\Hom(A,A\tpR A)$ and $\eps\iN\Hom(A,K)$ such that
$(\Delta \tpR \idA) \cir \Delta \eq (\alph_{A,A,A})\inv {\circ}\, 
(\idA \tpR \Delta) \cir \Delta$ as morphisms in $\Hom(A,(A\tpR A)\tpR A)$ and
$l_A \cir (\eps \tpR \idA) \cir \Delta \eq \idA
\eq r_A \cir (\idA \tpR \eps) \cir \Delta$ as morphisms in $\End(A)$.
\end{Definition}

Similarly one defines for an (co)algebra in $\calc$ the notion of an $A$-(co)module. 
It is well known that for any algebra $A \iN \calc$ the category $\moD A$ of right
$A$-modules is a left $\calc$-module category. The category $\comoD C$ of right 
$C$-comodules over a coalgebra $C \iN \calc$ is a left $\calc$-module category, too.


\section{GV-module categories}

In this section we develop the theory of module categories over GV-categories.
According to our guiding principle, module categories and module functors should 
exhibit similar features with respect to both monoidal structures of a GV-category.


\subsection{Module categories and module functors for GV-categories}

\begin{Definition} \label{def:gvmod}
Let $\calc$ be a GV-category. A \emph{left GV-module category} is a
left module category over $(\calc,\tpr)$ with action functor
$\actre \colon \calc \Times\calm \Rarr~ \calm$, such that the functors
  \be
  c \actrE {-} \Colon \calm \rarr~ \calm \qquad \text{and} \qquad
  {-} \actrE m \Colon \calc \rarr~ \calm
  \label{eq:action-functors-GVMod}
  \ee
have a right adjoint for all objects $c \iN \calc$ and $m \iN \calm$.
\end{Definition}

\begin{Remark}
This definition could be formulated for an arbitrary monoidal category $\calc$, as 
the GV-structure of $\calc$ is not used at all. Nevertheless the separate terminology
`GV-module category' is legitimate: it is justified by the results that we will
obtain below, such as Proposition \ref{proposition:actle}. Also note that a functor
that has a right adjoint preserves colimits and is thus in particular cocontinuous.
Demanding the existence of adequate adjoints allows us to work simultaneously in the
set-theoretic and linear setting, while in the latter case it generalizes 
certain exactness assumptions, see Lemma \ref{lem:mod=GVmod}.
\end{Remark}

\begin{Remark} \label{rem:twistbyG2}
Twisting an action of a monoidal category by a monoidal endofunctor yields again
an action. In the case of a GV-category we can in particular twist the action by
any even power of the functor $G$. Thus with any GV-module category over a 
GV-category $\calc$ there comes a whole family of GV-module categories over $\calc$.
\end{Remark}

The internal Hom functor $\iHomr$ is the right adjoint of the functor $- \actrE m$,
which exists by the assumptions in Definition \ref{def:gvmod}:
  \be
  \Hom_\calm(c\actrE m,m') \cong \Hom_\calc(c,\iHomr(m,m')) 
  \ee
for $c\iN\calc$ and $m,m'\iN\calm$. These
will be examined in detail in Section \ref{sec:internal-homs-cohoms}.

In case $\calc$ is a rigid finite category, then for a module category $\calm$
one requires exactness in the first variable (see Definition 7.3.1 of
\cite{EGno}); exactness in the second variable is then
automatic (Exercise 7.3.2 in \cite{EGno}). 
Indeed, existence of the respective adjoints follows in this case as well:

\begin{Lemma} \label{lem:mod=GVmod}
If $\calc$ is a rigid finite category, then a GV-module category over $\calc$ is
the same as a $\calc$-module category in the sense of \Cite{Def.\,7.3.1}{EGno}, i.e.\
the action functor $\act \colon \calc \Times \calm \Rarr~ \calm$ is exact in 
the first variable.
\end{Lemma}

\begin{proof}
Let $\calm$ be a module category, in the sense of \cite{EGno}, over a rigid finite
category. A functor between finite categories has a right/left adjoint if and only if
it is right/left exact (see e.g.\ Corollary 2.3 in \cite{fScS2}).
Thus if $-\actrE n \colon \calc \Rarr~ \calm$ is exact, then it
has in particular a right adjoint. Hence by \Cite{Exc.\,7.3.2}{EGno}, the functors 
$c \actrE - \colon \calm \Rarr~ \calm$ have a right adjoint.
Conversely, if the module category $\calm$ is a GV-module category, then the functors 
$-\actrE n$ have a right adjoint and are thus right exact. Further, the isomorphisms
  \be
  \Hom(m, c \actrE n) \cong \Hom(c^{\vee} \actrE m, n )
  \cong \Hom(c^{\vee},\iHomr(m,n)) \cong \Hom({}^{\vee}\iHomr(m,n),c) 
  \label{eq:alsoleftadj}
  \ee
for $m,n\iN\calm$ and $c\iN\calc$ show that $-\actrE n$
also has a left adjoint and is thus also left exact. 
(Note that in \eqref{eq:alsoleftadj} rigidity is only used in the first step.)
\end{proof}

For a GV-module category, by definition for every $c \iN \calc$ the endofunctor
$\rho^\tpr_\calm(c)\,{:=}\, \,c\actrE{-}$ of $\calm$ has a right adjoint. In the
special case $\calm=\calc$, i.e.\ the regular module category of $\calc$ acting 
on itself by $\tpr$, according to the isomorphisms \eqref{eq:196} we already have 
such a right adjoint, namely $G^{-1}(c) \tpL -$; by convention we will work with this
particular right adjoint. In all other cases we fix some right adjoint, which we
denote by $H_c\colon \calm\To\calm$, i.e.\ we have a natural family
  \be
  \Hom_\calm(c\actrE m,n) \tocong \Hom_\calm(m,H_c(n))
  \label{eq:definingadjunction}
  \ee
of isomorphisms for $m,n\iN\calm$. The endofunctor $H_c$
defines via $c \,{\xmapsto{~~}}\, H_{Gc}(-)$ a functor
  \be
  \rho^\tpl_\calm\Colon \calc \rarr~ \Fun_{\rm l.e.}(\calm,\calm)
  \ee
from $\calc$ to the category of left exact endofunctors of $\calm$, i.e.\
$\rho^\tpl_\calm(-)(m) \eq H_{Gc}(m) \colon \calc \Rarr~ \calm$.
The functor $H_{Gc}$ appearing here furnishes an action of the monoidal category 
$(\calc,\tpl)$ on $\calm$:

\begin{Proposition} \label{proposition:actle} 
Let $\calm$ be a left GV-module category over a GV-category $(\calc,\tpr,K)$. Then the
bifunctor $\actle \colon \calc\Times\calm \Rarr~ \calm$ defined by
  \be
  c\actlE m := H_{Gc}(m)
  \ee
for $c\iN\calc$ and $m\iN\calm$
has a left adjoint (and thus preserves limits, and is in particular continuous) 
in each variable and defines a left module category structure over $(\calc,\tpl)$.
\end{Proposition}

\begin{proof}
The left adjoint of $c \actlE {-}$ is by definition the endofunctor $Gc \actrE {-}$.
For the adjoint in the second variable we use the family of isomorphisms  
  \be
  \begin{aligned}
  \Hom_\calm(m,c\actlE n) & = \Hom_\calm(m,H_{Gc}(n)) \cong \Hom_\calm(Gc\actrE m,n)
  \Nxl1
  & \cong \Hom_\calc(Gc, \iHomr(m,n)) \cong  \Hom_\calc( G^{-1}\iHomr(m,n),c)
  \end{aligned}
  \label{eq:Ginv-iHom}
  \ee
to conclude that $- \actlE n$ has $ G^{-1}\iHomr(-,n)$ as left adjoint.
The module category structure is established as follows. The module constraint of
$\actle$ is the natural family of isomorphisms that are obtained by the 
Yoneda lemma from the isomorphisms
  \be
  \begin{aligned}
  \Hom_\calm(m,b\actlE (c\actlE n))
  & \cong \Hom_\calm(Gb\actrE m,c\actlE n)
  \cong \Hom_\calm(Gc\actrE (Gb\actrE m), n)
  \Nxl1
  & \cong \Hom_\calm((Gc\tpR Gb)\actrE m), n)
  \Nxl1
  & \cong \Hom_\calm( m,G^{-1}(Gc\tpr Gb)\actlE n)
  = \Hom_\calm( m,(b\tpL c)\actlE n)
  \end{aligned}
  \ee
for $m,n\iN\calm$ and $b,c\iN\calc$.
Moreover, the Yoneda embedding transports the pentagon relation for the module
constraint of $\actre$ to the pentagon relation for the module constraint of $\actle$.
\end{proof}

By construction we have the adjunction formula
  \be
  \Hom_\calm(m,G^{-1}c\actlE m')\cong \Hom_\calm(c\actrE m,m') \,.
  \label{eq:HomM.cm.n}
  \ee
The resulting module structure depends on the choice of right adjoints. However, we have

\begin{Lemma} \label{lemma:relating-choices}
The module categories  over $(\calc, \tpl)$ for different choices of right adjoints 
\eqref{eq:definingadjunction} are canonically equivalent. 
\end{Lemma}

\begin{proof}
Let $\calm$ be a left GV-module category over $\calc$, and let $H_{Gc}$ and
$\widetilde{H}_{Gc}$, for $c \iN \calc$, be two families of right adjoints
\eqref{eq:definingadjunction}, which result in two module categories by $\calm$ and
$\widetilde{\calm}$ with respective module structures $c\actlE m \eq H_{Gc}(m)$ and
$c \,{\widetilde{\actle}}\,m \eq \widetilde{H}_{Gc}(m)$. The following considerations
show that the identity functor $\id \colon \calm \Rarr~ \widetilde{ \calm}$
has a canonical structure of a module functor and thus provides an equivalence of
module categories:
By the uniqueness of right adjoints we get for every $c \iN \calc$ a unique natural
isomorphism $\varphi_{c} \colon H_{Gc} \Rarr~ \widetilde{H}_{Gc}$ with components
$\varphi_{c,m}\colon c \actlE m \eq H_{Gc}(m) \Rarr~ \widetilde{H}_{Gc}(m) \eq
c \,{\widetilde{\actle}}\, m$ for all $m \iN \calm$. These morphisms are natural in
$m \iN \calm$ as well as natural in $c \iN \calc$: For any morphism 
$f\colon c \Rarr~ d$ in $\calc$ we are given a natural transformation
$\psi_{f} \colon (G(f) \actrE {-}) \colon (G(d) \actrE {-}) \,{\xRightarrow{~\,}}\,
(G(c) \actrE {-})$, which defines two natural transformations 
  \be
  \psi_{f}^{*} \Colon H_{Gc} \xRightarrow{~~} H_{Gd} \qquad \text{and} \qquad
  \widetilde\psi_{f}^{*} \Colon \widetilde{H}_{Gc} \xRightarrow{~~} \widetilde{H}_{Gd}
  \ee
between the respective right adjoints. By a general fact about the isomorphism
relating different right adjoints,
for all $c, d \iN \calc$ and all $f\colon c \Rarr~ d$ we have the equality
  \be
  \widetilde{\psi_{f}}^{*} \cdot \varphi_{c} = \varphi_{d} \cdot \psi_{f}^{*}
  \ee
between vertical composites of natural transformations. This directly implies that
the isomorphims
$\varphi_{c,m}\colon c \actlE m \Rarr~ c \,{\widetilde{\actle}}\, m$ are also natural
in $c \iN \calc$. To conclude that $(\id, \varphi)$ is a module functor from $\calm$
to $\widetilde{\calm}$, it thus remains to verify the coherence diagrams for
$\varphi$. These readily follow from the uniqueness of the natural  isomorphism
that relates two different right adjoints. 
\end{proof}
  
For the opposite category of a GV-module category over $\calc$ we get analogously 
right $\calc$-ac\-ti\-ons:

\begin{Proposition} \label{prop:ractre,ractle}
Let $\calc$ be a GV-category and $\calm$ a left GV-module category over $\calc$. Then
$\calm\opp$ is a right GV-module category over $(\calc,\tpr)$
and a right module category over $(\calc,\tpl)$,
with left and right exact action bifunctors given by
  \be
  \begin{array}{rl}
  \ractre\, \Colon \calm\opp\Times\calc \!\!\! &\rarr~ \calm\opp \,,
  \Nxl1
  (m,c) \!\!\!& \xmapsto{~~~} G^{-1}c\actlE m
  \eear
  \qquad\text{and}\qquad
  \begin{array}{rl}
  \ractle\, \Colon \calm\opp\Times\calc \!\!\!& \rarr~ \calm\opp \,,
  \Nxl1
  (m,c) \!\!\!& \xmapsto{~~~} G\inv c\actrE m \,,
  \eear
  \label{eq:Mopp-actions}
  \ee
respectively, for $c\iN\calc$ and $m\iN\calm$.
\end{Proposition}

\begin{proof}
First note that we only need to define the $\ractre\,$-action; the $\ractle$-action
is then fixed by the adjunction
$\Hom_{\calm\opp}(m',m \ractlE Gc) \Cong \Hom_{\calm\opp}(m'\ractrE c,m)$
which follows directly from \eqref{eq:HomM.cm.n}.
Further, the proposed bifunctor $\ractre$ is well defined due to $G$ being
contravariant, and it inherits the asserted right and left exactness from $\actle$.
Again we determine the module constraint with the help of the Yoneda lemma. We have
  \be
  \begin{aligned}
  \Hom_\calm(m,n\ractrE(b\ractrE c))
  & \cong\Hom_\calm(m,G^{-1}c\actlE(G^{-1}b\actlE n))
  \Nxl1
  & \cong\Hom_\calm(m,(G^{-1}c\tpL G^{-1}b)\actlE n)
  \Nxl1
  & \cong\Hom_\calm(m,G^{-1}(b\tpR c)\actlE n)
  \cong\Hom_\calm(m,n\ractrE (b\tpR c))
  \end{aligned}
  \label{eq:def:ractre,ractle}
  \ee
for $b,c\iN\calc$ and $m,n\iN \calm$.
The pentagon axiom for this action constraint follows from the injectivity of the
Yoneda embedding and the fact that it holds for the constraints of $\actle$.
\end{proof}

\begin{Remark}
  \Enumerate
  \item
In accordance with our guiding principle, if $(\calc,\tpr,\TI, K)$ is a GV-category,
then $(\calc,\tpl,K,\TI)$ is an `op-GV-category', for which the primary monoidal 
structure is left exact and has $K$ as monoidal unit, and for which the defining 
isomorphisms are $\Hom(\TI,b \tpL c) \Tocong \Hom(Gc,b)$ in place of \eqref{eq:GV1}. 
Similarly there is a notion of op-GV-module category over $(\calc,\tpl)$, for
which the action functors that are analogous to \eqref{eq:action-functors-GVMod} are
required to admit left adjoints. Analogously to Proposition \ref{proposition:actle} 
there is then an associated GV-module category over $(\calc, \tpr)$. 
In this terminology,
$(\calm\opp,\ractle)$ is a right op-GV-module category over $(\calc,\tpl)$.
  \item
Recall from Remark \ref{rem:twistbyG2} that for a module category $(\calm,\actre)$
also $G^{2n}(-) \actrE m$ is a $\calc$-action on $\calm$. 
Since $G^{2}$ is also a monoidal autoequivalence of $(\calc,\tpl)$, an analogous
statement applies to $\actle$: for any integer $n$, also $G^{2n}(-) \actlE m$
defines a $(\calc,\tpl)$-action on $\calm$.
      
\end{enumerate}
\end{Remark}

\begin{Example}
$\calc$ is a $(\calc,\otimes)$-module -- the regular module -- with action given
by the right exact tensor product $\otimes$. The corresponding left exact action 
is given by the $\tpl$-tensor product.
\end{Example}

For $\calm$  a left module category over a GV-category $\calc$ we write the
associated endofunctors of $\calm$ that result from the left action as
  \be
  \begin{array}{rl}
  \L_{c}^{\actre} \Colon & \calm \rarr~ \calm \,, 
  \Nxl1
  & m \xmapsto{~~} c \actrE m 
  \eear
  \qquad\text{and}\qquad
  \begin{array}{rl}
  \L_{c}^{\actle} \Colon & \calm \rarr~ \calm \,, 
  \Nxl1
  & m \xmapsto{~~} c \actlE m \,.
  \eear
  \label{eq:obj-re-LR}
  \ee
The basic adjunction \eqref{eq:HomM.cm.n} for
the left $\calc$-GV-module category $\calm$ then reads
  \be
  \Hom_{\calm}(m, \L_{c}^{\actle}(n)) \cong \Hom_{\calm}(L_{Gc}^{\actre}(m),n) \,.
  \label{eq:actle-actre-adjunction}
  \ee
The counit and unit of this adjunction are
  \be
  \ev_{c,m} \Colon G(c) \actrE (c \actlE m) \rarr~ m \qquad \text{and} \qquad
  \coev_{c,m} \Colon m \rarr~ c \actlE (Gc \actrE m) \,,
  \label{eq:basic-ev-coev}
  \ee
while
  \be
  \begin{tikzcd}[column sep=4.0em]
  c \actlE m \ar{r}[swap]{\coev_{c,c \actle m}} \ar[bend left]{rr}{\id}
  & c \actlE (G(c) \actrE (c \actlE m)) \ar{r}[swap]{c \actle \ev_{c,m}}
  & c \actlE m
  \end{tikzcd}
  \label{eq:triangle-1}
  \ee
and
  \be
  \begin{tikzcd}[column sep=4.5em]
  Gc \actrE m \ar{r}[swap]{Gc \actrE \coev_{c,m}} \ar[bend left]{rr}{\id}
  & Gc \actrE (c \actle (Gc \actrE m)) \ar{r}[swap]{\ev_{c,Gc \actre m}}
  & Gc \actrE m
  \end{tikzcd}
  \label{eq:triangle-2}
  \ee
are the triangle identities of the adjunction.

The $\actle$-action provides us with distinguished natural isomorphisms
$\L_c^{\actle} \cir L_{d}^{\actle} \Tocong L_{c \tpl d}^{\actle}$ for $c,d \iN \calc$.
Likewise, the $\actre$-ac\-ti\-on gives distinguished natural isomorphisms 
$\L_{Gd}^{\actre} \cir L_{Gc}^{\actre} \Tocong L_{G(c \tpl d)}^{\actre}$. Using that
for composable functors $F_1$ and $F_2$ that have right adjoints there is a
commuting diagram $\hspace*{-0.5em}\begin{tikzcd}[column sep=1.4em]
(F_1 \cir F_2) \cir (F_1 \cir F_2)\ra\!\!  \ar[yshift=4pt]{r}{\,}
\ar[yshift=-3pt]{r}{\,} & \!\id \end{tikzcd} \hspace*{-0.4em}$
of counits, we thus have a commuting diagram
  \be
  \begin{tikzcd}[column sep=4.8em,row sep=2.4em]
  G(c \tpL d) \actrE ((c \tpL d) \actlE m) \ar{r}{\ev_{c\tpl d,m}} \ar{d}[swap]{\cong} & m
  \\
  Gd \actrE (Gc \actrE (c \actlE (d\actlE m))) \ar{r}{Gd\actre\ev_{c,d\actle m}}
  & Gd \actrE (d\actlE m) \ar{u}[swap]{\ev_{d,m}}
  \end{tikzcd}
  \label{eq:coher-ev}
  \ee
A similar coherence diagram holds for the collection of morphisms $\coev_{c,m}$.


\subsection{Conjugated pairs of lax and oplax module functors}

The following provides a meaningful notion of GV-module functors:

\begin{Definition} \label{def:laxmodulefunctors}
Let $\calc$ be a GV-category and let $(\calm,\actre)$ and $(\caln,\actrp)$ 
be left GV-module categories over $(\calc,\tpr)$. 
 \\[2pt]
A \emph{lax $\actre$-module functor} from $\calm$ to $\caln$ is a functor
$F\colon \calm\To\caln$ together with a family of morphisms
$f_{c,m}\colon c\actrP F(m) \Rarr~ F(c\actrE m)$ for $c\iN\calc$ and $m\iN\calm$
that obeys coherence conditions in the form of the commutativity of the 
pentagon and triangle diagrams
  \be
  \begin{tikzcd}[column sep=-0.7 em,row sep=2.6 em]
  ~& F((c\tpr d) \actrE m)
  \ar{dl}{\cong}[swap]{F(\alph^{}_{c.d.m})\!\!} \ar[leftarrow]{dr}[xshift=-4pt]{f^{}_{c\tpr d,m}} & ~
  \\
  F(c\actrE(d\actrE m)) \ar[leftarrow]{d}[swap]{f^{}_{c,d\act m}}
  & ~ & (c\tpR d) \,{\act'}\, F(m) \ar{d}{\alph'_{\!c,d,F(m)}}[swap]{\cong}
  \\[0.7em]
  c\actrP F(d\actrE m) \ar[leftarrow]{rr}{\id_c \actrP f^{}_{d,m}} & ~ & c\actrP (d\actrP F(m))
  \end{tikzcd}
  \ee 
and
  \be
  \begin{tikzcd}
  F(\TI\actrE m)
  \ar[leftarrow]{rr}{f^{}_{\TI,m}} \ar{dr}{\!\!\cong}[swap]{F(\lambda^{}_m)}
  & ~& \TI \,{\act'}\, F(m) \ar{dl}[xshift=-4pt]{\lambda'_{F(m)}}[swap]{\cong\!\!}
  \\
  ~ & \raisebox{8pt}{$F(m)$} & ~
  \end{tikzcd}
  \ee
$($with $\alph$, $\lambda$ and $\alph'$, $\lambda'$ the module constraints for the 
actions $\actre$ and $\actrp$, respectively$)$
with respect to the action of $(\calc,\tpr)$.
 \\[2pt]
An \emph{oplax $\actre$-module functor} is a functor $F\colon \calm\Rarr~\caln$,
together with an analogous coherent family of morphisms
$f_{c,m}\colon  F(c\actrE m) \Rarr~ c\actrE F(m)$.
 \\[2pt]
Similarly, a \emph{lax $($oplax$)$ $\actle$-module functor} is a functor
$F\colon \calm\To\caln$, together with an analogous family of morphisms that are
coherent with respect to the action $\actle$ of $(\calc,\tpl)$.
\end{Definition}

Note that no assumption on the existence of adjoints of such functors is made. 
To a monoidal category with exact tensor product one associates a bicategory with
module categories as objects and strong module functors as 1-morphisms. If $\calc$ is
instead a GV-category, we can naturally associate to it
four different bicategories
$\GVModrela(\calc)$, $\GVModreol(\calc)$, $\GVModlela(\calc)$ and $\GVModleol(\calc)$:

\begin{Lemma}
Let $\calc$ be a GV-category.
GV-module categories as objects, lax, respectively oplax, $\actre$-module functors 
as 1-morphisms, and $\actre$-module natural transformations as 2-mor\-phisms form a 
bicategory $\GVModrela(\calc)$, respectively $\GVModreol(\calc)$.
Analogously, GV-mo\-dule categories, lax, respectively oplax, $\actle$-module functors
and $\actle$-module natural transformations form a bicategory $\GVModlela(\calc)$,
respectively ($\GVModleol(\calc)$).
\end{Lemma}

\begin{proof}
The composition of lax $\actre$-module functors is canonically a lax $\actre$-module
functor; analogous arguments apply in the oplax case and for $\actle$-module
functors. The claim thus follows from the fact 
that, for $B\cald$ the delooping of any monoidal category $\cald$, the bicategory
$\Lax(B\cald, \Cat)$ of strong 2-functors, lax natural 2-transformations and
modifications is canonically isomorphic to the bicategory of $\cald$-module
categories, lax module functors and module natural transformations.
The case of oplax module functors is treated analogously.
\end{proof}

According to our guiding principle a functor that has a compatibility with the
$\actre$-action should have a compatibility with the corresponding $\actle$-action
as well. The details depend, however, crucially on whether the considered 
compatibility is strong, lax or oplax. In the situation at hand, 
for a functor $F \colon \calm \Rarr~ \caln$ between left GV-module categories $\calm$
and $\caln$, lax $\actre$-module functor structures and oplax $\actle$-module functor
structures on $F$ are in bijection with each other. This can be seen as follows.
Given a collection of morphisms $f_{c,m}^{\actre}\colon c\actrE F(m) \To F(c\actrE m)$
for $c\iN\calc$ and $m\iN\calm$,
define the \emph{conjugated} collection of morphisms $f^{\actle}_{c, m}\colon
F(c\actlE m) \Rarr~ c \actlE F(m)$ as the family whose members are the composites
  \be
  \begin{aligned}
  F(c \actlE m) \rarr{\coev_{c, F(c\actle m)}}
  c \actlE (Gc \actrE F(c \actlE m)) \rarr{c\actlE f^{\actre}_{Gc,c\actle m}}
  &\; c \actlE F(Gc \actrE (c \actlE m))
  \nxl2
  & \rarr{c\actlE F(\ev_{c,m})} c \actlE F(m) 
  \end{aligned}
  \label{eq:tranport-lax}
  \ee
with the morphisms $\ev$ and $\coev$ as defined in \eqref{eq:basic-ev-coev}.
Note that different choices of adjunction data lead to different conjugated functors;
however, these can be shown to be related by the functors from Lemma 
\ref{lemma:relating-choices}.
Under the adjunction
  \be
  \Hom(F(c \actlE m), c \actlE F(m)) \cong \Hom(Gc \actrE (F (c \actlE m)), F(m))
  \label{eq:transport-via-adj}
  \ee
the composite \eqref{eq:tranport-lax} corresponds to the morphism
  \be
  Gc \actrE (F (c \actlE m)) \rarr{f^{\actre}} F(Gc \actrE (c \actlE m))
  \rarr{F(\ev_{c,m})} F(m) \,.
  \label{eq:transport-under-adj}
  \ee
Conversely, for a collection of morphisms $g_{c,m}^{\actle }\colon  F(c\actlE m) 
\To c\actlE F(m) $, we define the conjugated collection of morphisms
$g^{\actre}_{Gc,m}\colon Gc \actrE F(m) \To F(Gc \actrE m)$ by
  \be
  \begin{aligned}
  Gc \actrE F(m) \rarr{Gc \actre F(\coev_{c,m})} Gc \actrE F(c \actlE (Gc \actrE m))
  & \rarr{g^{\actle}} Gc \actrE (c \actlE F(Gc \actrE m))
  \nxl2
  & \qquad \rarr{\ev_{c,F(Gc\actre m)}} F(Gc \actrE m) \,.
  \end{aligned}
  \label{eq:transport-rex}
  \ee
          
\begin{Proposition} \label{Prop:main-transport-part1}
Let $F\colon \calm \To \caln$ be a functor between left GV-module categories
$\calm$ and $\caln$.
 \Enumerate
  \item
A collection of morphisms $f_{c,m}^{\actre}\colon c\actrE F(m) \Rarr~ F(c\actrE m)$,
for $c\iN\calc$ and $m\iN\calm$,
defines the structure of a lax $\actre$-mo\-dule functor on $F$ if and only if the
conjugated collection $f_{c,m}^{\actle}$ as given in \eqref{eq:tranport-lax}
defines the structure of an oplax $\actle$-module functor on $F$.
 \\
We denote $F$ with the latter oplax structure by $F^{\retole}$.
  \item
A collection of morphisms $g_{c,m}^{\actle }\colon F(c\actlE m) \Rarr~ c\actlE F(m)$
defines the structure of an oplax $\actle$-mo\-dule functor on $F$ if and only if the
conjugated collection $g_{c,m}^{\actre}$ as given in \eqref{eq:transport-rex} defines
the structure of a lax $\actre$-module functor on $F$.
 \\
We denote $F$ with the latter oplax structure by $F^{\letore}$.
  \item
The prescriptions \eqref{eq:tranport-lax} and \eqref{eq:transport-rex} are mutually
inverse, in the sense that
  \be
  ((F,f^{\actre})^{\retole})^{\letore}=(F,f^{\actre}) \qquad \text{and} \qquad
  ((F,g^{\actle})^{\letore})^{\retole}=(F,g^{\actle}) \,.
  \ee
\end{enumerate}
\end{Proposition}

\begin{proof}
      Suppose that $f_{c,m}^{\actre}\colon c\actrE F(m) \To F(c\actrE m)$ is a lax
$\actre$-module functor structure on $F$. Let $f^{\actle}$ be the corresponding
structure under the adjunction given in Equation \eqref{eq:transport-under-adj}.
For objects $x,y \iN \calc$, the morphism
  \be
  G(x \tpL y) \actre F((x \tpL y) \actlE m) \rarr{f^{\actre}}
  F(G(x \tpL y) \actre ((x \tpL y )\actlE m)) \rarr{F(\ev)} F(m)
  \label{eq:direct-xy}
  \ee
corresponds under the adjunction \eqref{eq:transport-via-adj} to the morphism
$F((x \tpL y) \actlE m) \Rarr{f^{\actle}_{x \tpl y, m}} (x\tpL y) \actlE F(m)$,
while the composite
  \be
  \begin{aligned}
  G(x \tpL y) \actre F((x \tpL y) \actle m)
  & \rarr{\cong} G(y) \actre (G(x) \actrE F(x \actlE (y \actlE m)))
  \nxl1
  & \rarr{Gy \actrE f^{\actre}} G(y) \actre F(Gx \actrE (x \actlE (y \actlE m)))
  \nxl1
  & \Rarr{Gy \actrE F(\ev)} G(y) \actrE F(y \actlE m)
  \nxl1
  & \rarr{f^{\actre}} F(G(y) \actrE (y \actlE m)) \,\Rarr{F(\ev)}\, F(m)
  \end{aligned}
  \label{eq:sequence-xy}
  \ee
corresponds under the  adjunction to the composite
  \be
  \begin{aligned}
  F((x \tpL y) \actlE m) \Rarr{\cong} F(x \actlE (y \actlE m))
  \Rarr{f^{\actle}} x \actlE F(y \actlE m) \Rarr{x \actle f^{\actle}}
  & x \actlE (y \actlE F( m))
  \nxl1
  \Rarr{~\cong\,} & (x \tpL y) \actlE F(m) \,.
  \end{aligned}
  \ee
Using that $f^{\actre}$ is coherent and that the diagram \eqref{eq:coher-ev} 
commutes, it follows that the  morphisms \eqref{eq:direct-xy} and
\eqref{eq:sequence-xy} are equal, and thus that $f^{\actle}$ satisfies the pentagon 
diagram that is required for the structure of an oplax $\actle$-module functor.
For the corresponding triangle diagram, note that for $c \eq K$ we can choose as the
adjunction data $\ev_{K,m}$ and $\coev_{K,m}$ in \eqref{eq:basic-ev-coev} a
combination of the unitors.  When doing so, then clearly the triangle identity for 
$f^{\actre}$ implies the triangle identity for $f^{\actle}$. It follows that 
$f^{\actle}$ defines an oplax $\actle$-module functor structure on $F$. Starting with
an oplax $\actle$-module functor structure $g^{\actle}$, one sees analogously that
\eqref{eq:transport-rex} defines a lax $\actre$-module functor structure on $F$. Using
the triangle identities \eqref{eq:triangle-1} and \eqref{eq:triangle-2} we conclude 
that the two constructions are mutually inverse. 
\end{proof}

According to this proposition we can equip every lax $\actre$-module functor with a 
corresponding conjugated oplax $\actle$-module functor structure. This assignment is
functorial:

\begin{Theorem} \label{Thm:bicat-GV-mod}
Let $\calc$ be a GV-category. The assignments in Proposition 
\ref{Prop:main-transport-part1} are functorial in the following sense:
  \Enumerate
  \item \label{item:GV-module-nat}
Let $\eta\colon F \,{\xRightarrow{~\,}}\, H$ be a natural transformation between
functors $F, H\colon \calm \Rarr~ \caln$. Given lax $\actre$-mo\-du\-le functor 
structures on $F$ and $H$, $\eta$ is a $\actre$-module
natural transformation if and only if it is an oplax $\actle$-module
natural transformation for the conjugated  oplax $\actle$-structures on $F$ and $H$.
  \item
Let $K\colon \caln \Rarr~ \calo$ be a functor to another left GV-module category
$\calo$. Given lax $\actre$-module functor structures on $F$ and $K$, the composite 
of the conjugated oplax $\actle$-module functor structures of $F$ and $K$ equals the
conjugated oplax $\actle$-module functor structure of $K \cir F$, i.e.
  \be
  K^{\retole} \circ F^{\retole} = (K \cir F)^{\retole} .
  \ee
\end{enumerate}
The correspondence between lax $\actre$-module functors and oplax $\actle$-mo\-dule
functors extends to an equivalence 
  \be
  \GVModrela(\calc)\rarr\simeq \GVModleol(\calc),
  \ee
of bicategories which is the identity on objects and on 2-morphisms as well as 
the identity on the functors that underly the module functors.
\end{Theorem}

\begin{proof}
Assume that $\eta\colon F \,{\xRightarrow{~\,}}\, H$ is a
$\actre$-module natural transformation between lax $\actre$-module functors
$(F,f^{\actre})$ and $(H,h^{\actre})$. Then the diagram
  \be
  \begin{tikzcd}[row sep=2.3em]
  Gx \actrE F(x \actlE m) \ar{rr}{} \ar{dr}[swap,xshift=12pt]{f^{\actre}_{Gx, x\actle m}}
  \ar{ddd}[swap]{Gx \actrE \eta_{x \actle m}^{}}
  & & F(m) \ar{ddd}{\eta_{m}^{}}
  \\
  & F(Gx \actrE (x \actlE m)) \ar{ur}[swap,xshift=-3pt]{F(\ev_{x,m})}
  \ar{d}{\eta_{Gx \actre (x \actle m)}^{}} &
  \\[7pt]
  & H(Gx \actrE (x \actlE m)) \ar{dr}[yshift=-4pt]{H(\ev_{x,m})} &
  \\
  Gx \actrE H(x \actlE m) \ar{rr}{} \ar{ur}[xshift=5pt]{h^{\actre}_{Gx, x\actle m}}& & H(m)
  \end{tikzcd}
  \label{eq:26}
  \ee
commutes: The triangles define the horizontal arrows; the left quadrilateral commutes
because $\eta$ is a $\actre$-module natural transformation; and the right quadrilateral
commutes by naturality. By \eqref{eq:transport-via-adj}, the resulting commutativity 
of the outer square shows that $\eta$ is a $\actle$-module natural transformation 
for the oplax $\actle$-module functor structures on $F$ and $H$.
 \\[2pt]
To prove the second assertion, consider the oplax $\actle$-module functor structure on
$K^{\retole}{\circ}\,F^{\retole}$. If the composite \eqref{eq:tranport-lax} (for 
$F$ and for $K$, respectively) is inserted for the two arrows in 
$KF(x \actlE m) \Rarr{~~} K(x \actlE F(m)) \Rarr{~~} x \actlE KF(m)$ and the triangle
identity \eqref{eq:triangle-1} is used in the middle, we arrive at the $\actle$-module
functor structure of $(KF)^{\retole}$. Thus the two structures coincide. The equivalence of bicategories follows now
directly from these two statements and Proposition \ref{Prop:main-transport-part1}. 
                    \end{proof}

\begin{Remark}
The two bicategories in Theorem \ref{Thm:bicat-GV-mod}  are distinguished by the
existence of this correspondence: There is no analogous biequivalence between 
$\GVModreol(\calc)$ and $\GVModlela(\calc)$. This can be traced back to the fact 
that, generically, the functors $c \actrE-$ admit only right adjoints.
\end{Remark}
    
We can  give a symmetric characterization of the 1-morphisms in the 
bicategories of Theorem \ref{Thm:bicat-GV-mod}: 

\begin{Definition}
Let $\calm$ and $\caln$ be left GV-module categories over a GV-category $\calc$.
 \Enumerate
 \item 
A \emph{GV-module functor} $F\colon \calm \Rarr~ \caln$ is a functor with a 
conjugated pair of lax $\actre$- and oplax $\actle$-mo\-dule functor structures.
 \item 
For two GV-module functors $F,H\colon \calm \Rarr~ \caln$, a \emph{GV}-module
natural transformation $\eta\colon F \,{\Rightarrow}\, H$ is a natural transformation
that is equivalently a lax $\actre$-module natural transformation or an oplax
$\actle$-module natural transformation.
 \end{enumerate}
\end{Definition}

Clearly there is a bicategory $\GVMod_{\calc}$ of GV-module categories, GV-module
functors and GV-module natural transformations, which is equivalent to the two
bicategories in Theorem \ref{Thm:bicat-GV-mod}. 

Now recall the counit and unit \eqref{eq:basic-ev-coev} of the adjunction
\eqref{eq:actle-actre-adjunction}. For later use we record compatibilities between
the two weak module structures for a GV-functor:

\begin{Lemma} \label{Lemma:Comp-GV-Fun-coev}
Let $F \colon \calm \Rarr~ \caln$ be a GV-functor between left $\calc$-GV-module
categories. The diagrams
  \be
  \begin{tikzcd}[column sep=3.0em]
  F(m) \ar{r}{\coev_{c,F(m)}} \ar{d}[swap]{F(\coev_{c,m})}
  & c \actlE (Gc \actrE F(m)) \ar{d}{c \actlE f^{\actre}_{Gc, m}}
  \\
  F(c \actlE (Gc \actrE m)) \ar{r}{f^{\actle}_{c, Gc \actre m}}
  & c \actlE F(Gc \actrE m)
  \end{tikzcd}
  \label{eq:comp-GV-coev}
  \ee
and
  \be
  \begin{tikzcd}[column sep=3.0em]
  Gc \actrE F(c \actlE m) \ar{r}{Gc \actre f^{\actle}_{c, m}}
  \ar{d}[swap]{f^{\actre}_{Gc, c \actle m}}
  & Gc \actrE (c \actlE F(m)) \ar{d}{\ev_{c,Fm}}
  \\
  F(Gc \actrE (c \actlE m)) \ar{r}{F(\ev_{c,m})} & F(m)
  \end{tikzcd}
  \label{eq:comp-GV-coev-2}
  \ee
commute for all $m \iN \calm$ and $c \iN \calc$.
\end{Lemma}

\begin{proof}
Consider in \eqref{eq:comp-GV-coev} the composite $f^{\actle}_{c, Gc \actre m}
\cir F(\coev_{c,m}) \colon F(m) \Rarr~ c\actlE F(Gc \actrE m)$. Using the definition
\eqref{eq:tranport-lax} of $f^{\actle}$ in terms of $f^{\actre}$ and the naturality
of the coevaluation it follows that this is morphism equal to the composite
  \be
  \begin{aligned}
  F(m) \rarr{\coev_{c,Fm}}\, c \actlE & (Gc \actrE Fm)
  \rarr{c \actlE f^{\actrE}_{Gc, m}} c \actlE F(Gc \actre m)
  \nxl1
  & \rarr{c \actlE F(Gc\actrE\coev_{c,m})} c\actlE F(Gc\actrE (c\actlE (Gc\actrE m)))
  \nxl1
  &\rarr{c \actlE F(\ev_{c,Gc \actre m})} c \actlE F(Gc \actrE m) \,.
  \end{aligned}
  \ee
By the triangle identity \eqref{eq:triangle-2} this is, in turn, equal to the
composite of the
other two arrows in \eqref{eq:comp-GV-coev}. The commutativity of the 
diagram \eqref{eq:comp-GV-coev-2} is seen analogously.
\end{proof}

Next we note that, given a GV-module $\calm$ over a GV-category $\calc$, for every
object $m\iN\calm$ we have a strong $\actre$-module functor
  \be
  \R_{m}^{\actre}\Colon {}_{\calc}^\tpr\calc \,{\equiv}\, {}_{\calc}\calc
  \rarr~ {}_{\calc}\calm \,, \quad c \xmapsto{~~} c \actrE m ~
  \label{eq:R-module-actre}
  \ee
and a strong $\actle$-module functor
  \be
  \R_{m}^{\actle}\Colon {}_{\calc}^\tpl\calc \rarr~ {}_{\calc}\calm \,, \quad
  c \xmapsto{~~} c \actlE m \,. \hspace*{1.7em}
  \label{eq:R-module-actle} 
  \ee

We are now in a position to introduce distributors for general GV-module categories.
Combining the strong module functors $\R_{m}^{\actre}$ and $\R_{m}^{\actle}$
defined in \eqref{eq:R-module-actre} and \eqref{eq:R-module-actle} we get

\begin{Definition} \label{Definition:module-distrib}
Let $\calm$ be a left GV-module category over $\calc$. The oplax $\actle$-module 
structure of the functor $\R_{m}^{\actre}$ defines the \emph{right module distributor}
  \be
  \deltar_{x,y,m} \Colon (x \tpL y) \actrE m \rarr~ x \actlE (y \actrE m) \,,
  \label{eq:right-distri-M}
  \ee
and the lax $\actre$-module structure of $\R_{m}^{\actle}$ defines the
\emph{left module distributor}
  \be
  \deltal_{x,y,m} \Colon x \actrE ( y \actlE m) \rarr~ (x \tpR y) \actlE m \,,
  \label{eq:left-distri-M}
  \ee
for $x,y \iN \calc$ and $m \iN \calm$.
\end{Definition}

Module distributors for a right GV-module category are defined analogously.
We will examine the module distributors in detail in Section \ref{sec:app:dist,dual}.

With the help of the strong module functor $\R_{m}^{\actle}$ we can show that the two 
weak module functor 
structures for a GV-module functor are compatible in the following sense:

\begin{Proposition} \label{Proposition:coher-GV-module}
Let $F\colon \calm \Rarr~ \caln$ be a GV-module functor. The weak module functor
structures $f^{\actre}$ and $f^{\actle}$ obey four pentagon relations:
Commutativity of the ordinary pentagon diagrams
  \be
  \begin{tikzcd}
  x \actre (y \actrE F(m)) \ar[yshift=4pt]{r}{\,} \ar[yshift=-3pt]{r}{\,}
  & F(x \actrE (y \actrE m))
  \end{tikzcd}
  \label{eq:coher-F-lax}
  \ee
and
  \be
  \begin{tikzcd}
  F(x \actlE (y \actlE m)) \ar[yshift=4pt]{r}{\,} \ar[yshift=-3pt]{r}{\,}
  & x \actle (y \actlE F(m))
  \end{tikzcd}
  \label{eq:coher-F-opl}
  \ee
for $x,y\iN\calc$ and $m\iN\calm$, each of which involves just one of the two actions,
and commutativity of
  \be
  \begin{tikzcd}[row sep=2.1em]
  x \actrE F(y \actlE m) \ar{r}{} \ar{d}{} & F(x \actrE (y \actlE m)) \ar{r}{}
  & F((x \tpR y) \actlE m ) \ar{d}{}
  \\
  x \actrE (y \actlE F(m)) \ar{rr} && (x \tpR y) \actlE F(m)
  \end{tikzcd}
  \label{eq:coher-F-lax-opl-1}
  \ee
and
  \be
  \begin{tikzcd}[row sep=2.1em]
  (x \tpL y) \actrE F(m) \ar{r}{} \ar{d}{} & F((x \tpL y) \actrE m) \ar{r}{}
  & F(x \actlE (y \actrE m)) \ar{d}{} \\
  x \actlE (y \actrE F(m)) \ar{rr}{} && x \actlE F(y \actrE m)
  \end{tikzcd}
  \label{eq:eq:coher-F-lax-opl-2}
  \ee
which involve both of them.
\end{Proposition}

\begin{proof}
The two diagrams \eqref{eq:coher-F-lax} and \eqref{eq:coher-F-opl}
are just the coherence diagrams for the individual weak module functors.
For the third diagram consider, for $m \iN \calm$, the composite
${}_{\calc}\calc \Rarr{\R_{m}^{\actle}} {}_{\calc}\calm \Rarr{\,F} {}_{\calc}\caln$
of module functors, which maps $x \iN \calc$ to $F(x \actlE m)$. On the other hand,
the module functor $R^{\actle}_{F(m)}\colon {}_{\calc}\calc \Rarr~ {}_{\calc}\caln$
maps $c$ to $c \actlE F(m)$.
The oplax $\actle$-module structure of $F$ provides a collection of morphisms
$f^{\actle}_{c,m} \colon F\cir R_{m}^{\actle}(c) \eq F(c\actlE m) \Rarr~ c\actlE F(m)
\eq R^{\actle}_{F(m)}(c)$. These morphisms are coherent with respect to the
$\actle$-module structure, which is equivalent to the statement that
$f^{\actle}_{-,m} \colon F \cir \R_{m}^{\actle} \Rarr~ \R_{F(m)}^{\actle}$ is a
$\actle$-module natural transformation.
By Theorem \ref{Thm:bicat-GV-mod}.\ref{item:GV-module-nat},  
$f^{\actle}_{-,m}$
is a also a $\actre$-module natural transformation. The  corresponding coherence
diagram is  precisely the diagram \eqref{eq:coher-F-lax-opl-1}.
Analogously, the lax $\actre$-structure  $f^{\actre}_{x,m}$ of $F$ provides an oplax
$\actle$-module natural transformation $\R^\actre_{F(m)}\Rarr~ F \cir R^\actre_m$,
thus proving the commutativity of the diagram \eqref{eq:eq:coher-F-lax-opl-2}.
\end{proof}

Results analogous to those above hold for $\calc$-right module
categories and functors between them. Indeed, a right module category $\calm$
can be seen as a left $\calctop$-module category, where $\calctop$ is the category
$\calc$ with reversed monoidal structure. 
Consider finally the case of a GV-bimodule category ${}_{\calc}\calm_{\cald}$: The
categories $\calc$ and $\cald$ are GV-categories and $\calm$ is a bimodule category
that is a left $\calc$- and right $\cald$-GV-module category. Equivalently, $\calm$
is a left $\calc{\times}\cald\mopp$-module category, where $\cald\mopp$ 
has the opposite monoidal product.
It is easily seen that $\cald\mopp$ as well as $\calc \Times \cald\mopp$ are
GV-categories and that $\calm$ is a left $\calc{\times}\cald\mopp$-GV-module category.
It therefore follows from Proposition \ref{proposition:actle} that $\calm$ has the
structure of a $(\calc, \tpl)\text-(\cald,\tpl)$-bimodule category.
In particular we thus obtain

\begin{Lemma} \label{Lemma:bimodule-functo-transport}
Let $F \colon \calm \rightarrow \caln$ be a functor between $(\calc,\cald)$-bimodule
categories. For each lax $\actre$-bimodule functor structure on $F$ there is a
conjugated oplax $\actle$-bimodule functor structure, and vice versa, such that
analogous statements as those in Proposition \ref{Prop:main-transport-part1} and 
Theorem \ref{Thm:bicat-GV-mod} are valid.
\end{Lemma}


\subsection{Adjoints of GV-module functors} \label{sec:adjoints}

We now discuss aspects of adjoints of GV-module functors. The structures we find will
be used in the next subsection to obtain snake relations for GV-duality.
In general the adjoint of a GV-module functor is not a GV-module functor. However,
there are interesting transports of (weak) module functor structures to the adjoints.
The first of these, which we call \ftransport\ is as follows:

\begin{Lemma} \label{Lemma:Adj-module-fun}
$($\emph{\Ftransport}.$)$
 \\[2pt]
Let $F \colon \calm\Rarr~\caln$ be a left (op)lax $\actre$-module functor between 
module categories over $(\calc,\tpr)$ which admits a right adjoint 
$F\ra \colon \caln\Rarr~\calm$. Then $F\ra$ is canonically a left (op)lax
$\actle$-module functor with respect to the left exact action $(\calc,\actle)$.
Analogously, if $F$ is an (op)lax $\actle$-module functor admitting a left adjoint
$F\la$, then $F\la$ is an (op)lax $\actre$-module functor.
Moreover, the $\actre$-module functor structure of $F$ is strong if and only if 
the $\actle$-module functor structure of $F\ra$ is strong.
\end{Lemma}                                        

\begin{proof}
Given a functor $F \colon \calm\Rarr~\caln$ with right adjoint $F\ra$, 
for $m,n\iN\calm$ and $c\iN\calc$ we have the composite isomorphisms
  \be
  \begin{aligned}
  \Hom_\caln(c\actrE F(m),n) & \tocong \Hom_\caln(F(m),G^{-1}c\actlE n)
  \tocong \Hom_\calm(m,F\ra(G^{-1}c\actlE n)) 
  \end{aligned}
  \ee
and
  \be
  \begin{aligned}
  \Hom_\caln(F(c\actrE m),n) & \tocong \Hom_\calm(c\actrE m,F\ra(n))
  \tocong \Hom_\calm(m,G^{-1} c\actlE F\ra(n)) \,.
  \end{aligned}
  \ee
Using them as horizontal arrows in the diagram 
  \be
  \begin{tikzcd}
  \Hom_{\caln}(F(c \actrE m),n) \ar{r}{\cong} \ar{d}{}
  & \Hom_{\calm}(m, G^{-1}(c) \actlE F\ra(n)) \ar{d}{}
  \\
  \Hom_{\caln}(c \actrE F(m),n) \ar{r}{\cong} & \Hom_{\calm}(m,F\ra(G^{-1}(c) \actlE n))
  \end{tikzcd}
  \label{eq:comm-tpl-tpr-transport}
  \ee
any one of the vertical arrows defines the other one.
Moreover, the Yoneda lemma provides the corresponding structure on $F$ or on its 
adjoint: In case that $F$ is a lax $\actre$-module functor, the module constraint
$f_{c,m}$ provides the left vertical arrow, so that the right vertical arrow is
a natural morphism $c\actlE F\ra(n) \Rarr~ F\ra(c\actlE n)$. By the injectivity of 
the Yoneda embedding, this satisfies its required coherence conditions because
$f_{c,m}$ does. The oplax case is analogous, with the vertical arrows now pointing
upwards.
The last statement follows from the fact that in the commuting diagram 
\eqref{eq:comm-tpl-tpr-transport} any one vertical arrow is an isomorphism if and 
only if the other one is an isomorphism.
\end{proof}

The second transport mechanism, to be called profunctor-transport, is
obtained with the help of $\calc$-module profunctors; it does not require the
setting of GV-categories. Following \Cite{Def.\,2.1}{shimi20} we give

\begin{Definition}
$($\emph{Profunctor-transport}.$)$
 \\[2pt]
Let $\calm$ and $\caln$ be left module categories over a linear monoidal category
$\calc$. A \emph{$\calc$-module profunctor} from $\calm$ to $\caln$ is a bilinear
functor $H\colon \calm\opp \Times \caln \Rarr~ \vect$ together with a family
  \be
  \theta_{m,n,c}\Colon H(m,n) \rarr~ H(c \actrE m, c \actrE n) 
  \label{eq:29}
  \ee
of morphisms that is natural in $m \iN \calm$ and in $n \iN \caln$, is
dinatural in $c \iN \calc$, and is coherent with respect to 
the monoidal structure and the monoidal unit $\TI$ of $\calc$, i.e.\ satisfies
  \be
  \theta_{m,n,c\otimes d} = \theta_{d\actre m,d\actre n,c} \circ \theta_{m,n,d}
  \qquad \text{and} \qquad \theta_{m,n,\TI} = \id_{H(m,n)}
  \ee
for all $m\iN\calm$, $n\iN\caln$ and $c,d\iN\calc$.
  \\[2pt]
Given two $\calc$-module profunctors  $H_{1}$ and $H_{2}$ from $\calm$ to $\caln$,
a \emph{morphism} $\varphi\colon H_{1} \Rarr~ H_{2}$ \emph{of module profunctors}
is a natural transformation that commutes with the respective $($di$)$natural 
transformations.
\end{Definition}

Our main use of this notion is

\begin{Proposition} \Cite{Lemma\,2.3}{shimi20} \label{proposition:transport}
 \\
Let $F\colon \calm \Rarr~ \caln$ and $K\colon \caln \Rarr~ \calm$ be linear 
functors. There are canonical bijections
 \Enumerate
 \item 
between oplax module functor structures on $F$ and $\calc$-module profunctor structures on
  \be
  \Hom_{\caln}(F(-),-) \Colon \calm\opp \Times \caln \rarr~ \vect \,;
  \label{eq:31}
  \ee
 \item 
between lax module functor structures on $K$ and $\calc$-module profunctor structures on
  \be
  \Hom_{\calm}(-,K(-)) \Colon \calm\opp \Times \caln \rarr~ \vect \,.
  \label{eq:32}
  \ee
\end{enumerate}
\end{Proposition} 

\begin{proof}
The proof is constructive: If $f_{c,m}\colon F(c \actrE m) \Rarr~ c \actrE F(m)$
is the
oplax module structure on $F$, the $\calc$-module profunctor structure on
$\Hom_{\caln}(F(-),-)$ is the composite
  \be
  \theta_{m,n,c}\Colon \Hom(F(m),n) \xrightarrow{\,c\actrE{-}\,}
  \Hom(c\actrE F(m),c \actrE n)
  \xrightarrow{\,f_{c,m}^{*}\,} \Hom(F(c\actrE m),c \actrE n) \,.
  \label{eq:32-2}
  \ee
Conversely, given $\theta$, the oplax module structure is the image of the 
identity under the map
  \be
  \theta_{m,Fm,c}\Colon \Hom(Fm,Fm) \rarr~ \Hom(F(c \actrE m), c \actrE F(m)) \,.
  \label{eq:30}
  \ee
These two constructions are clearly inverse, and it is straightforward to check 
that the respective coherence and naturality diagrams correspond to each other. The 
second bijection follows from the first by considering opposite categories.
\end{proof}

As an immediate consequence we obtain the profunctor-transport of weak module 
functor structures:

\begin{Cor} \label{Corollary:adj-lax}
Let $F\colon \calm \Rarr~ \caln$ be a linear functor with right adjoint 
$F\ra\colon \caln \Rarr~ \calm$.
There is a canonical bijection between oplax $\calc$-module functor structures 
on $F$ and lax $\calc$-module structures on $F\ra$, such that the adjunction
  \be
  \varphi_{m,n} \Colon \Hom_\caln(F(m),n) \rarr~\Hom_\calm(m, F\ra(n))
  \label{eq:C-pro-adju}
  \ee
  is an isomorphism of $\calc$-module profunctors.
\end{Cor}

In this situation, it is natural to expect a compatibility between the unit and
counit of the adjunction and the module action. Note, however, that in general the 
weak module structures of $F$ and $F\ra$ go in opposite directions, so that it 
does not make sense to require the unit and counit to be module natural
transformations. Instead, they have the following coherence properties with respect 
to the module action:

\begin{Lemma} \label{lemma:coherence-counit}
Let $F\colon \calm \Rarr~ \caln$ be an oplax module functor having a right adjoint,
and let $\eta\colon \Id \,{\xRightarrow{~\,}}\, F\ra F$ and 
$\eps\colon FF\ra \,{\xRightarrow{~\,}}\, \Id$ be the 
unit and counit of the adjunction \eqref{eq:C-pro-adju}. Then the diagrams
  \be
  \begin{tikzcd}
  c \actrE m \ar{r}{\eta_{c \actre m}\,} \ar{dd}[swap]{=} & F\ra F(c \actrE m)
  \ar{d} {} 
  && c \actrE FF\ra (n) \ar[<-]{d}{} \ar{r}{c \actrE \eps_n\,} & c \actrE n \ar{dd}{=}
  \\
  ~ & F\ra (c \actrE F(m)) & \text{and}
  & F(c \actrE F\ra (n)) & ~
  \\
  c \actrE m \ar{r}{c \actrE \eta_m\,} & c \actrE F\ra F(m) 
  \ar{u}{} 
  && F(F\ra (c \actrE n)) \ar{r}{\eps_{c \actre n}\,} \ar[<-]{u}{} & c \actrE n
  \end{tikzcd}
  \label{eq:coher-eta-eps}
  \ee
commute for all $m \iN \calm$, $n \iN \caln$ and $c \iN \calc$, using the 
oplax and lax module functor structures of $F$ and $F\ra$, respectively.
\end{Lemma}

\begin{proof}
Denote the lax module functor structure of $F$ by $f$ and the oplax structure of 
$F\ra$ by $g$.  
Since $\varphi$ is an isomorphism of  $\calc$-module profunctors, the diagram
  \be
  \begin{tikzcd}
  \Hom(Fm,n) \ar{r}{c \actrE} \ar{d}{\varphi_{m,n}}
  & \Hom(c \actrE F(m), c \actrE n) \ar{r}{f^{*}}
  & \Hom(F(c \actrE m), c \actrE n) \ar{d}{\varphi_{c \actre m, c \actre n}}
  \\
  \Hom(m, F\ra(n) \ar{r}{c\actrE}) & \Hom(c\actrE m, c\actrE F\ra (n)) \ar{r}{g_{*}}
  & \Hom(c \actrE m, F\ra (c\actrE n))
  \end{tikzcd}
  \label{eq:34}
  \ee
commutes for all $m,n$. Evaluating \eqref{eq:34} at $n\eq F(m)$ gives
the commutative diagram
  \be
  \hspace*{-0.8em}
  \begin{tikzcd}[column sep=1.9em]
  \Hom(Fm,F(m)) \ar{r}{c \actrE} \ar{d}{\varphi_{m,F(m)}}
  & \Hom(c \actrE F(m), c \actrE F(m)) \ar{r}{f^{*}}
  & \Hom(F(c \actrE m), c \actrE F(m)) \ar{d}{\varphi_{c \actre m, c \actre F(m)}}
  \\
  \Hom(m,F\ra(F(m)) \ar{r}{c\actrE}) & \Hom(c\actrE m, c\actrE F\ra(F(m))) \ar{r}{g_*}
  & \Hom(c \actrE m, F\ra (c\actrE F(m))) 
  \end{tikzcd}
  \label{eq:34-2}
  \ee
Taking the identity morphism in the upper left corner, this yields the left 
diagram in \eqref{eq:coher-eta-eps}. The other diagram in \eqref{eq:coher-eta-eps}
follows by duality. 
\end{proof}

Profunctor transport is compatible with composition:

\begin{Lemma} \label{lemma:prof-comp}
Let $F_{1}\colon \calm \Rarr~ \caln$ and  $F_{2}\colon \caln \Rarr~ \calo$ be oplax 
$\calc$-module functors between left $\calc$-mo\-du\-le categories. The canonical 
natural isomorphism $(F_{2} F_{1}) \ra \Cong F_{1}\ra F_{2}\ra$ is an isomorphism 
of lax module functors.
\end{Lemma}

\begin{proof}
The claim follows from the definition of the transport: all morphisms in the diagram
  \be
  \begin{tikzcd}[row sep=2.4em]
  \Hom_{\calo}(F_{2}F_{1}m,o) \ar{r}{\cong} \ar{d}[swap]{\cong}
  & \Hom_{\calm}(m,(F_{2}F_{1})\ra o) \ar{d}{\varphi}
  \\
  \Hom_{\caln}(F_{1}m,F_{2}\ra o) \ar{r}{\cong} & \Hom_{\calm}(m,F_{1}\ra F_{2}\ra o),
  \end{tikzcd}
  \ee
are isomorphisms, and all of them except for possibly the one labeled $\varphi$ are
isomorphisms of $\calc$-module profunctors. But since the diagram commutes, $\varphi$
is in fact an isomorphism of $\calc$-module profunctors, too. By the Yoneda Lemma,
$\varphi$ provides the canonical isomorphism $(F_{2} F_{1})\ra 
\Cong F_{1}\ra F_{2}\ra$, which is thus an isomorphism of lax module functors.
\end{proof}

We can summarize the constructions given above as follows. Let $\calc$ be a
GV-category, let $\calm$ and $\caln$ be left GV-module categories over $\calc$,
and let $F \colon \calm \Rarr~ \caln$ be a functor. The 
transport constructions are the rows in the following table:
  \be
  \begin{tabular}{|lc|c|c|}
    \hline &&& \\[-0.9em] 
  (T1) & profunctor transport & $F$ oplax-$\actre$ &$F\ra$ lax-$\actre$ 
    \\ &&& \\[-0.9em] \hline &&& \\[-0.9em]
  (T2) & profunctor transport & $F$ oplax-$\actle$ &$F\ra$ lax-$\actle$ 
    \\ &&& \\[-0.9em] \hline \hline && \\[-0.9em]
  (T3) & \ftransport &  $F$ oplax-$\actre$ &$F\ra$ oplax-$\actle$
    \\ &&& \\[-0.9em] \hline &&& \\[-0.9em]
  (T4) & \ftransport &  $F$ lax-$\actre$ &$F\ra$ lax-$\actle$
    \\ &&& \\[-0.9em] \hline &&& \\[-0.9em]
  (T5) & \ftransport &  $F$ strong-$\actre$ &$F\ra$ strong-$\actle$ 
    \\[-0.9em] &&& \\ \hline 
\end{tabular}
\label{eq:table}
\ee
Each of these transport constructions can be applied from left to right as well as
from right to left, i.e.\ if $F$ has a left adjoint, then $F\la$ acquires  
transported  structure from $F$.
In parallel with Proposition \ref{Prop:main-transport-part1}, for a functor
$F \colon \calm \Rarr~ \caln$ this provides a pairing between lax $\actre$-
and oplax $\actle$-module structures:

\begin{Cor} \label{cor:transport}
Let $\calm$ and $\caln$ be left-module categories over a GV-category $\calc$, and let
$F \colon \calm \Rarr~ \caln$ be a functor admitting a right adjoint. Then the 
transports via the right adjoint provide a bijection between lax $\actre$-module 
functor structures and oplax $\actle$-module functor structures on $F$.
In case that $F$ admits a left adjoint, the transports via the left adjoint 
provide a (potentially different) bijection between lax $\actre$-module functor
structures and oplax $\actle$-module functor structures on $F$.
\end{Cor}

\begin{proof}
Suppose that $F$ is a lax $\actre$-module functor with a right adjoint $F\ra$. Then
by the transport (T4), $F\ra$ is a lax $\actle$-module functor, while by (T2) $F$ 
acquires an oplax $\actle$-module functor structure. Since both steps are bijections
of structures, the structures are in bijection.
 \\
If $F$ has a left adjoint, then the transports (T1) and (T3) provide
again a bijection between the same structures: If $F$ is a lax $\actre$-module
functor, then by (T1) $F\la$ is a oplax $\actre$-module functor, and by (T3)
$F$ is an oplax $\actle$-module functor. 
\end{proof}

In case that a lax $\actre$-module functor $F$ admits both a right and a left adjoint,
the transported structures coincide:

\begin{Proposition} \label{proposition:transports-agree}
Let $F\colon \calm \Rarr~ \caln$ be a lax $\actre$-module functor admitting a right
adjoint. Then the oplax $\actle$-module functor structure on $F$ given in Proposition
\ref{Prop:main-transport-part1} coincides with the transported oplax $\actle$-module
functor structure given in Corollary \ref{cor:transport}. Analogously, if $F$ admits
a left adjoint, then both oplax $\actle$-module functor structures coincide.
In particular, if $F$ has both a right and a left adjoint, then both transported 
module structures from Corollary \ref{cor:transport} are the same. 
\end{Proposition}

\begin{proof}
Assume that $F$ is a lax $\actre$-module functor and has a right adjoint $F\ra$. Under
the adjunction $\Hom_{\calm}(F(c \actlE m), c \actlE F(m)) \Cong 
\Hom(Gc \actrE F(c \actlE m), F(m))$, the conjugated oplax $\actle$-mo\-du\-le functor
structure $f^{\actle}$ of $F$ is the composite
  \be
  Gc \actrE F(c \actlE m) \rarr{f^{\actre}_{Gc, c \actlE m}}
  F(Gc\actrE(c{\actle}m)) \rarr{F(\ev_{c,m})} F(m) \,. 
  \label{eq:conj-adju}
  \ee
On the other hand, the transported module functor structure $\widetilde{f}^{\actle}$
is the image of the identity under the $(\calc,\tpl)$-profunctor structure
$\Hom_{\calm}(Fm,Fm) \Rarr{\theta_{m,Fm,c}} \Hom(F(c\actlE m),c \actlE F(m))$
that is obtained by the diagram
  \be
  \begin{tikzcd}[column sep=1.9em]
  \Hom(Fm,n) \ar{d}{} \ar{rr}{\theta_{m,n,c}} & & \Hom(F(c \actlE m), c \actlE n)
  \\
  \Hom(m,F\ra(n)) \ar{d}{c\actlE{-}} & & \Hom(Gc \actrE F(c \actlE m), n) \ar{u}{}
  \\
  \Hom(c \actlE m, c \actlE F\ra(n)) \ar{r}{} & \Hom(Gc\actrE (c\actlE m), F\ra(n)) \ar{r}{}
  & \Hom(F (Gc \actrE (c \actlE m)), n) \ar{u}{f^{\actrE}}
  \end{tikzcd}
  \ee
Taking  $n \eq F(m)$ and the identity morphism in the upper left corner yields 
the mappings
  \be
  \begin{tikzcd}
  \id_{Fm} \ar[mapsto]{d} & & \widetilde{f}^{\actlE}_{c,m} 
  \\
  \eps_{F}(m) \ar[mapsto]{d} & & \ev_{c,m} \cir f^{\actrE}_{Gc,c \actlE m}
  \ar[mapsto]{u}
  \\
  c \actlE \eps_{F}(m) \ar[mapsto]{r} & \eps_{F}(m) \cir \ev_{c,m} \ar[mapsto]{r}
  & F(\ev_{c,m}) \ar[mapsto]{u}
  \end{tikzcd}
  \ee
($\eps_{F}\colon \id \To F\ra F$ is the unit of the adjunction). Thus in view of
Equation \eqref{eq:conj-adju} the two oplax $\actle$-module functor structures
coincide. The remaining statements follow by duality. 
\end{proof}


\subsection{Applications: Distributors and duality} \label{sec:app:dist,dual}

Recall the right and left module distributors from Definition
\ref{Definition:module-distrib}.
The fact that the module distributors come from weak module structures implies the
commutativity of two pentagon diagrams for the right module distributor and of 
two pentagons for the left module distributor, while the two
coherence diagrams involving both the right and the left module distributors
follow from the two mixed pentagons in Proposition \ref{Proposition:coher-GV-module}.
Altogether we have the following coherences for the module distributors:

\begin{Proposition} \label{Proposition:module-distri-coherence}
Let $\calm$ be a left GV-module category over $\calc$. The following six 
pentagon diagrams commute for all $x,y,z\iN\calc$ and $m\iN\calm$:
 \begin{enumerate}
 \item
The pentagons
  \be
  \begin{tikzcd}
  ((x \tpL y) \tpL z ) \actre m \ar[yshift=4pt]{r}{\,} \ar[yshift=-4pt]{r}{\,}
  & x \actlE (y \actlE (z \actrE m))
  \end{tikzcd}
  \label{eq:right-dist-coh-1}
  \ee
and
  \be
  \begin{tikzcd}
  (x \tpL y ) \actre (z \actrE m) \ar[yshift=4pt]{r}{\,} \ar[yshift=-4pt]{r}{\,}
  & x \actlE (y \actrE (z \actrE m))
  \end{tikzcd}
  \label{eq:right-distr-coh-2}
  \ee
for the distributor $\deltar$;
 \item
the pentagons
  \be
  \begin{tikzcd}
  (x \tpR y ) \tpr (z \actlE m) \ar[yshift=4pt]{r}{\,} \ar[yshift=-4pt]{r}{\,}
  & (x \tpR (y \tpR z)) \actlE m
  \end{tikzcd}
  \label{eq:left-distr-coh-1}
  \ee
and
  \be
  \begin{tikzcd}
  x \actre ((y \tpL z) \actlE m) \ar[yshift=4pt]{r}{\,} \ar[yshift=-4pt]{r}{\,}
  & (x \tpR y) \actle (z \actlE m)
  \end{tikzcd}
  \label{eq:left-distr-coh-2}
  \ee
for the distributor $\deltal$;
 \item
the two mixed pentagons
  \be
  \begin{tikzcd}
  (x \actrE (y \tpL z)) \actre m \ar[yshift=4pt]{r}{\,} \ar[yshift=-4pt]{r}{\,}
  & (x \tpR y) \actle (z \actrE m)
  \end{tikzcd}
  \label{eq:mixed-1}
  \ee
and
  \be
  \begin{tikzcd}
  (x \tpL y) \actre (z \actlE m) \ar[yshift=4pt]{r}{\,} \ar[yshift=-4pt]{r}{\,}
  & x \tpl ((y \tpR z) \actlE m) \,.
  \end{tikzcd}
  \label{eq:mixed-2}
  \ee
 \end{enumerate}
\end{Proposition}

\begin{proof}
The commutativity of the diagram \eqref{eq:right-dist-coh-1} expresses the fact
that $\deltar$ is an oplax $\actre$-mo\-du\-le structure on the functor 
$\R_{m}^{\actre}$; commutativity of \eqref{eq:right-distr-coh-2}
follows from \eqref{eq:eq:coher-F-lax-opl-2} for $F \eq \R_{m}^{\actre}$.
The commutativity of \eqref{eq:left-distr-coh-1} expresses the fact that 
$\deltal$ is a lax $\actre$-module structure on $\R^{\actle}_{m}$, while
commutativity of \eqref{eq:left-distr-coh-2} follows from 
\eqref{eq:coher-F-lax-opl-1} for $F \eq \R^{\actle}_{m}$.
Finally, the diagram \eqref{eq:mixed-1} is \eqref{eq:coher-F-lax-opl-1} for 
$F \eq \R^{\actre}_{m}$, while the diagram \eqref{eq:mixed-2} is 
\eqref{eq:eq:coher-F-lax-opl-2} for $F \eq \R^{\actle}_{m}$.
\end{proof}

We now specialize to the case of the regular $\calc$-module category
${}_{\calc}\calc$. From the module distributors we then obtain natural isomorphisms
  \be
  \deltar_{x,y,z} \Colon (x \tpL y) \tpR z \rarr~ x \tpL (y \tpR z)
  \qquad \text{and} \qquad
  \deltal_{x,y,z} \Colon x \tpR ( y \tpL z) \rarr~ (x \tpR y) \tpL z \,.
  \label{eq:deltar,deltal}
  \ee
By Proposition \ref{Proposition:module-distri-coherence} these morphisms obey all
coherence conditions that the left and right distributors in a linearly 
distributive category (see e.g.\ \cite{coSe3,past}) have to satisfy.
Thus our result is stronger than Proposition 4.11 of \cite{fssw}, in which 
it was shown that all these conditions except for the two pentagon equations 
that involve both left and right distributors can be derived from the definition 
of distributors in terms of weak module functors.

We can equally well treat $\calc$ as a right GV-module category $\calc_{\calc}$.
The lax $\actre$-module functor structure of $\L_{y}^{\ractle} \colon \calc_{\calc}
\Rarr~ \calc_{\calc}$ consists again of morphisms
$\widetilde{\delta}^{\ihr}_{x,y,z} \colon (x \tpL y) \tpR z \Rarr~ x \tpL (y \tpR z)$,
and analogously we obtain morphisms 
$\widetilde{\delta}^{\ihl}_{x,y,z} \colon x \tpR (y \tpL z) \Rarr~ (x \tpR y) \tpL z$.
Again it follows that these satisfy all coherence diagrams for distributors. However,
they do not constitute new structures:

\begin{Lemma}
For all objects $x,y,z$ in a GV-category $\calc$ we have
$\widetilde{\delta}^{\ihr}_{x,y,z} \eq \deltar_{x,y,z}$ and
$\widetilde{\delta}^{\ihl}_{x,y,z} \eq \deltal_{x,y,z}$.
\end{Lemma}

\begin{proof}
By construction, $\deltar_{x,y,z}$ is the image of the morphism 
$\ev_{x,y}\tpR \id_{z}$ under the adjunction
  \be
  \Hom(G(x) \tpR ((x\tpL y) \tpr z), y \tpR z)
  \cong \Hom((x \tpL y) \tpR z, x \tpl (y \tpR z)) \,.
  \ee
while under the adjunction
  \be
  \Hom((G(x) \tpR (x\tpL y)) \tpR z, y \tpR z)
  \cong \Hom(G(x) \tpR (x \tpL y),  (y \tpR z) \tpL G(z))
  \ee
it is mapped to $\coevp_{y,Gz}\cir\ev_{x,y}$. On the other hand,
$\widetilde{\delta}^{\ihr}_{x,y,z}$ is the image of $\id_{x} \tpL \coevp_{y,Gz}$ under
  \be
  \Hom(x \tpL y, (x \tpL (y \tpR z)) \tpL G(z))
  \cong \Hom((x \tpL y) \tpR z, x \tpL (y \tpR z))) \,,
  \ee
and under the adjunction
  \be
  \Hom(x \tpL y, x \tpL ((y \tpR z) \tpl G(z)))
  \cong \Hom(G(x) \tpR (x \tpL y), (y \tpR z) \tpL G(z))
  \ee
it gets mapped to the morphism $ \coevp_{y,Gz}\circ\ev_{x,y}$ as well. It thus
follows that $\widetilde{\delta}^{\ihr}_{x,y,z} \eq\deltar_{x,y,z}$.
The second statement is shown analogously.
\end{proof}

The distributors also possess another symmetry: applying the anti-equivalence $G$
to $\deltar_{x,y,z}$, we obtain a morphism
  \be
  (Gz \tpL Gy) \tpR Gx \cong G(y{\tpr}z) \tpR Gx \cong G(x \tpL (y{\tpr}z))
  \rarr~ G((x{\tpl}y) \tpR z) \cong Gz \tpL (Gy \tpR Gx) \,.
  \ee

\begin{Proposition} \label{prop:G.dist}
The distributors of a GV-category $\calc$ satisfy
  \be
  G(\deltar_{x,y,z}) = \deltar_{Gz,Gy,Gx} \qquad \text{and} \qquad
  G(\deltal_{x,y,z}) = \deltal_{Gz,Gy,Gx}
  \label{eq:Gdelta=deltaGGG}
  \ee
for all $x,y,z \iN \calc$.
\end{Proposition}

\begin{proof}
We equip $G$ with the structure of a module functor,
so that we can use the general results on lax/oplax module functor structures.
The category $\calc\opp$ is a left $\calc$-module category with action
$x \actrE \overline{y} \eq \overline{y \tpL G(x)}$ for $x \iN \calc$ and
$\overline{y} \iN \calc\opp$, and similarly it is a right module category with 
action $\overline{y} \ractrE x \eq \overline{Gx \tpL y}$. By a direct computation 
one sees that the ${}_{\calc}\calc$-module endofunctor 
$G\cir \R^{\ractre}_{G^{-1}(x)}\cir G^{-1}$, with $\R^{\ractre}_{G^{-1}(x)}$ the 
left module endofunctor that maps $\overline{y}$ to $\overline{x\tpL y} $, is 
isomorphic to $R^{\actre}_{G(x)}$ as a $\actre$-module functor. These two functors 
are then also equivalent as oplax $\actle$-module functors. This translates 
directly into the equality $G(\deltar_{x,y,z}) \eq \deltar_{Gz,Gy,Gx}$. 
The second equality in \eqref{eq:Gdelta=deltaGGG} follows analogously.
\end{proof}

Module functors preserve the distributors in the following sense:

\begin{Proposition} \label{proposition:module-fun-distr}
Let $F \colon \calm \Rarr~ \caln$ be a strong $\actre$-module functor. Then $F$
preserves the distributor $\deltar$ in the sense that the pentagon
  \be
  \begin{tikzcd}[row sep=2.5em,column sep=3.8em]
  F((d \tpL c) \actrE m)  \ar{r}{F(\deltar_{d,c,m})\,}
  \ar{d}[swap,xshift=1pt]{f^{\actre}_{d\tpl c,m}}
  & F(d \actlE (c \actrE m)) \ar{r}{f^{\actle}_{d,c\actre m}}
  & d \actlE F(c \actrE m) \ar{d}[xshift=-1pt]{d\actle f^{\actre}_{c,m}}
  \\
  (d \tpL c) \actrE F(m) \ar{rr}{\deltar_{d,c,F(m)}} && d \actlE (c \actrE F(m))
  \label{eq:penta-pres-dist}
  \end{tikzcd}
  \ee
commutes for all $c, d \iN \calc$ and all $m \iN \calm$.  A strong $\actle$-module
functor preserves the distributor $\deltal$ in an analogous manner. 
\end{Proposition}

\begin{proof}
With the strong $\actre$-module functor structure on
$R^{\actre}_{m} \colon {}_{\calc}^\tpl\calc \Rarr~ {}_{\calc}\calm$, the module
functor constraint of $F$ provides an isomorphism
$R_{Fm}^{\actre} \Cong F \cir R_{m}^{\actre}$ of module functors from $\calc$ to
$\caln$. By Theorem \ref{Thm:bicat-GV-mod} this isomorphism is also an isomorphism of
the conjugated oplax $\actle$-structures on these functors. This is precisely the
commutativity of the diagram \eqref{eq:penta-pres-dist}. 
\end{proof}
  
Let now ${}_{\calc}\calm_{\cald}$ be a GV-bimodule category over GV-categories $\calc$
and $\cald$.
Since for every $c \iN \calc$ the functor $\L_{c}^{\actre} \colon \calm \Rarr~ \calm$ 
with $\L_{c}^{\actre}(m) \eq c \actrE m$ is a $\cald$-$\actrer$-module functor, we
obtain from its oplax $\actler$-module structure the distributor
  \be
  \deltal_{c,m,d} \Colon  c \actrE (m \actleR d)\rarr~ (c \actrE m) \actleR d \,.
  \ee
Analogously we obtain the natural morphisms
  \be
  \deltar_{c,m,d} \Colon (c \actlE m) \actreR d \rarr~ c \actlE (m \actreR d) \,.
  \ee
There are in total $8 \,{\cdot}\, 4 \eq 32$ pentagon diagrams for the distributors
of a GV-bimodule category $\calm$; they can be schematically symbolized as
$\calm \Box \Box \Box$, $\Box \calm \Box \Box$, $\Box \Box \calm \Box$ and
$\Box \Box \Box \calm$, with $\Box \iN \{\tpr, \tpl\}$.
For example, the pentagon of type $\tpr \calm \tpl \tpr$ is assembled from the 
two possible composite morphisms $(c \actrE (m\actleR d)) \actreR d'
\,{\xrightrightarrows{~\phantom~}}\, (c\actrE m) \actleR (d\tpR d')$
for all $c\iN\calc$, $d,d'\iN\cald$ and $m\iN\calm$.

\begin{Cor}
Let $\calm$ be a GV-bimodule category. Each of the $32$ pentagon diagrams for the
distributors of $\calm$ commutes. 
\end{Cor}

\begin{proof}
The commutativity of the $16$ pentagons of type $\calm \Box \Box \Box$ or
$\Box \Box \Box \calm$ follows by regarding $\calm$ just as a left or right module 
category. The commutativity of all $16$ pentagons of type $\Box \calm \Box \Box $ or
$\Box \Box \calm \Box$ follows from  Proposition \ref{Proposition:coher-GV-module}. 
The pentagon of type $\tpr \calm \tpl \tpr$, for example, is obtained by applying
Proposition \ref{Proposition:coher-GV-module} to the $\cald$-$\actrer$-module functor
$\L_{c}^{\actre}$.
\end{proof}

Next we introduce GV-analogues of the evaluation and coevaluation morphisms of a
rigid duality. To this end we make use of the distinguished isomorphisms
  \be
  \Hom(Gy\tpR y,K) \rarr\cong \Hom(Gy,Gy) 
  \label{eq:disti-iso-1}
  \ee
and
  \be
  \Hom_\calc(y \tpR G^{-1}y,K) \cong \Hom_\calc(y,GG^{-1}y) =\Hom (y,y) \,,
  \label{eq:disti-iso-2}
  \ee  
which are special cases of the defining isomorphisms \eqref{eq:GV1} of the 
GV-structure.

\begin{Definition} 
The right and left \emph{GV-evaluation morphisms}
  \be
  \tevr_y \Colon Gy\tpR y \rarr~ K \qquad \text{and} \qquad 
  \tevl_y \Colon y\tpR G^{-1}y \rarr~ K
  \label{eq:def:tevr,tevl}
  \ee
are the pre-image of the identity morphisms $\id_{Gy}$ and $\id_y$ under the
isomorphisms \eqref{eq:disti-iso-1} and \eqref{eq:disti-iso-2}, respectively,
Analogously, the right and left \emph{GV-coevaluation morphisms}
  \be
  \coevr_x := G(\evl_x) \Colon \TI \rarr~ x\tpL G(x) 
  \label{eq:def:coevr}
  \ee
and
  \be        
  \coevl_x := G^{-1}(\evr_x) \Colon \TI \rarr~ G^{-1}(x)\tpL x \,.
  \label{eq:def:coevl}
  \ee
are obtained from the isomorphisms
  \be
  \Hom_\calc(1,x\tpL Gy) \cong \Hom_\calc(GK, G(y\tpR G^{-1} x))
  \cong \Hom_\calc(y\tpR G^{-1}x, K) \,.
  \ee
\end{Definition}

We now show that when complemented with the distributors \eqref{eq:deltar,deltal},
the GV-evaluation and GV-coevaluation morphisms $\tevl,\, \tevr$ and 
$\coevl,\, \coevr$ obey snake relations which involve distributors and
generalize the familiar duality structure of a rigid monoidal category.
For the monoidal case this was already observed in Theorem 4.5 of (the revised
version of) \cite{coSe3}. The proof given there involves ``a straightforward
verification to check that *-autonomous categories are weakly distributive, though
the diagrams can be pretty horrid.'' Our approach generalizes the statement to the
case of module categories and provides a conceptual understanding of the diagrams.
We work in the setting of a left $\calc$-module category $\calm$. Recall from
Equation \eqref{eq:basic-ev-coev} the collection of morphisms
$\coev_{c,m} \colon m \Rarr~ c \actlE (Gc \actrE m)$. It follows from our convention
concerning the adjunction \eqref{eq:definingadjunction} for the case $\calm \eq \calc$
that in this case the composite
$\TI \Rarr{\coev_{c,\TI}} c \tpL (Gc \tpR \TI) \Cong c \tpL Gc$ coincides with
the GV-coevaluation $\coevr_{c} \colon \TI \Rarr~ c \tpL Gc$,
and analogously for $\ev_{c,\TI}$:
  \be
  \coev_{\!c,\TI}^{} = \coevr_{\!c} \qquad\text{and}\qquad
  \ev_{\!c,K}^{} = \evr_{\!c} \,.
  \ee
The GV-(co)evaluations and the module distributors determine all (co)evaluations:

\begin{Lemma}\label{Lemma:GV-coev-gives-all}
Let $\calm$ be a GV-module category. For all $m \iN \calm$ and $c \iN \calc$
the diagrams
  \be
  \begin{tikzcd}[column sep=3.5em]
  m \ar{r}{\coev_{c,m}} \ar{d}[swap]{\cong} & c \actlE (Gc \actrE m)
  \\
  \TI \actrE m \ar{r}{\coev_{c} \actre m}
  & (c \tpL Gc) \actrE m \ar{u}[swap]{\deltar_{c,Gc,m}}
  \end{tikzcd}
  \qquad \text{and} \qquad
  \begin{tikzcd}[column sep=3.2em]
  Gc \actrE (c \actlE m) \ar{r}{\ev_{c,m}} \ar{d}[swap]{\deltal_{Gc,c,m}} & m
  \\
  (Gc \tpR c) \actlE m \ar{r}{\ev_{c} \actlE m} & K \actlE m \ar{u}[swap]{\cong}
  \end{tikzcd}
  \label{eq:ev-from-GVev}
  \ee
  commute.
\end{Lemma}

\begin{proof}
Applying Lemma \ref{Lemma:Comp-GV-Fun-coev} to the GV-module functor 
$R_{m} \colon \calc \Rarr~ \calm$ and the object $\TI \iN \calc$ directly shows
the commutativity of the diagram
  \be
  \begin{tikzcd}[column sep=3.5em, row sep=2.3em]
  R_{m}(\TI) \eq \TI \actrE m \ar{r}{\coev_{c,\TI \actre m}}
  \ar{d}[swap]{R_{m}(\coev_{c,\TI})}
  & c \actlE (Gc \actrE (\TI \actrE m)) \ar{d}{\cong} \ar{d}{\cong}
  \\
  (c \tpL (Gc \tpR \TI)) \actrE m \ar{r}{\deltar_{c,Gc\tpr \TI,m}}
  & c \actlE ((Gc \tpR 1) \actrE m)
  \end{tikzcd}
  \ee
Inserting the unitors, this implies the commutativity of the first of the
diagrams \eqref{eq:ev-from-GVev}.
The commutativity of the second diagram is shown analogously.
\end{proof}

It is worth pointing out that in view of the fact that the distributors are
generically not isomorphisms, it is not obvious that the GV-(co)evalu\-a\-tions obey
appropriate snake relations. However, with the help of Lemma
\ref{Lemma:GV-coev-gives-all} we can show that this is indeed the case:

\begin{Proposition}\label{prop:GVsnakes}
\emph{$[$Snake relations.$]$}
 \\
Let $\calc$ be a GV-category. For every $c \iN \calc$ the diagrams
  \be
  \begin{tikzcd}[column sep=3.9em]
  c \ar{r}[swap]{\coevr_c \tpR \id} \ar[bend left]{rrr}{\id_c}
  & (c \tpL Gc) \tpr c \ar{r}[swap]{\deltar_{c,Gc,c}}
  & c \tpL (Gc \tpr c) \ar{r}[swap]{\id \tpL \evr_{c}} & c
  \end{tikzcd}
  \label{eq:GV-triangle-1}
  \ee
and
  \be
  \begin{tikzcd}[column sep=3.9em]
  Gc \ar{r}[swap]{\id \tpR \coevr_c} \ar[bend left]{rrr}{\id_c}
  & Gc \tpR (c \tpL Gc) \ar{r}[swap]{\deltal_{Gc,c,Gc}}
  &(Gc \tpR c) \tpL Gc \ar{r}[swap]{ \evr_{c} \tpL \id} & Gc
  \end{tikzcd}
  \label{eq:GV-triangle-2}
  \ee
commute.
\end{Proposition}

\begin{proof}
By Lemma \ref{Lemma:GV-coev-gives-all} the diagram
  \be
  \begin{tikzcd}[column sep=3.4em, row sep=2.3em]
  c \tpL K \cong c \ar{r}{\coev_{c,c \tpL K}}
  \ar[xshift=16pt]{d}[swap]{\coevr_{c} \tpR 1}
  & c \tpL (Gc \actrE (c \tpL K)) \ar[xshift=-6pt]{d}{\cong} \ar{r}{c \tpL \ev_{c,K}}
  & c \tpL K \cong c \ar[xshift=-6pt]{d}{=}
  \\
  (c \tpL Gc) \tpR c \ar{r}{\deltar} & c \tpL (Gc \tpR c) \ar{r}{c \tpL \evr_{c}}
  & c \tpL K
  \end{tikzcd}
  \ee
commutes. By Equation \eqref{eq:triangle-1} the horizontal arrows in the first row
compose to $\id_c$. This yields commutativity of \eqref{eq:GV-triangle-1}.
Commutativity of \eqref{eq:GV-triangle-2} is shown analogously.
\end{proof}

\begin{Remark}
In particular, the evaluation and coevaluation morphisms endow every object of a
GV-category with the structure of a (left and right) \emph{nuclear object}
in the sens of Definition A.1 of \cite{coSe4}.
\end{Remark}

As a consequence of these triangle identities, we can express the defining 
isomorphisms \eqref{eq:GV1} of $\calc$ also as follows:

\begin{Proposition} \label{Prop:rules-adj}
Let $\calc$ be a GV-category, and let $\pik_{x,y}^{} \colon 
\Hom(x\tpR y,K) \Rarr\cong \Hom(x,Gy)$ be defined as in \eqref{eq:GV1}. 
 \Enumerate
  \item 
The image of a morphism $f \colon x \Rarr~ Gy$ under the isomorphism 
$\pik_{x,y}\inv$ equals $\evr_{y} \cir (f \tpR \id)$.
  \item  \label{item:apply-G-adj-coev}
The image of a morphism $\xi \colon x \tpR y \Rarr~ K$ under the isomorphism
$\pik_{x,y}^{}$ equals the composite (omitting unitors)
  \be
  x \rarr{\id \tpR \coevr_{y}} x \tpR (y \tpL Gy)
  \rarr{\deltal} (x \tpR y) \tpL Gy \rarr{\xi \tpL \id} Gy \,.
  \label{eq:apply-triangle}
  \ee
\item \label{item:apply-G-adj}
For a morphism $f \colon x \Rarr~ y$, the morphism $Gf$ is equal to the composite
  \be
  \begin{aligned}
  Gy \rarr{\id \tpR \coevr_{x}} Gy \tpR (x \tpL Gx)
  & \rarr{~~\deltal~~} (Gy \tpR x) \tpL Gx
  \Nxl1
  & \rarr{(\id \tpR f) \tpL \id} (Gy \tpR y) \tpL Gx \rarr{\evr_{y}\tpL \id} Gy \,.
  \end{aligned}
  \ee
\end{enumerate}
\end{Proposition}
 
\begin{proof}
The first statement is just the usual expression of the adjunction in terms of the
unit. The second statement follows by combining the first with the triangle identity
\eqref{eq:GV-triangle-1}: Composing \eqref{eq:apply-triangle} with the isomorphism 
from the first part gives, by \eqref{eq:GV-triangle-1}, the identity. The last
statement follows by combining the first two with the definition of the functor $G$ 
via the commuting diagram \eqref{eq:def:Gonmorphisms}.
\end{proof}

\begin{Remark} \label{rem:K-z}
It follows from Lemma \ref{Lemma:GV-coev-gives-all} that we can extend part 2 of 
Proposition \ref{Prop:rules-adj} as follows: for all $x,y,z \iN \calc$ the adjunction
$\Hom(x\tpR y,z) \Rarr\cong \Hom(x,z \tpL Gy)$ can be described as a
composition with the coevaluation $\coevr_{y}$ and an appropriate distributor.
\end{Remark}

We finally have a compatibility of the evaluation and coevaluation morphisms with the monoidal structure of $\calc$:

\begin{Proposition} \label{Prop:comp-coev-tpr}
Let $\calc$ be a GV-category and $x,y\iN\calc$. The right coevaluation
$\coevr_{x \tpr y}$ of the object $x \tpR y$ is equal to the composite
  \be
  \begin{aligned}
  \TI \rarr{\coevr_{x}} x \tpL Gx 
  \rarr{(\id \tpR \coevr_{y}) \tpL \id} (x \tpR (y \tpL Gy)) \tpL Gx
  & \rarr{\deltal \tpL \id} (x \tpR y) \tpL Gy \tpL Gx
  \Nxl1
  & \rarr{~~\cong~~} (x \tpR y) \tpL G(x \tpR y) \,.
  \end{aligned}
 \label{eq:comp-coevs}
  \ee
Analogous expressions are valid for the left coevaluation and for the evaluations. 
\end{Proposition}

\begin{proof}
To show the claim for the right coevaluation, notice that the diagram
  \be
  \begin{tikzcd}[row sep=2.2em]
  \Hom(x \tpR y, z ) \ar{r}{\cong} \ar{d}[swap]{\cong}
  & \Hom(\TI, z \tpL G(x \tpR y)) 
  \\
  \Hom(x, z \tpL Gy) \ar{d}[swap]{\cong} & ~
  \\
  \Hom(\TI, (z \tpL Gy) \tpL Gx) \ar{r}{\cong} 
  & \Hom(\TI, z \tpL (Gy \tpL Gx)) \ar{uu}[swap]{\cong}
  \end{tikzcd}
  \label{eq:coev-compat-tpr}
  \ee
commutes: The isomorphism from the upper to the middle row can be expressed as the
composite
  \be
  \begin{aligned}
  \Hom(x \tpR y,z) & \rarr{\pik\inv_{x\tpr y,G\inv\!z}}
  \Hom((x \tpR y) \tpR G\inv z,K) \rarr{\,\cong\,} \Hom(x \tpR (y \tpR G\inv z),K)
  \\
  & \rarr{\pik^{}_{x,y\tpr G\inv\!z}} \Hom(x,G(y \tpR G\inv z))
  = \Hom(x,z \tpL Gy)
  \end{aligned}
  \ee
with the isomorphisms $\pik$ as defined in \eqref{eq:GV1}; commutativity of 
\eqref{eq:coev-compat-tpr} is a consequence of the resulting cancellation in the 
composition of the two downwards arrows. Now when taking $z \eq x \tpR y$ and the 
identity morphism in the upper left corner, the isomorphism in the upper row yields 
the coevaluation $\coevr_{x \tpr y}$, while by Remark \ref{rem:K-z} the other path 
gives the composite \eqref{eq:comp-coevs}. The other cases are shown analogously. 
\end{proof}


\subsection{Internal Homs and weak module functors} \label{sec:internal-homs-cohoms}

The definition of GV-module category ensures the existence of internal Homs and
coHoms as suitable adjoints of action functors. We now discuss various aspects of
these functors.

\begin{Definition}
Let $\calc$ be a GV-category and $(\calm,\actre)$ a left GV-module category over
$(\calc,\tpr)$, and let $m \iN \calm$. Then the \emph{internal Hom}
$\iHomr(m,-)$ is the right adjoint of the functor from $\calc$ to $\calm$ that
maps objects as $c \,{\xmapsto{~~}}\, c\actrE m$. The \emph{internal coHom} 
$\icoHomr(m,-)$ is the left adjoint of the functor from $\calc$ to $\calm$ that maps
objects as $c \,{\xmapsto{~~}}\, c\actlE m$. In more detail, we have isomorphisms
  \be
  \Hom_\calm(c\actrE m,m') \cong \Hom_\calc(c,\iHomr(m,m'))
  \label{eq:def:iHomr}
  \ee
and
  \be
  \Hom_\calm(m',c\actlE m) \cong \Hom_\calc(\icoHomr(m,m'),c)
  \label{eq:iso:icoHom}
  \ee
for all $c\iN\calc$ and $m,m'\iN \calm$.
\end{Definition}

In view of the definition of the action $\actle$ in terms of $\actre$ (see
Proposition \ref{proposition:actle}), it is not surprising that the internal coHom 
can be expressed in terms of the internal Hom:

\begin{Lemma} \label{lemma:cohom-asHom}
Let $\calm$ be a GV-module category. Then the internal coHom can be expressed as
  \be
  \icoHomr(m,m') \cong G^{-1}(\iHomr(m',m))
  \label{eq:int-coHom-viaHom}
  \ee
for all $m,m' \iN \calm$. 
\end{Lemma}

\begin{proof} 
This follows directly from the definition of $\actle$, 
compare the calculation in \eqref{eq:Ginv-iHom}.
\end{proof}

Being defined via the adjunctions \eqref{eq:def:iHomr} and \eqref{eq:iso:icoHom}, the
internal (co)Hom functors come with canonical (co)evaluation morphisms 
  \be
  \iev_{m,n} \Colon \iHomr(m,n) \actrE m \rarr~ n
  \label{eq:iev}
  \ee
and
  \be
  \icoev_{m,n} \Colon m \rarr~ \icoHomr(n,m) \actlE n \,.
  \ee
These, in turn, give rise to canonical multiplications
  \be
  \imu_{m,n,l}^{} \Colon \iHomr(n,l) \tpR \iHomr(m,n) \rarr~ \iHomr(m,l)
  \label{eq:def:imu}
  \ee
and comultiplications
  \be
  \iDelta_{m,n,l} \Colon \icoHomr(l,m) \rarr~ \icoHomr(n,m) \tpL \icoHomr(l,n)
  \label{eq:def:iDelta}
  \ee
via the compositions (suppressing the associator)
  \be
  \iev_{n,l} \circ (\id_{\iHomr(n,l)} \actrE \iev_{m,n}) \Colon
  \iHomr(n,l) \actrE \iHomr(m,n) \actrE m \rarr{} \iHomr(n,l) \actrE  n \rarr{} l
  \ee
and
  \be
  \begin{aligned}
  (\id \actlE \icoev_{n,l} ) \cir & \icoev_{m,n} \Colon
  \nxl1
  & m \rarr{}
  \icoHomr(n,m) \actlE n \rarr{} \icoHomr(n,m) \actlE \icoHomr(l,n) \actlE l \,,
  \end{aligned}
  \ee
respectively. The following result is standard:

\begin{Lemma} \label{lem:iHom-alg}
The (co)multiplications \eqref{eq:def:imu} and \eqref{eq:def:iDelta} equip 
$\iHomr(m,m)$ with the structure of an associative algebra in $(\calc,\tpr)$ and 
$\icoHomr(m,m)$ with the structure of a coassociative coalgebra in $(\calc, \tpl)$. 
Moreover, for all
$m,l \iN \calm$ the objects $\iHomr(m,l)$ are canonically right $\iHomr(m,m)$-modules
and the objects $\icoHomr(l,m)$ are canonically right $\icoHomr(m,m)$-comodules.
\end{Lemma}

In the case of the regular module category 
we can compute the internal Homs by making use of the isomorphisms
$\Hom(c\tpR d,K) \Cong \Hom(c,Gd) \Cong \Hom(d,G^{-1}c)$. This gives 

\begin{Lemma} 
The internal (co)Homs of the regular module category $\calc$ are given by
  \be
  \iHomr(c,d)\cong G(c\tpR G^{-1} d) = d \tpL Gc \qquad\text{and}\qquad
  \icoHomr(c,d) \cong d \tpr G^{-1} c\,.
  \label{eq:iHom=}
  \ee
\end{Lemma}

Further we have

\begin{Lemma} \label{Lemma:cohom-intHom}
Let $\calm$ be a left GV-module category over a GV-category $\calc$.  
Then the internal Hom satisfies
  \be
  \iHomr(b\actrE m,c\actlE m')\cong c \tpL \iHomr(m,m')\tpL Gb
  \label{eq:iHom-actre-actle}
  \ee
for $b,c \iN \calc$ and $m,m' \iN \calm$.
In particular, the internal Hom is a strong module functor with respect to $\tpl\,$.
Similarly, the internal coHom satisfies
  \be
  \icoHomr(c\actlE m,b\actrE m')\cong b\tpr\icoHomr(m,m')\tpr G^{-1}c
  \label{eq:icoHom-actle-actre}
  \ee
for $b,c \iN \calc$ and $m,m' \iN \calm$.
\end{Lemma}
 
\begin{proof}
Making use of \eqref{eq:HomM.cm.n} and of the easily established adjunction formulas
$ \Hom_\calc(Ga \,{\tpr}
	$\linebreak[0]$
	b,c) \Cong \Hom_\calc(b,a\tpL c)$ and
  \be
  \Hom_\calc(a\tpR G^{-1}b,c) \cong \Hom_\calc(a,c\tpL b)
  \cong \Hom_\calc(a,c\tpL b) \,,
  \ee
for any $b,c,d\iN\calc$ and $m,n\iN \calm$ we obtain the sequence
  \be
  \begin{aligned}
  \HomC(d,\iHomr & (b\actrE m,c\actlE n))
  \cong \Hom_\calm(d\actre b\actrE m,c\actlE n)
  \nxl2
  & \cong \Hom_\calm(G c\actrE d\actrE b\actrE m,n)
  \cong \Hom_\calm((Gc\tpR d\tpR b)\actre m,n)
  \nxl2
  & \cong \Hom_\calc(Gc\tpR d\tpR b,\iHomr(m,n))
  \cong \Hom_\calc(d, c \tpL \iHomr(m,n)\tpL Gb)
  \end{aligned}
  \ee
of isomorphisms. Then \eqref{eq:iHom-actre-actle} follows by the Yoneda lemma.
Similarly, for internal coHoms we get
  \be
  \begin{aligned}
  \Hom_\calc(\icoHomr & (c\actlE m,b\actrE n),d)
  \cong \Hom_\calm(b\actrE n,d\actlE c\actlE m)
  \nxl2
  & \cong \Hom_\calm(n, G^{-1} b\actlE d\actlE c\actlE m)
  \cong \Hom_\calm(n,(G^{-1}b\tpL d\tpL c)\actlE m)
  \nxl2
  &\cong \Hom_\calc(\icoHomr(m,n),G^{-1}b\tpL d\tpL c)
  \nxl2
  & \cong \Hom_\calc(b\tpR\icoHomr(m,n)\tpR G^{-1} c,d) \,,
  \end{aligned}
  \ee
which shows \eqref{eq:icoHom-actle-actre}.
\end{proof}

In the following we use the symbol $\cc{m}$ for the object in $\calm\opp$ that
corresponds to $m \iN \calm$, whereby
e.g.\ the module structures on $\calm\opp$ from Proposition \ref{prop:ractre,ractle}
read
  \be
  \cc{m} \ractre c = \cc{G^{-1}c \actlE m} \qquad\text{and}\qquad
  \cc{m} \ractle c = \cc{G^{-1}c \actrE m} \,.
  \ee

\begin{Proposition} \label{Proposition:weak-str-internal-hom}
Let $\calm$ be a left GV-module category over a GV-category $\calc$. The internal Hom
$\iHomr \colon \calm\opp \Times \calm \Rarr~ \calc$ carries the following (weak)
module functor structures. For $m,n\iN\calm$ and $c\iN\calc$ there are
  \Enumerate
 \item
coherent morphisms
  \be
  c \tpr \iHomr(\cc{m},n)\rarr~ \iHomr(\cc{m},c \actrE n)
  \ee
which endow the functors $\iHomr(\cc{m},-)\colon \calm \To \calc$ with the
structure of lax $\actre$-module functors;
  \item
coherent isomorphisms
  \be
 \iHomr(\cc{m}, c \actlE n) \rarr\cong c \tpL \iHomr(\cc{m},n)
  \ee
which endow the functors $\iHomr(\cc{m},-)\colon \calm \To \calc$ with the
structure of strong $\actle$-module functors;
 \item
coherent morphisms
  \be
  \iHomr(\cc{m},n) \tpr c \rarr~ \iHomr(\cc{m} \ractrE c,n)
  \ee
which endow the functors $\iHomr(-,n)\colon \calm\opp \To \calc$ with the
structure of lax $\ractre$-module functors;
 \item
coherent isomorphisms
  \be
  \iHomr(\cc{m} \ractlE c,n) \rarr\cong \iHomr(m,n) \tpL c
  \ee
which endow the functors  $\iHomr(- ,n)\colon \calm\opp \To \calc$ with the
structure of strong $\ractle$-module functors.
   \end{enumerate}

\noindent
In addition there are the following compatibilities: The corresponding lax and
strong module functor structures in each argument form a 
conjugated pair of lax and oplax mo\-dule functor structures.
In particular the module functors $\iHomr(m,-)\colon \calm \Rarr~ \calc$ and
$\icoHomr(m,-)\colon \calm \Rarr~ \calc$ are GV-module functors.

\smallskip

\noindent
Furthermore we have, for $m,n\iN\calm$ and $c,d\iN\calc$:
 \begin{enumerate} \addtocounter{enumi}{4}
 \item \label{item:comp-both-variables}
The strong module functor structures commute, in the sense that the two obvious
isomorphisms
  \be
  \iHomr(\cc{m} \ractlE c, d \actlE n) \xrightrightarrows[\,~\cong~]{\,~\cong~}
  d \tpL \iHomr(\cc{m},n) \tpL c
  \ee
are equal. Thus $\iHomr \colon \calm\opp \Times \calm \rarr~ \calc$ is a strong
$\actle$-bimodule functor.
 \item
The lax module functor structures commute, in the sense that the two obvious morphisms
  \be
  d \tpr \iHomr(\cc{m},n) \tpr c \xrightrightarrows{~~~~}
  \iHomr(\cc{m} \ractrE c, d \actrE n)
  \ee
are equal. Thus $\iHomr \colon \calm\opp \Times \calm \rarr~ \calc$ is a lax
$\actle$-bimodule functor.
 \item \label{item:lax-iso-module-fun}
The natural isomorphism $\iHomr_{\calm}(c \actrE m, -) 
\Cong \iHomr_{\calc}(c,\iHomr_{\calm}(m,-)) \Cong \iHomr(m,-) \tpL Gc$ obtained from
the isomorphisms \eqref{eq:iHom-actre-actle} is an isomorphism of lax 
$\actre$-module functors.
 \\
Similarly, the natural isomorphism
$\icoHomr_{\calm}(c \actlE m, -) \Cong \icoHomr_{\calc}(c, \icoHomr_{\calm}(m,-))$
is an isomorphism of oplax $\actle$-module functors.
  \end{enumerate}

\noindent
Dual statements hold for the internal coHom.

\end{Proposition}

\begin{proof}
Recall that the functor $\R_{m}^{\actre}\colon {}_{\calc}\calc \Rarr~ {}_\calc\calm$
given in \eqref{eq:R-module-actre} is a strong $\actre$-module functor; its 
right adjoint is the functor $\iHomr(\cc m,-)\colon {}_\calc\calm \To {}_\calc\calc$.
The profunctor transport from Corollary \ref{Corollary:adj-lax} thus provides us
with the lax $\actre$-module functor structure of
$\iHomr(\cc{m},-)$. The strong $\actle$-mo\-du\-le functor structure is obtained by
the \ftransport\ from Lemma \ref{Lemma:Adj-module-fun}. The remaining
structures are obtained with the help of Lemma \ref{Lemma:cohom-intHom} from the
corresponding structures of the internal coHom. By Proposition
\ref{proposition:transports-agree}, the resulting strong and lax module functor
structures correspond to the  lax/oplax functor structure pair of a GV-functor. By
the proof of Lemma \ref{Lemma:cohom-intHom}, the two strong module functor
structures are compatible in the form of statement \ref{item:comp-both-variables}.
 \\
By the transport of structures in a GV-bimodule functor, as described in Lemma
\ref{Lemma:bimodule-functo-transport}, it follows that the lax module functor
structures are compatible as well.
For statement \ref{item:lax-iso-module-fun}, consider the strong $\actre$-module
functors $R_{c}^{\actre} \colon \calc \Rarr~ \calc$ and
$R_{m}^{\actre} \colon \calc \Rarr~ \calm$. The right adjoint lax module functor of
their composite is $(R^{\actre}_{m} R^{\actre}_{c})\ra \Cong \iHomr(c \actrE m,-)$.
By Lemma \ref{lemma:prof-comp}, the isomorphism
$(R^{\actre}_{m} R^{\actre}_{c})\ra \Cong (R^{\actre}_{c})\ra (R^{\actre}_{m})\ra$ is
an isomorphism of lax module functors. This proves the claim concerning the internal
Hom; the statement about the internal coHom follows analogously.
\end{proof}

In the case of $\calc$ as a left module category over itself, the internal Hom
is canonically identified with $\iHomr(x,y) \Cong y \tpL Gx \eq R^{\tpl}_{Gx}(y)$.
The lax $\actre$-module functor structure of $\iHomr(x, -)$ is given in Proposition
\ref{Proposition:weak-str-internal-hom}. On the other hand,
$R^{\tpl}_{Gx}\colon \calc \Rarr~ \calc$ is a strong $\actle$-module functor with
respect to the left action of $(\calc,\tpl)$, hence by
Proposition \ref{Prop:main-transport-part1}
it is also a lax $\actre$-module functor with the conjugated module functor structure.
In this case the lax $\actre$-module functor structures can also serve as a means for
defining the distributors of the GV-category $\calc$ \Cite{Sect.\,4}{fssw}. The so
obtained definition of the distributors is equivalent to the one given in Definition
\ref{Definition:module-distrib}:

\begin{Cor}
The natural isomorphism  $\iHomr(x,-)\Cong R^{\tpl}_{Gx}$ is an isomorphism of
lax $\actre$-mo\-du\-le functors.
\end{Cor}

\begin{proof}
The left adjoint of $R^{\tpl}_{Gx}$ is the functor
$R^{\tpr}_{x} \colon \calc \Rarr~ \calc$. Under the \ftransport, the
strong $\actle$-mo\-du\-le functor structure of $R^{\tpl}_{Gx}$ corresponds to the
strong $\actre$-module functor structure of $R^{\tpr}_{x}$. The profunctor transport
of this structure yields, in turn, the lax $\actre$-module functor structure of
$\iHomr(x,-)$. By Proposition \ref{proposition:transports-agree}, this structure
coincides with the conjugated  lax $\actre$-module functor structure of $R^{\tpl}_{Gx}$.
\end{proof}

Besides the internal Hom $\iHomr \,{\equiv}\, \iH^{\ihr}$ for the regular left module
category, for which in the adjunction $\Hom(c\tpR d,d') \Cong \Hom(c,\iHomr(d,d'))$
the right tensor factor changes place, analogously there is an internal Hom 
$\iHoml$ for the regular right module category\,%
 \footnote{~As compared to \cite{brlv}, in our notation the use of the 
 superscripts $\ihl$ and $\ihr$ is interchanged.}
for which the adjunction keeps the right tensor factor,
  \be
  \Hom(c\tpR d,c') \cong \Hom(d,\iHoml(c,c')) \,.
  \label{eq:def:iHoml}
  \ee
It computes as $\iHoml(c,d)\Cong G^{-1}(Gd \tpR c) \eq G^{-1}c \tpL d$.
Analogously there is a second internal coHom $\icoHoml$.
In accordance with our guiding principle, for every statement involving $\iHomr$ or
$\icoHomr$ there is an analogous statement for $\iHoml$, respectively $\icoHoml$.
This observation applies not just to the regular right module, but likewise to 
every right GV-module category over $\calc$.
In particular there are evaluations
  \be
  \iev'_{m,n} \Colon m \tpR \iHoml(m,n) \rarr{} n \,,
  \label{eq:iev'}
  \ee
and multiplication morphisms
  \be
  \imu'_{m,n,l} \Colon \iHoml(m,n) \tpR \iHoml(n,l) \rarr{} \iHoml(m,l)
  \ee
via
  \be
  \iev'_{n,l} \cir (\iev'_{m,n} \tpR \id_{\iHoml(n,l)}) \Colon
  m \tpR \iHoml(m,n) \tpR \iHoml(n,l) \rarr{} m \tpR \iHoml(m,l) \rarr{} l
  \ee
analogously to to \eqref{eq:iev} and \eqref {eq:def:imu}.


\section{Frobenius algebras and admissible objects} \label{sec:Frob}

\subsection{Algebras and Frobenius algebras in GV-categories}

Throughout this section $(\calc,\tpr) \eq (\calc,\tpr,\one,\alph,l,r,K)$
is a GV-category. Recall the notions of algebras and coalgebras in $\calc$ as 
given in Definition \ref{def:alg,coalg}.
  
\begin{Example}
As seen in Lemma \ref{lem:iHom-alg}, for any object $m$ in a GV-module category
over $\calc$, the internal End $\iHomr(m,m)$ carries the structure of an algebra
in $(\calc, \tpr, \TI)$, while the internal coEnd $\icoHomr(m,m)$ carries the 
structure of a coalgebra in $(\calc, \tpl, K)$.
\end{Example}

\begin{Definition}
In a GV-category $\calc$ a \emph{GV-algebra} is an algebra
in the monoidal category $(\calc,\tpr,1,\alph^\tpr,l^\tpr,r^\tpr)$. 
A \emph{GV-coalgebra} in $\calc$ is a coalgebra
in the monoidal category $(\calc,\tpl,K,\alph^\tpl,
          $\linebreak[0]$
	l^\tpl,r^\tpl)$.
\end{Definition}
   
\begin{Lemma} \label{Lemma:Dual-coalg}
In a GV-category $\calc$, GV-algebras and GV-coalgebras are in bijection
under both the dualizing functor \(G\) and its inverse \(G^{-1}\).
\end{Lemma}

\begin{proof} 
We show that if \((A,\mu,\eta)\) is a GV-algebra (that is an algebra in
$(\calc,\tpr,1)$), then \(G(A)\) naturally inherits the structure of a
coalgebra in \((\calc,\tpl,K)\);
the argument for $G\inv(A)$ is analogous. The reverse statement for
GV-coalgebras then immediately follows from \(G\) being an antiequivalence.
One checks that the morphism $\Delta \colon GA \Rarr{G(\mu)} G(A \tpR A) 
\Rarr{~\cong~} GA \tpL GA$ obtained via the composite isomorphism
  \be
  \Hom(A \tpR A,A) \tocong \Hom(GA,G(A \tpR A)) \tocong \Hom(GA,GA \tpL GA)
  \label{eq:defDelta-viaG}
  \ee
is a coassociative comultiplication.
Furthermore, the morphism $\eps \colon GA \Rarr{\id\tpR\eta} GA \tpR A
\Rarr{\ev_\AA^{}} K$ is a counit for this comultiplication.
\end{proof}

Given the notions of algebras and coalgebras it is natural to also consider Frobenius
algebras, which we do following Definition 2.3.2 of \cite{egge3}.

\begin{Definition} \label{def:GV-Frobenius}
A \emph{GV-Frobenius algebra} in $\calc$ is a quintuple $(A,\mu,\eta,\Delta,\eps)$
such that $(A,\mu,\eta)$ is an algebra in $\calc$, $(A,\Delta,\eps)$ is a GV-coalgebra
in $\calc$, and
  \be
  (\mu \tpL \idA) \circ \distl_{\!A,A,A} \circ (\idA \tpR \Delta) = \Delta \circ \mu
  = (\idA \tpL \mu) \circ \distr_{\!A,A,A} \circ (\Delta \tpR \idA)
  \label{eq:GVfrob}
  \ee
as morphisms in $\Hom(A\tpR A,A\tpL A)$.
 \\[2pt]
A \emph{morphism $f$ of GV-Frobenius algebras} between GV-Frobenius algebras
$(A,\mu^A,\eta^A,\Delta^{\!A},\eps^A)$ and $(B,\mu^B,\eta^B,\Delta^{\!B},\eps^B)$ 
is a morphism $f\colon A \Rarr~ B$ that is compatible with the algebra and 
co\-al\-gebra structures, i.e.\ satisfies $\mu^B \cir (f \tpR f) \eq f \cir \mu^A$,
$f \cir \eta^A \eq \eta^B$, $(f \tpL f) \cir \Delta^{\!A} \eq \Delta^{\!B} \cir f$ 
and $\eps^B \cir f \eq \eps^A$.
\end{Definition}

\begin{Remark}
The two equalities postulated in Equation \eqref{eq:GVfrob} are not independent.
Indeed, it suffices to require that the left and right most expression are the same
morphism. This fact is well known for Frobenius algebras in rigid monoidal
categories. In the present setting the proof is considerably more subtle, since
non-invertible distributors enter. A proof has been given in \cite{demir3}
with the help of a three-dimensional graphical calculus.
\end{Remark}

Note that this definition does not assume that the distributors
$\distl_{\!A,A,A}$ and $\distr_{\!A,A,A}$ are isomorphisms.
Just as with Frobenius algebras in a monoidal category we can also consider
their morphisms.

\begin{Lemma} \label{lem:groupoid}
Every morphism of GV-Frobenius algebras is an isomorphism.
\end{Lemma}

\begin{proof}
Let $f\colon A \Rarr~ B$ be a morphism of GV-Frobenius algebras.
It is not hard to check that the morphism
 \be
  f^- := l_\AA^\tpl \circ \big( (\eps^B\cir\mu^B) \tpL \idA \big) \circ \distr_{B,B,A}
  \circ (\id_B \tpR (f \tpL \idA)) \circ \big( \id_B \tpR (\Delta^A \cir \eta^A))
  \circ (r^\tpr_{\!B})^{-1}
  \ee
is both left and right inverse to $f$. As compared to the calculation for
`ordinary' Frobenius algebras, the only new ingredient is the use of the naturality
of the distributor $\distr$ in place of the naturality of the associator (which
one would normally just suppress).
\end{proof}

Next we establish alternative characterizations of GV-Frobenius algebras that are
equivalent to Definition \ref{def:GV-Frobenius} by preparing some suitable notions.

\begin{Definition}
Let $c$ be an object in a GV-category $\calc$.
A \emph{GV-pairing} on $c$ is a morphism $\kapp_c \iN \Hom(c \tpR c,K)$;
a \emph{GV-co\-pairing} on $c$ is a morphism $\kappo_c \iN \Hom(\one,c \tpL c)$.
A GV-pairing on $c$ is called \emph{non-degenerate} if and only if there exists
a GV-co\-pairing $\kappo_c$ that is side-inverse to $\kapp$, i.e.\ such that
  \be
  \begin{aligned}
  l_c^\tpl \circ (\kapp_c \tpL \id_c) \circ \distl_{c,c,c} & \circ (\id_c \tpr \kappo_c)
  \circ r_c\inv &
  \Nxl1
  & = \id_c = r_c^\tpl \circ
  (\id_c\tpL \kapp_c) \circ \distr_{c,c,c} \circ (\kappo_c \tpr \id_c) \circ l_c\inv
  \end{aligned}
  \label{kapp-kappo}
  \ee
as morphisms in $\End(c)$.
\end{Definition}

It readily follows that the adjunctions $\Hom(c \tpR c, K) \Cong \Hom(c, Gc) \Cong
\Hom(\TI, Gc \tpL Gc)$
provide a bijection between the sets of GV-pairings for $c$ and of GV-copairings for 
$Gc$. Under this bijection, non-degenerate (co)pairings correspond to each other.

\begin{Definition}
Let $A$ be an algebra in a GV-category $\calc$.
An \emph{invariant} GV-pairing on $A$ is a GV-pairing $\kapp_A$ on $A$ such that
  \be
  \kapp_A \circ (\mu \tpr \idA) = \kapp_A \circ (\idA \tpr \mu) \circ \alph_{A,A,A}
  \label{eq:kapp-invariance}
  \ee
as morphisms in $\Hom((A\tpR A) \tpR A,K)$.
\end{Definition}

As is familiar from ordinary Frobenius algebras, invariant GV-pairings
provide an alternative means of characterization:

\begin{Proposition} \label{prop:asscoassFrob=pairFrob}
For an algebra  $A$ in $\calc$ there is a bijection between the 
GV-Fro\-benius algebra structures on $A$ (in the sense of Definition 
\ref{def:GV-Frobenius}) and the  invariant non-de\-ge\-nerate pairings for $A$.
\end{Proposition}

\begin{proof}
Assume that $(A,\mu,\eta,\Delta,\eps)$ is a GV-Frobenius algebra. Define a GV-pairing
$\kapp$ and a GV-copairing $\kappo$ on $A$ by $\kapp \,{:=}\, \eps \cir \mu$ and 
$\kappo \,{:=}\, \Delta \cir \eta$. By the associativity of $\mu$, the GV-pairing
$\kapp$ is invariant. Moreover, the calculation
  \be
  \begin{aligned}
  l_A^\tpl \circ (\kapp \tpL \idA) & \circ \distl_{\!A,A,A} \circ (\idA \tpR \kappo)
  \circ r_\AA\inv &
  \nxl1
  & ~=\, l_A^\tpl \circ (\eps \tpL \idA) \circ (\mu \tpL \idA) \circ \distl_{\AA,A,A}
  \circ (\idA \tpR \Delta) \circ (\idA \tpR \eta) \circ r_\AA\inv
  \nxl1
  & \equ{eq:GVfrob}
  l_A^\tpl \circ (\eps \tpL \idA) \circ \Delta \circ \mu
  \circ (\idA \tpR \eta) \circ r_\AA\inv = \idA \,.
  \end{aligned}
  \ee
shows that $\kapp \eq \kapp_A$ and $\kappo \eq \kappo_A$ satisfy the first of
the equalities \eqref{kapp-kappo} for $c \eq A$. The second of those
equalities follows analogously. Thus the GV-pairing $\kapp$ is non-degenerate.
 \\[2pt]
To show the converse, define, for $(A,\mu,\eta)$  an algebra in $\calc$ and $\kapp$
an invariant non-degenerate GV-pairing on $A$ with side-inverse $\kappo$,
  \be
  \Delta := (\mu \tpL \idA) \circ \distl_{\AA,A,A} \circ (\idA \tpR \kappo) \circ
  r_\AA\inv
  \qquad \text{and} \qquad \eps := \kapp \circ (\idA \tpR \eta) \circ r_\AA\inv .
  \label{eq:Delta1}
  \ee
Then with the help of the invariance property \eqref{eq:kapp-invariance} of 
$\kapp$ one sees that $\eps$ can alternatively be written as
$\eps \eq \kapp \cir (\eta \tpR \idA) \cir l_\AA\inv$. Similarly, for brevity
suppressing for the moment unitors as well as associators, one has
  \be
  \mu = (\kapp \tpL \idA)
  \circ \big( (\idA \tpR \mu) \tpL \idA \big) \circ \distl_{\AA,A\tpr A,A}
  \circ (\idA \tpR \distl_{\AA,A,A}) \circ (\idA \oti \idA \oti \kappo)
  \ee
which, in turn, implies that $\Delta$ can alternatively be written as
  \be
  \Delta = (\idA \tpL \mu) \circ \distr_{\AA,A,A} \circ (\kappo \tpR \idA) \,.
  \label{eq:Delta2}
  \ee
The calculation 
  \be
  \begin{aligned}
  (\idA \tpR \eps) \circ \Delta
  & \equiv (\idA \tpL \kapp) \circ \distr_{\AA,A,A} \circ \big( (\mu \tpL \idA)
  \tpR \eta \big) \circ \distl_{\AA,A,A} \circ (\idA \tpR \kappo)
  \nxl1
  & = \mu \circ (\idA \tpR \eta) = \idA
  \end{aligned}
  \ee
then proves one of the counit properties.
Analogously one shows that $(\eps \tpR \idA) \cir \Delta \eq \idA$.
Next we calculate, denoting for better distinction the two expressions 
\eqref{eq:Delta1} and \eqref{eq:Delta2} for $\Delta$ by different symbols
$\Delta_1$ and $\Delta_2$, respectively,
  \be
  \begin{aligned}
  (\idA \tpL \Delta_1) \circ \Delta_2
  & = (\idA \tpL (\mu \tpL \idA)) \circ (\idA \tpL \distl_{\AA,A,A})
  \circ (\idA \tpL (\idA \tpR \kappo))
  \nxl1
  & ~~~ \circ (\idA \tpL r_A\inv) \circ (\idA \tpL \mu) \circ \distr_{\AA,A,A}
  \circ (\kappo \tpR \idA) \circ l_A\inv 
  \Nxl1
  & = \big[\, \idA \tpL \big[ \big(\mu \cir (\mu \tpR \idA)\big) \tpL \idA \big] \big]
  \circ (\idA \tpL \distl_{A\tpr A,A,A}) \circ \distr_{A,A\tpr A,A\tpl A}
  \nxl1
  & ~~~~ \circ (\distr_{\AA,A,A} \tpR \id_{\AA\tpl A}) \circ ((\kappo \tpR \idA)
  \tpr \kappo) \circ ((\l_\AA\inv \tpR \idA) \tpr \id_\TI) \circ r_\AA\inv ,
  \end{aligned}
  \ee
where we use naturality of the distributors (as well as of the unitors).
Similarly we obtain
  \be
  \begin{aligned}
  (\Delta_2 \tpL \idA) \circ \Delta_1
  & = ((\idA \tpL \mu) \tpL \idA) \circ (\distr_{\AA,A,A} \tpL \idA)
  \circ ((\kappo \tpR \idA) \tpL \idA)
  \nxl1
  & ~~~ \circ (l_A\inv \tpL \idA) \circ (\mu \tpL \idA) \circ \distl_{\AA,A,A}
  \circ (\idA \tpR \kappo) \circ r_\AA\inv
  \Nxl1
  & = \big[ \big[\, \idA \tpL \big( \mu \cir ( \idA \tpR \mu) \big] \tpL \idA \big]
  \circ \big( (\idA \tpL \alph^\tpr_{\AA,A,A}) \tpL \idA \big)
  \nxl1
  & ~~~~ \circ \alph^\tpl_{A,(A\tpr A)\tpr A,A} \circ (\distr_{\AA,A\tpr A,A}\tpL\idA)
  \circ \distl_{\AA\tpl(A\tpr A),A,A}
  \nxl1
  & ~~~~ \circ (\distr_{\AA,A,A} \tpR \id_{\AA\tpl A}) \circ ((\kappo \tpR \idA)
  \tpr \kappo) \circ ((\l_\AA\inv\tpR\idA) \tpr \id_\TI) \circ r_\AA\inv ,
  \end{aligned}
  \ee
where now besides naturality we also make use of the pentagon identity
\eqref{eq:left-distr-coh-1} for the left distributor, specialized to the regular
GV-module category and $x \eq A \tpL A$ and $y\eq z\eq m \eq A$,
followed by the pentagon identity \eqref{eq:right-distr-coh-2} for the right 
distributor with $x \eq y\eq z\eq m \eq A$ ($\tpl$-mul\-ti\-plied with $\idA$).
Invoking now the mixed pentagon identity \eqref{eq:mixed-2} with $x\eq z\eq m \eq A$ 
and $y \eq A\tpR A$, and associativity of $\mu$, it follows that
  \be
  (\idA \tpL \Delta) \cir \Delta
  = \alph^\tpl_{A,A,A} \circ (\Delta \tpL \idA) \cir \Delta \,,
  \ee
i.e.\ that $\Delta$ is a coassociative comultiplication for the $\tpl$-tensor product.
 \nxl2
The validity of the Frobenius relations \eqref{eq:GVfrob} follows similarly as
coassociativity, by expressing the coproduct $\Delta$ suitably either as in
\eqref{eq:Delta1} or as in \eqref{eq:Delta2} and then making use of associativity
of $\mu$ and of the properties of the distributors.
\end{proof}

The result of Proposition \ref{prop:asscoassFrob=pairFrob} can be promoted to
an equivalence of categories (in fact, by Lemma \ref{lem:groupoid}, of groupoids):

\begin{Definition} \label{def:Algkappa}
For $\calc$ a GV-category, $\Alg^{\kappa}_{\calc}$ is the following category: An
object in $\Alg^{\kappa}_{\calc}$ is a pair $(A,\kappa)$ consisting of an algebra $A$
in $\calc$ and an invariant non-degenerate pairing $\kapp$ on $A$. A morphism
$f\colon (A,\kapp_{A}) \Rarr~ (B,\kapp_{B})$ in $\Alg^{\kappa}_{\calc}$ is an algebra
morphism $f\colon A\Rarr~ B$ such that $\kapp_{B} \cir (f \tpR f) \eq \kapp_{A}$.
\end{Definition}

\begin{Proposition} \label{prop:Frob1=2}
The category $\Alg^{\kappa}_{\calc}$ is equivalent to the category of GV-Frobenius 
algebras in $\calc$. In particular, $\Alg^{\kappa}_{\calc}$ is a groupoid.
\end{Proposition}

\begin{proof}
Using the formulas relating the two equivalent definitions of Frobenius algebra
one directly checks that an algebra morphism $f$ in $\Alg^{\kappa}_{\calc}$ is a 
coalgebra morphism (and thus a morphism of Frobenius algebras)
iff $\kapp_{B} \cir (f \tpR f) \eq \kapp_{A}$.
\end{proof}

A third equivalent definition of the notion of a GV-Frobenius algebra is obtained
with the help of the following structure:

\begin{Definition} \label{def:Frobeniusform}
Let $A \iN \calc$ be an algebra in a GV-category $\calc$. A \emph{Frobenius form
for $A$} is a morphism $\lambda\colon A \Rarr~ K$ such that the morphism
  \be
  \begin{aligned}
  \Psil \equiv \Psil_\lambda \Colon A
  & \rarr{r_\AA\inv} A \tpr \TI \rarr{\idA \tpR \coevr_\AA} A \tpr (A \tpL G(A))
  \nxl1
  & \rarr{\distl_{\AA,A,GA}} (A \tpR A) \tpL G(A)
  \rarr{(\lambda \circ \mu) \tpL \id_{G(A)}} K \tpL G(A) \rarr{l^\tpl_{GA}} G(A)
  \end{aligned}
  \label{eq:def:Psil}
  \ee
is invertible.
\end{Definition}

\begin{Proposition} \label{prop:lambda2kapp}
Let $A \iN \calc$ be an algebra in a GV-category $\calc$. If $\lambda$ is a
Frobenius form for $A$, then the morphism
  \be
  \kapp := \lambda \circ \mu
  \label{eq:kapp=lambda.mu}
  \ee
is a non-degenerate invariant GV-pairing on $A$.
Conversely, for $\kapp$ a non-degenerate invariant GV-pairing on $A$, the morphism
  \be
  \lambda := \kapp \circ (\idA \oti \eta)
  \label{eq:lambda=kapp.eta}
  \ee
is a Frobenius form for $A$.
\end{Proposition}

\begin{proof}
Let $\lambda$ be a Frobenius form and define $\kapp$ by \eqref{eq:kapp=lambda.mu}.
By associativity of $\mu$, $\kapp$ is invariant. Further, set
  \be
  \kappo := (\idA \tpr \Psili) \circ \coevr_\AA \Colon \TI \Rarr~ A \tpR A
  \ee
with $\Psili$ the inverse of $\Psil \,{\equiv}\, \Psil_\lambda$.
By naturality of $l^\tpl$ and of $\distl$ we have
  \be
  \idA \equiv \Psili \cir \Psil
  = l_A^\tpl \circ (\kapp \tpL \idA) \circ \distl_{\AA,A,A}
  \circ \big( \idA \oti [ (\idA \tpL \Psili) \cir \coevr_\AA) ] \big) \cir r_\AA\inv ,
  \ee
which means that $\kappo$ satisfies the first of the side-inverseness equalities 
\eqref{kapp-kappo}.
Similarly, using naturality of $r\inv$ and $\distl$ one sees that the identity
$\Psil \cir \Psili \eq \idA$ amounts to $\kappo$ satisfying the second of those
equalities. Thus  $\kappo$ is non-degenerate.
 \\[2pt]
Conversely, let $\kapp$ be a non-degenerate invariant GV-pairing on $A$ and define
$\lambda$ by \eqref{eq:lambda=kapp.eta}. Then we have
$\Psil \eq (\kapp \tpL \id_{GA}) \cir \distl_{\AA,A.GA} \cir (\idA \tpR \coevr_\AA)$.
We claim that $\Psil$ has a two-sided inverse given by 
  \be
  \Psili = (\evr_\AA \tpL \idA) \circ \distr_{GA,A,A} \circ (\id_{GA} \tpR \kappo) \,.
  \label{eq:def:Psili}
  \ee
That this morphism is a right inverse is seen by computing
  \be
  \begin{aligned}
  \Psil \cir \Psili = (\evr_\AA \tpL \idA) \circ \distl_{GA,A,GA} \circ \big( \id_{GA}
  \tpR \big[ (\idA \tpL \kapp) \cir \distr_{\AA,A,A} \cir (\kappo \tpR \idA) \big]
  \tpL \id_{GA} \big) &
  \nxl1
  \circ\, (\id_{GA} \tpR \coevr_\AA) &
  \end{aligned}
  \ee
and noticing that the term in square brackets in this expression equals $\idA$
by the second of the equalities \eqref{kapp-kappo}, so that after invoking the 
snake identity \eqref{eq:GV-triangle-2} we end up with $\id_{GA}$.
That the morphism \eqref{eq:def:Psili} is also a left inverse of $\Psi$ follows
similarly with the help of the snake identity \eqref{eq:GV-triangle-1} and the 
first of the equalities \eqref{kapp-kappo}.
\end{proof}

Together with Proposition \ref{prop:Frob1=2} we have thus arrived at:

\begin{Theorem} \label{thm:4.14}
For a GV-category \C\ the following three groupoids are equivalent: 
 \Enumerate
 \item
The category of Frobenius algebras in \C\ as described in Definition 
\ref{def:GV-Frobenius}.
 \item
The category $\Alg^{\kappa}_{\calc}$ introduced in Definition \ref{def:Algkappa}.
 \item
The category of pairs $(A,\lambda)$ consisting of an algebra $A$ in \C\ and a 
Frobenius form $\lambda$ on $A$. 
 \end{enumerate}
\end{Theorem}

\begin{Remark}
Replacing $\Psil$ in the Definition \ref{def:Frobeniusform} of a Frobenius form by
the morphism
  \be
  \Psir := r^\tpl_{G\inv\AA} \circ (\id_{G\inv\AA} \tpL \kapp)
  \circ \distr_{G\inv\AA,A,A} \cir (\coevl_\AA \tpR \idA) \circ l_\AA\inv
  \label{eq:def:Psir}
  \ee
from $A$ to $G\inv\AA$ one obtains statements analogous to Proposition
\ref{prop:lambda2kapp} in which the equality \eqref{eq:lambda=kapp.eta} is replaced 
by $\lambda \eq \kapp \cir (\eta \oti \idA)$.
\end{Remark}


\subsection{Symmetric Frobenius algebras} \label{sec:symmFrob}

Like in the rigid case there is a notion of symmetric Frobenius algebra,
provided that the GV-category $\calc$ is endowed with a pivotal structure.

\begin{Definition} \label{def:pivotalGV} {\rm \Cite{Def.\,6.1}{boDr}} \\
A \emph{pivotal structure} on a GV-category $\calc$ is a natural family of
isomorphisms
  \be
  \psi_{c,d} \Colon \Hom(c\tpR d,K) \rarr\cong \Hom(d\tpR c, K)
  \ee
for $c,d\iN\calc$ such that $($suppressing associators$)$
  \be
  \psi_{d,c}\circ\psi_{c,d} = \id \qquad \text{and} \qquad
  \psi_{b\tpr c,d}\circ \psi_{c\tpr d,b} \circ \psi_{d\tpr b,c} = \id 
  \label{eq:pivax}
  \ee
for $b,c,d\iN\calc$. A \emph{pivotal GV-category} is a GV-category 
together with a choice of a pivotal structure.
\end{Definition}

As shown by the transformation behavior \eqref{eq:widetildeG**2} of the functor $G^2$
under a change of dualizing object, the notion of a pivotal structure depends on the
choice of dualizing object. Being pivotal is thus not a property of the underlying
monoidal category only. In particular,
the existence of a pivotal structure for one choice of dualizing
object does not guarantee that a pivotal structure also exists for a different choice.

Pivotal structures on $\calc$ are in bijection
with monoidal isomorphisms $\pi\colon \Id_\calc\,{\Rightarrow}\,G^2$ for which
$\pi_K \colon K \Rarr\cong G^2(K)$ coincides with the canonical isomorphism 
$K \Rarr\cong G 1 \eq G^2G^{-1}\TI \Rarr\cong G^2 K$ \Cite{Prop.\,6.7}{boDr}.
Conveniently, this monoidal isomorphism $\pi$ allows us to construct Nakayama
automorphisms:

\begin{Definition} \label{def:symmFrob}
Let $A$ be a Frobenius algebra in a pivotal GV-category $\calc$. Then the left 
and right \emph{Nakayama automorphisms} of $A$ are the invertible endomorphisms
  \be
  \NAl := \Psiri \cir \piv_{G\inv \AA}\inv \cir \Psil \qquad\text{and}\qquad
  \NAr := (\NAl)\inv = \Psili \cir \piv_{G\inv \AA}^{} \cir \Psir
  \ee
of $A$, respectively, with $\Psil$ and $\Psir$ as defined in \eqref{eq:def:Psil}
and \eqref{eq:def:Psir}.
 \\
A \emph{symmetric Frobenius algebra} $A$ in $\calc$ is a Frobenius algebra in $\calc$
for which the Nakayama automorphisms are identities, $\NAl \eq \idA$.
\end{Definition}

\begin{Proposition} 
A Frobenius algebra $A$ in a pivotal GV-category is symmetric if and only if the
corresponding invariant GV-pairing $\kapp$ is symmetric, in the sense that the
equality
  \be
  \ev_\AA \circ \big[ \idA \tpr \big(
  \piv_{G\inv\AA}\inv \cir r^\tpl_{GA} \cir (\kapp \tpL \id_{GA})
  \cir \deltal_{\AA,A,GA} \cir (\idA \tpR \coev_\AA) \big) \big] = \kapp
  \ee
holds.
\end{Proposition}

\begin{proof}
The claim follows directly by using the relation \eqref{eq:def:Psili}
between $\Psil$ and $\kapp$ and the corresponding relation for $\Psir$.
\end{proof}

\begin{Remark}
(i) The proof is fully analogous to the rigid case, which is treated in Section 4 of
\cite{fuSt}.
 \nxl1
(ii) Symmetric Frobenius algebras in \emph{symmetric} linearly distributive categories
have been considered, under the name \emph{Girard monoids}, in \cite{egge3}.
In that case the definition of the algebra being symmetric can be expressed in terms
of the (symmetric) braiding of the category.
\end{Remark}


\subsection{GV-module categories versus categories of modules}

A natural question to ask is under which conditions the category of modules
over an algebra in a GV-category is a GV-module category, and conversely, under which 
conditions a GV-module category comes from an algebra. To investigate the first
issue, recall from Lemma \ref{Lemma:Dual-coalg} that, given an algebra 
$A \iN (\calc,\tpr)$, the object $G(A)$ has a canonical structure of a coalgebra 
in $(\calc,\tpl)$. The analogous result interrelating modules and comodules holds
as well:

\begin{Lemma}
Let \(\calc\) be a GV-category, \(A\iN \calc\) a GV-algebra and \(m\iN \calc\).
There are canonical bijections of
  \Enumerate
  \item \label{itm:ract}
right \(A\)-actions on \(m\) and right \(G(A)\)-coactions on \(m\);
  \item \label{itm:lact}
left \(A\)-actions on \(m\) and left \(G\inv(A)\)-coactions on \(m\);
  \item \label{itm:gract}
right \(A\)-actions on \(m\), left \(G(A)\)-coactions on \(G(m)\) and left
\(G\inv(A)\)-coactions on \(G\inv(m)\);
  \item \label{itm:glact}
left \(A\)-actions on \(m\), right \(G(A)\)-coactions on \(G(m)\) and right
\(G\inv(A)\)-coactions on \(G\inv(m)\).
\end{enumerate}
\end{Lemma}

\begin{proof}
The bijection in Part \ref{itm:ract} is given by the adjunction
\(\Hom(m\tpR A,m)\Cong \Hom(m,m\tpL G(A))\). To see this, first note that the 
comultiplication $\Delta$ considered in the proof of 
Lemma \ref{Lemma:Dual-coalg} coincides with the image of $\mu$ under the sequence 
  \be
  \begin{aligned}
  \Hom(A \tpR A,A) \tocong \Hom(A,A \tpL GA) & \tocong \Hom(\TI,A \tpL GA \tpL GA)
  \Nxl1
  & \tocong \Hom(GA, GA\tpL GA)
  \end{aligned}
  \label{eq:tripleadjunctin}
  \ee
of adjunctions. This follows directly from Part ref{item:apply-G-adj} of Proposition
\ref{Prop:rules-adj} together with Proposition \ref{Prop:comp-coev-tpr}.
Then one sees that the associativity constraint on one side is mapped to the
coassociativity constraint on the other. Part \ref{itm:lact} follows analogously.
Parts \ref{itm:gract} and \ref{itm:glact} are obtained by applying \(G\) or \(G\inv\)
to (co)actions and using the fact that \(G\) is an antiequivalence.
\end{proof}

\begin{Proposition} \label{Proposition:IsGVMod}
Let $\calc$ be a GV-category with equalizers. For any algebra $A \iN \calc$
the category $\calm \eq \moD A$ is a left GV-module category over $\calc$.
The module distributors of $\moD A$ are the  distributors of $\calc$.
\end{Proposition}

\begin{proof}
According to Definition \ref{def:gvmod} we must show that the functors
${-} \actrE m$ and $c \actrE {-}$ admit right adjoints, for $m\iN\calm$ and
$c\iN\calc$, respectively. To see that ${-} \actrE m$ has a right adjoint, define
for $(n,\ohr)$ a right $GA$-comodule and $(x,\rho)$ a left $GA$-comodule
the object $n\oloA x$ as the equalizer
  \be
  \begin{tikzcd}[column sep=2.9em]
  n\oloA x \ar{r}{q_{n,x}} & n\tpl x \ar[yshift=3pt]{r}{\,\ohr\tpL\id_x~}
  \ar[yshift=-3pt]{r}[swap]{\,\id_n\tpL\rho~}
  & n \tpL GA \tpL x \,.
  \end{tikzcd}
  \label{eq:equ-tpl}
  \ee
Then for $m,n \iN \moD A$ consider the diagram
  \be
  \begin{tikzcd}
  \Hom_{A}(c \tpR m, n) \ar{r} 
  \ar{d}
  & \Hom_{\calc}(c, n \oloA G(m)) \ar{d}{\Hom_{\calc}(c,q_{n,G(m)})}
  \\
  \Hom_{\calc}(c \tpR m, n) \ar{r}{\cong} & \Hom_{\calc}(c, n \tpL G(m))
  \end{tikzcd}
  \ee
where the left vertical arrow is obtained from the forgetful functor
and the isomorphism in the bottom row is the adjunction \eqref{eq:195}.
Using the universal property of the equalizer one sees that this diagram
commutes, with the horizontal arrow in the top row an isomorphism. Thus we have
  \be
  \Hom_{A}(c \tpR m, n) \cong \Hom_{\calc}(c, n \oloA G(m))
  \label{eq:HomAvsHomC}
  \ee
for $m,n \iN \moD A$, which shows that ${-} \actrE m$ has 
$\iHomr(m,n) \eq n \oloA G(m)$  as a right adjoint.
 \\[2pt]
It remains to show  that the  action of $c \iN \calc$ on $\moD A$ has a right adjoint.
To this end we first note that for $c \iN \calc$ and 
$(m,\ohr) \iN \moD A$ the composite
  \be
  (c \tpL m) \tpR A \rarr{\distr} c \tpL (m \tpR A) \rarr{\id\tpL\ohr} c \tpL m
  \ee
furnishes a right $A$-module structure on the object $c \tpL m$. {}From the
basic adjunction \eqref{eq:196} we then directly get the desired adjunction \be
  \Hom_{A}(c \tpR m, n) \cong \Hom_{A}(m, G^{-1}c \tpL n) \,.
  \ee
Thus indeed $\calm \eq \moD A$ is a GV-module category.
The statement about the module distributors is obtained by applying Proposition
\ref{proposition:module-fun-distr} to the strong $\actre$-module functor 
$U \colon \moD A \Rarr~ \calc$ that forgets the $A$-module structure.
\end{proof}

Next we ask, conversely, when a GV-module category over a GV-category $\calc$ comes
from an algebra in $\calc$. To analyze this issue, first recall that for an algebra
$A\iN(\calc,\tpr)$ the $\calc$-module category $\moD A$ of right $A$-mo\-du\-les
comes with a forgetful functor
  \be
  U \Colon \calm\simeq \moD A \rarr~ \calc \,.
  \ee
The left adjoint of $U$ is the induction functor $\Ind_A\colon \calc\Rarr~ \calm$,
mapping objects as $x \,{\mapsto}\, x\tpR A$. 
Dual statements hold for comodules and coinduction.

We also note that a $\calc$-\emph{generator} \Cite{Lemma\,2.22}{doSs-v1} of a
GV-module category $\calm$ over $\calc$ is an object $m_0\iN\calm$ such that for
every $m\iN\calm$ 
there exists an object $c\iN\calc$ with an epimorphism $c\actrE m_0\Rarr~ m$.
Analogously we call an object $n_0\iN\calc$ a $\calc$-\emph{cogenerator} of $\calm$
if for every $m\iN\calm$ there exists an object $d\iN\calc$ with 
a monomorphism $m \Rarr~ d\actlE m_0$.

\begin{Lemma}\label{lem:little}
Any algebra $A \eq (A,\mu_A,\eta_A)$ in a monoidal category $(\calc,\tpr)$
is a $\calc$-generator of the $\calc$-module category $\moD A$ of right $A$-modules.
\end{Lemma}

\begin{proof}
For any object $m \iN \moD A$ the $A$-action $m\tpR A \Rarr{\,\rho} m$ on $m$ amounts
to a morphism $\rho \colon \Ind_{A}(U(m)) \Rarr~ m$ which is a morphism in $\moD A$. 
This morphism has the unit as a right inverse and is thus an epimorphism.
\end{proof}
  
Now let $\calm$ be a GV-module category. For any $m_0\iN\calm$,
the object $A_{m_0} \,{:=}\, \iHom(m_0,m_0)\iN\calc$ has a
natural structure of an algebra in $\calc$, and there is a natural functor 
  \be
  \begin{aligned}
  \iHom(m_0,-) \Colon & \calm \rarr~ \moD A_{m_0} \,, 
  \\
  & ~m \xmapsto{~~~} \iHom(m_0,m) \,.
  \end{aligned}
  \ee
Similarly there is a functor $\icoHom(m_0,-)$ from $\calm$ to the category 
$\comoD C_{m_0}$ of $C_{m_0}$-co\-mo\-du\-les.
Our goal is now to specify conditions under which these functors give us
equivalences of GV-module categories.

To get in a position to do so we introduce a particular subclass of objects in 
a GV-module category $\calm$, namely those for which the internal Hom is a strong 
module functor. 
This allows us to generalize the notion of dualizability to the module setting.
We will use the abbreviations
  \be
  \iHomr(m,-) =: Y_m \qquad\text{and}\qquad \icoHomr(m,-) =: W_m \,.
  \label{eq:def:Y,W}
  \ee

\begin{Definition} \label{definition:Admissible}
Let $\calm$ be a GV-module category over a GV-category $\calc$. 
  \Enumerate
  \item
An object $m \iN \calm$ is called \emph{$\tpr$-admissible} if the lax $\actre$-module
functor $Y_m\colon \calm \Rarr~ \calc$ is in fact a strong $\actre$-module functor
and has a right adjoint.
  \item
An object $m \iN \calm$ is called \emph{$\tpl$-admissible} if the oplax 
$\actle$-module functor $W_m\colon \calm \Rarr~ \calc$ is in fact
a strong $\actle$-module functor and has a left adjoint.
  \item 
We denote by $\widehat\calm^\tpr$ and $\widehat\calm^\tpl$ the full
subcategories of $\calm$ on the $\tpr$- and $\tpl\,$-admissible objects, respectively.
  \item
In particular, the subcategories $\widehat{\calc}^{\tpr}$ and $\widehat{\calc}^{\tpl}$
of $\calc$ are the full subcategories on the $\tpr$- and $\tpl\,$-admissible objects
that are obtained by regarding $\calc$ as a left GV-module category over itself.
  \end{enumerate}
\end{Definition}

Note that the functor $Y_m \eq \iHomr(m,-)$ always has a left adjoint and that
$W_m \eq \icoHomr(m,-)$ has a right adjoint.
In general, a strong module functor need not have an adjoint. For instance, a linear 
functor $F \colon \calm \Rarr~ \vect$ from a linear category $\calm$ is a strong 
$\vect$-module functor, but having a left adjoint requires $F$ to be representable.
 
Also note that by \eqref{eq:iHom=} we have
  \be
  Y_\TI = \iHomr_\calc(\TI,-) \cong \Id_\calc \cong \icoHomr_\calc(K,-) = W_K 
  \label{eq:iH1=id=coiHK}
  \ee
and hence
  \be
  \TI \in \widehat{\calc}^{\tpr} \qquad \text{and} \qquad
  K \in \widehat{\calc}^{\tpl} .
  \label{eq:1inCtpr}
  \ee

\begin{Proposition} \label{proposition:admissible-module} 
  \Enumerate
  \item
The subcategories  $\widehat{\calc}^{\tpr}$ and $\widehat{\calc}^{\tpl}$ of a
GV-category $\calc$ are monoidal subcategories.
  \item
Let $\calm$ be a left GV-module category over $\calc$.
By restriction of the action of $\calc$ on $\calm$, the category
$\widehat{\calm}^{\tpr}$ is a left $\widehat{\calc}^{\tpr}$-module category, while
$\widehat{\calm}^{\tpl}$ is a left $\widehat{\calc}^{\tpl}$-module category.
  \end{enumerate}
\end{Proposition}

\begin{proof}
According to Proposition
\ref{Proposition:weak-str-internal-hom}.\ref{item:lax-iso-module-fun},
for all $b,c \iN \calc$ and all $m \iN \calm$ there are natural isomorphisms
$\iHomr_{\calc}(b \actrE c, -) \Cong \iHomr_{\calc}(b, \iHomr_{\calc}(c,-))$ and
$\iHomr_{\calm}(b \actrE m, -) \Cong \iHomr_{\calc}(b, \iHomr_{\calm}(m,-))$
which are isomorphisms of lax $\actre$-module functors. The composition of strong
$\actre$-module functors having a right adjoint is again a strong $\actre$-module
functor having a right adjoint. Hence together with \eqref{eq:1inCtpr}, the first of
those isomorphisms implies that $\widehat{\calc}^{\tpr}$ is monoidal, while the second
isomorphism shows that $\widehat{\calm}^{\tpr}$ is a left
$\widehat{\calc}^{\tpr}$-module category.
The statements about $\widehat{\calc}^{\tpl}$ and $\widehat{\calm}^{\tpl}$ follow in 
an analogous manner from Proposition \ref{Proposition:weak-str-internal-hom} as well.
\end{proof}

We are now in a position to give sufficient conditions for a
GV-module category to be equivalent to a category of modules or comodules.

\begin{Proposition} \label{Prop:i(co)Hom-equivalence}
Let $\calc$ be a finite abelian GV-category and $\calm$ a finite abelian GV-module
over $\calc$.
  \Enumerate
 \item
For $m_0\iN\calm$, the functor $F \,{:=}\, \iHom(m_0,-) \colon \calm \Rarr~ \moD A_{m_0}$
is an equivalence of  $(\calc,\tpr)$-mo\-du\-le categories if and only if
$m_0 \iN \widehat\calm^\tpr$ and $m_0$ is a $\calc$-generator of $\calm$. 
 \item
For $n_0\iN\calm$, the functor $\icoHom(n_0,-) \colon \calm \Rarr~ \comoD C_{m_{0}}$
is an equivalence of $(\calc,\tpl)$-mo\-du\-le categories if and only if 
$n_0 \iN \widehat\calm^\tpl$ and $n_0$ is a $\calc$-cogenerator of $\calm$. 
\end{enumerate}
\end{Proposition}

\begin{proof}
We employ the strategy of \cite[Thm\,7.10.1]{EGno}. We only show the first statement,
the second follows by duality.
Suppose that $m_{0} \iN \widehat\calm^\tpr$ is a $\calc$-generator.
  \Enumerate
 \item
Assume first that $m\iN\calm$ is of the form $m \eq c\actrE m_0$ with some
$c\iN\calc$. Using the fact that owing to $m_0\iN\widehat\calm^\tpr$ the functor
$\iHom(m_0,-)$ is a strong module functor, we then have
  \be
  F(m) = \iHom(m_0,c\actrE m_0) \cong c \tpr \iHom(m_0,m_0)
  = c\actrE A_{m_0} = \Ind_{A_{m_0}}(c) \,.
  \ee
Thus any object of $\calm$ of the form $c\actrE m_0$ is mapped by $F$ to an induced
$A_{m_0}$-module.
Henceforth we drop the index $m_0$ and write $A \,{:=}\, A_{m_0}$.
 \item
For all $m_1\iN\calm$ of the form $c\actrE m_0$ with $c\iN\calc$
and all $m_2\iN\calm$ we have
  \be
  \begin{aligned}
  \Hom_A(F(m_1),F(m_2)) & \cong \Hom_A(\Ind_A(c),F(m_2))
  \cong \Hom_\calc(c,UF(m_2))
  \nxl1
  & = \Hom_\calc(c,\iHom(m_0,m_2)) \cong \Hom_\calm(c\actrE m_0,m_2)
  \nxl1
  & = \Hom_\calm(m_1,m_2) \,.
  \end{aligned}
  \ee
Hence for such objects $m_1$ and $m_2$ the map
$F\colon \Hom_\calm(m_1,m_2) \Rarr~ \Hom_A(F(m_1),F(m_2))$ is an isomorphism.
 \item
By assumption, for every $m_1 \iN \calm$ there is an exact sequence
$c_1\actrE m_0 \Rarr~ c_2 \actrE m_0 \Rarr~ m_1 \Rarr~ 0$ for some $c_1,c_2\iN\calc$.
Since $F$ is exact, the sequence
$F(c_1\actrE m_0) \Rarr~ F(c_2\actrE m_0) \Rarr~ F(m_1) \Rarr~ 0$
is exact as well. Moreover, since for any $m\iN\calm$, the functor $\Hom_\calm(-,m)$
is left exact, the rows of the diagram
  \be
  \hspace*{-0.8em}
  \begin{tikzcd}[row sep=2.5em,column sep=1.3em]
  0 \ar{r}&\Hom_\calm (m_1,m_2) \ar{r} \ar{r}\ar{d}{F}
  & \Hom_\calm (c_1\actrE m_0,m_2) \ar{r} \ar{r}\ar{d}{F}
  & \Hom_\calm (c_2\actrE m_0,m_2) \ar{d}{F}
  \\
  0 \ar{r} & \Hom_A (F(m_1),F(m_2)) \ar{r} &\Hom_A (F(c_1\actrE m_0),F(m_2)) \ar{r}
  & \Hom_A (F(c_2\actrE m_0),F(m_2))
  \end{tikzcd}
  \ee
are exact. By step 2 the second and the third vertical arrow in this diagram are
isomorphisms. Thus (the four-version of) the five-lemma implies that the first
vertical arrow is an isomorphism as well. Thus the map $F\colon \Hom_\calm(m_1,m_2)
\Rarr~ \Hom_A(F(m_1),F(m_2))$ is an isomorphism for arbitrary $m_1,m_2\iN\calm$.
 \item
For any $L\iN\moD A$ there exists an object $c_1\iN\calc$ with a surjection 
$\Ind_A(c_1) \Rarr~ L$ (e.g.\ one can take $c_1\eq U(L)$). Thus there is an exact 
sequence
  \be
  \Ind_A(c_2) \rarr{\,f_A^{}} \Ind_A(c_1) \rarr~ L \rarr~ 0 
  \label{eq:fA}
  \ee
for some $c_2\iN\calc$, by which
$L$ is written as the cokernel of a morphism of induced modules. Further,
using first step 2 and then step 1 we obtain a composite isomorphism
  \be
  \Hom_A(\Ind_A(c_2),\Ind_A(c_1))
  \tocong \Hom_\calm(c_2\actrE m_0,c_1\actrE m_0) \,.
  \label{eq:HomA-HomM}
  \ee
Denote by $f_\calm\iN \Hom_\calm(c_2\actrE m_0,c_1\actrE m_0)$ the image of the
morphism $f_A$ from \eqref{eq:fA} under the linear map \eqref{eq:HomA-HomM}. Let 
$\tilde m\iN \calm$
be the cokernel of $f_\calm$. Since $F$ is exact, we have $F(\tilde m)\Cong L$.
We have thus shown that the functor $F \eq \iHom(m_0,-)$ is essentially surjective.
\end{enumerate}
Suppose now conversely that $m_{0} \iN \calm$ is such that the internal Hom functor
$F \eq \iHom(m_{0},-)$ is an equivalence of $\calc$-module categories. Then $F$ is,
by definition, a strong module functor and thus, as it is an equivalence, $F$ is 
exact. Since the forgetful functor $U \colon \moD A_{m_{0}} \Rarr~ \calc$ is exact,
too, it follows that the composite $UF$ is exact as well, and hence that
$m_{0} \iN \widehat\calm^\tpr$. Further, the object $A_{m_{0}} \eq F(m_{0})$ is 
clearly a $\calc$-generator in $\moD A_{m_0}$; thus, since $F$ is an equivalence, 
$m_0 \iN \calm$ is a $\calc$-generator in $\calm$.
\end{proof}

It can happen that the subcategory $\widehat\calm^\otimes$ does not contain a
$\calc$-generator, e.g.\ it can be zero. In this case the $\calc$-module category
$\calm$ cannot be written as the category of modules over any algebra in $\calc$. As
an illustration, consider the following example:

\begin{Example} 
Consider the semisimple category $\calc$ with two isomorphism classes of simple 
objects, represented by objects $\TI$ and $x$, and with tensor product $x\tpR x \eq 0$
\cite[Example\,2.20]{doSs-v1}.
The category $\calm \eq \vect_\ko$ becomes a $\calc$-module by setting
$x \actrE \ko \,{:=}\, 0$ for the one-dimensional vector space $\ko$. 
Now suppose that we had $\vect_\ko \eq \moD A$ for some unital associative algebra
$A$ in $\calc$. The underlying object of $A$ is a direct sum
$\TI^{\oplus n_1} \,{\oplus}\, x^{\oplus n_x}$ with $n_1 \,{\geq}\, 1$. For the 
regular right $A$-module $A_A$ we thus have $x \actrE A_A \eq x \tpR 
(\TI^{\oplus n_1} \,{\oplus}\, x^{\oplus n_x}) \eq x^{\oplus n_1} \,{\ne}\, 0$, which
contradicts the action $x \actrE 1 \eq 0$ on objects of $\vect_\ko$. Hence $\vect_\ko$
cannot be written as a category of right modules over an algebra in $\calc$.
 \\ 
Next we determine the subcategory $\widehat{\vect}_\ko^\tpr$. Since the relevant
categories are semisimple, the linear functor $\iHom(m,-)$ is exact for every
$m\iN\vect_\ko$.  The adjunctions
  \be
  \begin{aligned}
  & \Hom_\calc(\TI,\iHom(\ko,\ko)) \cong \Hom_\calm(\ko,\ko) \cong \ko
  \nxl2
  \text{and}\qquad & \Hom_\calc(x,\iHom(\ko,\ko)) \cong \Hom_\calm(x\actrE\ko,\ko) = 0
  \end{aligned}
  \ee
show that $\iHom(\ko,\ko) \Cong \TI$. Finally, by comparing the equalities
$\iHom(\ko,x\actrE\ko) \eq \iHom(\ko,0) \eq 0$ and
$x\tpR \iHom(\ko,\ko) \Cong x\tpR \TI \eq x$ we conclude that $\iHom(\ko,-)$
cannot be a strong module functor, and hence that $\widehat\vect_\ko^\tpr$ is zero.
\end{Example}

The following result supplies us with specific 
objects in the subcategory $\widehat\calm^\tpr$:

\begin{Lemma} \label{lem:iHomInd-iHomU}
Let $\calm$ be a GV-module category over $\calc$. Assume that the subcategory
$\widehat\calm^\tpr$ of $\calm$ contains a $\calc$-generator, and thus in
particular is not zero, so that we have an equivalence $\calm \,{\simeq}\, \moD A$
for some algebra $A$ in $\calc$. Let $U$ be the forgetful functor
$U \colon \calm\,{\simeq}\, \moD A \Rarr~ \calc$,
with left adjoint $\Ind_A\colon \calc\Rarr~\calm$.
For any $x\iN \calc$ we have an isomorphism
  \be
  \iHom_\calm(\Ind_A(x),- ) \cong \iHom_\calc(x,U-)
  \ee
of $\calc$-module functors.
\end{Lemma}

\begin{proof}
This follows by direct calculation: we have
  \be
  \Hom_\calc(c,\iHom(\Ind_A(x),m)) \cong \Hom_\calc(c\tpR x, U(m))
  \cong \Hom_\calc(c,\iHom(x,U(m)))
  \ee
for all $x,c \iN \calm$ and all $m \iN \calm$.
\end{proof}

\begin{Cor} \label{cor:IxinMh-if-xinCh}
In the situation of Lemma \ref{lem:iHomInd-iHomU}
we have $\Ind_A(x)\iN\widehat\calm^\tpr$ for every $x\iN\widehat\calc^\tpr$.
\end{Cor}

\begin{proof}
Let $x\iN\widehat\calc^\tpr$. We must show that
$\iHom(\Ind_A(x),-)\colon \calm\Rarr~ \calc$ is an exact strong $\calc$-module
functor. Since the forgetful functor $U\colon \calm \Rarr~ \calc$ is an exact strong
$\calc$-module functor, this follows from Lemma \ref{lem:iHomInd-iHomU}.
\end{proof}

In view of the discussion above, it is desirable to have ``enough'' objects in 
the subcategory $\widehat\calm^\tpr$, and there is the following option for making 
this idea precise:

\begin{Definition} \label{def:algebraic}
We call a GV-module category $\calm$ \emph{algebraic} if it is equivalent to a 
GV-mo\-du\-le category of the form $\moD A$ for some algebra $A$ in $\calc$.
\end{Definition}

In terms of this notion, Proposition \ref{Prop:i(co)Hom-equivalence} amounts to

\begin{Theorem} \label{Thm:i(co)Hom-equivalence}
Let $\calm$ be a GV-module over a GV-ca\-te\-go\-ry $\calc$. $\calm$ is algebraic if 
and only if the subcategory $\widehat\calm^\tpr$ of $\calm$ contains a 
$\calc$-generator $($and is thus in particular not zero$)$.
If in addition $\calc$ and $\calm$ 
are finite abelian, then the objects in $m\iN\calm$ for which $\moD{\iHom(m,m)}$ is 
equivalent to $\calm$ as a $(\calc,\tpr)$-mo\-du\-le are precisely the
$\calc$-generators in $\widehat\calm^\tpr$, while the objects $n\iN\calm$ for which
$\comoD{\icoHom(n,n)}$ is equivalent to $\calm$ as a 
$(\calc,\tpl)$-module are precisely the $\calc$-cogenerators in $\widehat\calm^\tpl$.
\end{Theorem}

Working with algebraic GV-module categories, the statement in Corollary
\ref{cor:IxinMh-if-xinCh}
generalizes to the following sufficient criterion for objects in $\widehat\calm^\tpr$:

\begin{Lemma} \label{lemma:criterion-adm}
Let $\calm\eq \moD A$ be an additive idempotent-complete algebraic module category 
over an additive idempotent-complete GV-category $\calc$ with $A$ an algebra in 
$\calc$.  If for $m \iN \calm$ there exist objects $m' \iN \calm$ and
$x\iN\widehat\calc^\tpr$ with an isomorphism
  \be
  m \oplus m' \cong \Ind_{A}(x) = x \tpR A
  \ee
in $\calm$, then we have $m \iN \widehat\calm^\tpr$. 
\end{Lemma}

\begin{proof}
As is well-known (compare e.g.\ Lemmas 12.17.2 and 12.17.3 of \cite{stacks}),
the direct sum $F \,{\oplus}\, G$ of functors $F$ and $G$ has a right adjoint 
if and only if both $F$ and $G$ have a right adjoint.
Moreover, if $F$ and $G$ are weak module functors, then $F \,{\oplus}\, G$ is a
strong module functor if and only if both $F$ and $G$ are strong module functors.
Thus the statement follows.  
\end{proof}

Note that for $\calc\eq\vect_\ko$, this result reproduces the familiar
condition for $m$ to be projective. 

 \medskip

In the remainder of this section we treat the case of the GV-category 
of finite-dimensional right modules over a commutative algebra
$A$ from Example \ref{Example:comm-alg} in detail and classify the algebraic module
categories over $\moD A$.

\begin{Proposition} \label{Proposition:adm-is-proj}
Let $A$ be a commutative \ko-algebra. For any abelian module category $\calm$ over 
the category $\calc \eq \moD A$ of finite-dimensional $A$-modules we have:
 \\[3pt]
(i) Any $m \iN \walmr$ is projective in $\calm$.
 \\[3pt]
(ii) 
Any $m \iN \walml$ is injective in $\calm$.
\end{Proposition}

\begin{proof}
(i) The functor $\HomC(\TI,-)\colon \calc\Rarr~\calc$ is exact because 
$\TI \eq A_\AA$ is free and thus projective in $\calc \eq \moD A$, and the
functor $\iHom(m,-)\colon \calm\Rarr~\calc$ is exact because $m \iN \walmr$.
Now we have
  \be
  \HomM(m,-) \cong \HomM(\TI \actrE m,-) \cong \HomC(\TI,\iHom(m,-)) \,,
  \ee
i.e.\ $\HomM(m,-)$ is a composite of two exact functors and thus exact.
 \\[2pt]
(ii) 
In view of Remark \ref{rem:K=injective}, $K$ is injective because $\TI$ is
projective. As a consequence,
we can use the same type of argument as in (i), now based on the chain
  \be
  \HomM(-,m) \cong \HomM(-,K\actlE m)  \cong \HomC(\icoHom(m,-),K) 
  \ee
of isomorphisms.
\end{proof}

\begin{Remark}
The argument would fail for $\calc \eq A$-Bimod, since the monoidal unit
${}_AA_\AA$ is, in general, not projective as a bimodule.
\end{Remark}

\begin{Proposition} \label{Proposition:classify-comm-alg-mod}
Let $A$ be a finite-dimensional commutative \ko-algebra. The algebraic module
categories over the category $\moD A$ of finite-dimensional $A$-modules are,
up to equivalence, given by $\moD B$, where $B$ is an algebra extension of $A$:
$B$ is a \ko-algebra together with a unital algebra morphism 
$\iota \colon A \Rarr~ Z(B)$ from $A$ to the center of $B$. 
\end{Proposition}

\begin{proof}
Recall that for a commutative ring $R$, an algebra object in the category $\moD R$
is the same as an $R$-algebra \Cite{Ch.\,VII.3}{MAcl}. Analogously,
an algebra object $B$ in $\calc$ is the same as a \ko-algebra $U(B)$ together with 
a unital algebra morphism $\iota \colon A \Rarr~ Z(U(B))$.
Given such an algebra $B$, the category $\modC B$
of right $B$-modules in $\calc \eq \moD A$ is equivalent to the category $\moD U(B)$ 
of $B$-modules in $\vect$. Indeed, every $m \iN \modC B$ defines an object
$U(m) \iN \moD U(B)$ as follows. As a vector space, $U(m) \eq m$, and the action
of $B$ is the composite $U(m) \oti B \Rarr~ m \otA B \Rarr~ m $ of the
canonical projection and the $B$-action on
$m \iN \modC B$. Conversely, given an object $n \iN \moD U(B)$, by restriction
$n$ is also an $A$-module, hence $n \iN \calc$. Moreover, the action
$n \oti B \Rarr~ n$ is $A$-balanced so that $n$ naturally has the structure of an
object in $\modC B$. The two constructions are readily seen to be functorial and 
to provide an equivalence of categories.
\end{proof}

In the case at hand
we have the following description of the objects in  $ \widehat\calm^\tpr$.

\begin{Proposition}
Let $\calm \eq \moD B$ be a module category over $\calc \eq \moD A$ as in Proposition
\ref{Proposition:classify-comm-alg-mod}. The objects in $\widehat\calm^\tpr$ are
precisely the projective $B$-modules. 
\end{Proposition}

\begin{proof}
If $m$ is projective, then there are objects $n \iN \calm$ and $V \iN \vect$ 
with an isomorphism $m \,{\oplus}\, n \Cong V \otik B \Cong (V \otik A) \otA B$.
Since we have $V \otik A \iN \widehat\calc^\tpr$, it follows from Lemma
\ref{lemma:criterion-adm} that $m \iN \widehat\calm^\tpr$. Conversely, if 
$m \iN \widehat\calm^\tpr$, then by Proposition \ref{Proposition:adm-is-proj},
$m$ is projective. 
\end{proof}


\section{Relative Serre functors}

The prototypical idea behind relative Serre functors is the following. Let $\calc$ be
a monoidal category with a duality functor $D \colon \calc \Rarr~ \calcopp$ --
where $D$ could for instance arise from a rigid duality or from a GV-duality on 
$\calc$. Let furthermore $H \colon \calm\opp \Times \calm \Rarr~ \calc$ be a functor,
such as a Hom functor or an internal Hom. Then a \emph{relative Serre functor}
$S \colon \calm \Rarr~ \calm$ (for $\calc$ with respect to the functor $H$)
is an endofunctor together with a natural family
  \be
  H(n,S(m)) \cong D(H(m,n))
  \label{eq:Serre-cond}
  \ee
of isomorphisms for all $m,n \iN \calm$.
In the situation of our interest we will encounter the weaker variant
\Cite{Def.\,3.2}{oppS} of a \emph{partially defined relative Serre functor}
  \be
  S \Colon \widehat{\calm} \rarr~ \calm 
  \label{eq:partial-Serre}
  \ee
with an isomorphism \eqref{eq:Serre-cond} for $n \iN \calm$  and 
$m \iN \widehat{\calm}$, for $\widehat{\calm} \,{\subset}\, \calm$ a subcategory. 
As we will see, the subcategory $\widehat{\calm}$ relevant to us is the subcategory 
$\widehat\calm^\tpr$ of a GV-module category $\calm$ as introduced in Definition 
\ref{definition:Admissible}. Note that $\widehat\calm^\tpr$ was only designed for
objects whose internal Ends provide algebras representing the module category $\calm$.
It is therefore quite remarkable that the same subcategories appear naturally in the 
definition of relative Serre functors.


\subsection{Internal Homs and representable functors} \label{sec:iHom+repfun}

The internal (co)Homs can be used to discuss (co)representability of module functors.
Recall for a GV-module category $\calm$ over $\calc$ the module functors 
$\R_{m}^{\actre}\colon {}_{\calc}\calc \Rarr~ {}_{\calc}\calm$ and $\R_{m}^{\actle}$
from \eqref{eq:R-module-actre} and \eqref{eq:R-module-actle}.

\begin{Lemma} \label{Lemma-represent-module-fun}
Let $\calc$ be a GV-category and let $\calm$ be a left GV-module category over $\calc$. 
   \Enumerate
  \item
Let $F \colon \calc \To \calm$ be a strong $\actre$-module functor. Then there exists
an object $m \iN \calm$ with a module natural isomorphism $F\Tocong R_{m}^{\actre}$.
The object $m$ is unique up to unique isomorphism. Furthermore, $F$  has a right adjoint.
  \item 
Let $G \colon \calc \To \calm$ be a strong $\actle$-module functor. Then there exists
an object $m \iN \calm$ with a module natural isomorphism $F\Tocong R_{m}^{\actle}$.
The object $m$ is unique up to unique isomorphism. Furthermore, $F$  has a left adjoint.
  \item 
Let $H \colon \calm \To \calc$ be a strong $\actle$-module functor that admits a
left adjoint. Then there exists a (unique up to unique isomorphism)
object $m \iN \calm$ with a module natural isomorphism $H \Tocong \iHomr(m,-)$.
  \item
Let $H \colon \calm \To \calc$ be a strong $\actre$-module functor admitting a
right adjoint. Then there exists a (unique up to unique isomorphism)
object $m \iN \calm$ with a module natural isomorphism $H \Tocong \icoHomr(m,-)$.
   \end{enumerate}
\end{Lemma}

\begin{proof}
To obtain the first statement, set $m \,{:=}\, F(\TI)$. From the module constraint 
$f$ of $F$ we obtain a natural isomorphism 
$f_{c,\TI} \colon F(c) \eq F(c \tpR \TI) \Tocong c \actrE F(\TI) \eq c \actrE m$
for every $c \iN \calc$.
The pentagon axiom for the module constraint implies that this family of 
isomorphisms constitutes a module natural isomorphism $F \Tocong R_{m}^{\actre}$.
Since $R_{m}^{\actre}$ has a right adjoint, we conclude that $F$ has a right adjoint
as well. The second statement follows analogously. 
 \\[2pt]
To show the third statement, we invoke Lemma \ref{Lemma:Adj-module-fun} to conclude
that $H\la \colon \calc \Rarr~ \calm$ is a strong $\actre$-mo\-dule functor, so that 
by the first statement there is a module natural isomorphism 
$H\la \Tocong R_{m}^{\actre}$ for some $m \iN \calm$. But then from the isomorphisms
  \be
  \Hom_{\calc}(c, H(n)) \cong \Hom_{\calm}(H\la(c),n)
  \cong \Hom_{\calm}(c \actrE m,n) \cong\Hom_{\calc}(c, \iHomr(m,n))
  \label{eq:F-internally-rep}
  \ee
we conclude with the Yoneda Lemma that there is a natural isomorphism 
$H \Tocong \iHomr(m,-)$. 
Since all isomorphisms in Equation \eqref{eq:F-internally-rep} are module 
natural isomorphisms, we thus have $H \Cong \iHomr(m,-)$ as module functors.
The last statement follows analogously.  
\end{proof}

In view of Lemma \ref{Lemma-represent-module-fun} the following terminology makes
sense, for $\calc$ a finite GV-category and $\calm$ a GV-module category:

\begin{Definition} \label{Def:irepresentable}
A left exact $(\actle$-module$)$ functor $F\colon \calm\Rarr~ \calc$ is called 
\emph{internally representable} if it is isomorphic $($as a $\actle$-module functor$)$
to a functor of the form
  \be \label{int-rep-actle}   
  \begin{aligned}
  \calm & \rarr~\calc \,,
  \\
  m & \longmapsto \iHomr(m_0,m)
  \end{aligned}
  \ee
for some $m_0\iN \calm$.
 \\[2pt] 
We call a right exact $(\actre$-module$)$ functor $F\colon \calm\Rarr~ \calc$
\emph{internally representable} if it is isomorphic $($as a $\actre$-module functor$)$
to a functor of the form
  \be \label{int-rep-actre}
  \begin{aligned}
  \calm & \rarr~ \calc \,,
  \\
  m & \longmapsto \icoHomr(m_0,m)
  \end{aligned}
  \ee
for some $m_0\in \calm$.
\end{Definition}

Note that a necessary condition for a lax module functor to be internally 
representable is that the module structure is strong.

\begin{Remark}	
Dually to Definition \ref{Def:irepresentable} we call a left exact $(\actle$-module$)$
functor $F\colon \calm\opp\Rarr~ \calc$ \emph{internally corepresentable} if it is
isomorphic $($as a $\actle$-module functor$)$ to a functor from $\calm\opp$ to $\calc$
that maps objects as $m \,{\xmapsto{~~}}\, \iHomr(m,m_0)$ for some $m_0\iN \calm$.
Analogously we call a right exact $(\actre$-module$)$ functor 
$F\colon \calm\opp\Rarr~ \calc$ internally corepresentable if it is isomorphic 
$($as a $\actre$-module functor$)$ to a functor that maps objects as 
$m \,{\xmapsto{~~}}\, \icoHomr(m,m_0)$ for some $m_0\iN \calm$.
The statements about internal representability shown below have obvious analogues
for corepresentable functors.
\end{Remark}

We next examine a suitable subcategory of objects, on which the internal Hom functor 
is internally representable. Recall from Definition \ref{definition:Admissible}
the subcategories $\widehat{\calm}^{\tpr}$ and $\widehat{\calm}^{\tpl}$ of 
$\tpr$- and $\tpl$-admissible objects of a GV-module category $\calm$. 

\begin{Lemma} \label{lemma:ex-Serre}
Let $m \iN \calm$, and let $Y_m$ and $W_m$ be the internal Hom and coHom functors 
as defined in \eqref{eq:def:Y,W}.
  \Enumerate
  \item 
The object $m$ is $\tpr$-admissible if and only if the $\actre$-module functor 
$\iHomr(m,-)$ is internally representable in the form of \eqref{int-rep-actre}, 
i.e.\ if and only if there exist an object $S(m) \iN \calm$ 
and an isomorphism $\phi_{m}\colon \iHomr(m,-) \Tocong \icoHomr(Sm,-)$ of
$\actre$-module functors. If this is the case, then we can choose 
  \be
  S(m) = Y_{m}\ra(K) \,.
  \label{eq:def:Sm}
  \ee
  \item 
The object $m$ is $\tpl$-admissible if and only if the $\actle$-module functor
$\icoHomr(m,-)$ is internally representable as in \eqref{int-rep-actle}, i.e\
if and only if there exist an object $\widetilde{S}(m) \iN \calm$ and an isomorphism 
$\psi_{m}\colon \icoHomr(m,-) \Tocong \iHomr(\widetilde{S}m,-)$ of $\actle$-module
functors. If this is the case, then we can choose 
  \be   
  \widetilde{S}(m) = W_{m}\la(\TI) \,.
  \label{eq:def:tSm}
  \ee		  
    \end{enumerate}
The objects $S(m)$ and $\widetilde{S}(m)$ are unique up to unique isomorphism. 
\end{Lemma}

\begin{proof}
This is a direct consequence of Lemma \ref{Lemma-represent-module-fun}. 
\end{proof}

In particular we obtain natural isomorphisms 
  \be
  \rho_{m}(n)\Colon \Hom(n, S(m)) \cong \Hom(\iHomr(m,n),K)
  \label{eq:Serre-explicit}
  \ee
for $m \iN \widehat{\calm}^{\tpr}$ and $n\iN\calm$.
Since the objects $S(m)$ for $m \iN \widehat{\calm}^{\tpr}$ internally represent
functors, the assignment $m \,{\xmapsto{~}}\, S(m)$ extends to a functor
  \be
  S \Colon  \widehat{\calm}^{\tpr} \rarr{} \calm \,.
  \label{eq:Serre-fun-def-tpr}
  \ee
Analogously we obtain a functor
  \be
  \widetilde{S} \Colon  \widehat{\calm}^{\tpl} \rarr{} \calm \,.
  \label{eq:Serre-fun-def-tpl}
  \ee                           

The functors $S$ and $\widetilde{S}$ are indeed relative Serre functors:

\begin{Cor}
For every $m \iN \widehat{\calm}^{\tpr}$ and every $n \iN \calm$ there is a natural 
isomorphism
  \be
  \widetilde{\phi}_{m}(n) \Colon    \iHomr(n,Sm)  \tocong G (\iHomr(m,n)) \,.
  \label{eq:tildephi-Serre}
  \ee
Analogously, for every $m \iN \widehat{\calm}^{\tpl}$ and every $n \iN \calm$
there is a natural isomorphism
  \be
  \iHomr(\widetilde{S}m,n) \tocong G (\iHomr(n,m)) \,.
  \label{eq:Is-rel-Serre}
  \ee
\end{Cor}

\begin{proof}
By Lemma \ref{lemma:cohom-asHom} we have an isomorphism $\icoHomr(Sm,n) \Cong 
G^{-1}\iHomr(n,Sm)$. The isomorphisms \eqref{eq:tildephi-Serre} and
\eqref{eq:Is-rel-Serre} thus follow directly 
by invoking the isomorphisms from the definitions of $S$ and $\widetilde{S}$.
\end{proof}
  
\begin{Definition}
The functors $S \colon \widehat{\calm}^{\tpr} \Rarr~ \calm$ and
$\widetilde{S} \colon \widehat{\calm}^{\tpl} \Rarr~ \calm$ defined by the formulas
\eqref{eq:def:Sm} and \eqref{eq:def:tSm} are called the \emph{relative Serre functor
$S$ of $\calm$} and the \emph{inverse relative Serre functor $\widetilde{S}$},
respectively. Considering $\calc$ as a left $\calc$-module category defines the 
relative Serre functor $S$ of $\calc$.
\end{Definition}

There is also a variant of a relative Serre functor for $\calc$ when considering 
$\calc$ as right $\calc$-module category.

\begin{Theorem} \label{thm:Serre=equiv}
The relative Serre functors provide an equivalence of categories between the 
subcategories $\widehat{\calm}^\tpr$ and $\widehat{\calm}^{\tpl}$ of $\calm$, i.e.,
slightly abusing notation by keeping the symbols $S$ and $\widetilde{S}$, we have
  \be
  S \Colon \widehat{\calm}^\tpr \rarr\simeq \widehat{\calm}^\tpl \qquad \text{and} \qquad
  \widetilde{S} \Colon \widehat{\calm}^{\tpl} \rarr\simeq \widehat{\calm}^{\tpr} .
  \ee
In particular, for every $m\iN \widehat{\calm}^{\tpr}$ the object $S(m)$ is in the 
subcategory $\widehat{\calm}^{\tpl}$ of $\calm$, and analogously, 
$\widetilde{S}(n) \iN \widehat{\calm}^{\tpr}$ for every $n\iN \widehat{\calm}^{\tpl}$.
\end{Theorem}

\begin{proof}
Let $m$ be $\tpr$-admissible.
According to Proposition \ref{Proposition:weak-str-internal-hom} the two functors
$\iHomr(m,-)$ and $\icoHomr(Sm,-)$ are GV-module functors. Hence it follows from 
Theorem \ref{Thm:bicat-GV-mod} that the $\actre$-module natural isomorphism
$\phi_{m}\colon \iHomr(m,-) \Rarr\cong \icoHomr(Sm,-)$ is also a natural isomorphism 
of $\actle$-module functors.
As a consequence,
for $m \iN \widehat{\calm}^{\tpr}$ one has $S(m) \iN \widehat{\calm}^{\tpl}$,
using that $\icoHomr(Sm,-) \,{\cong}\, \iHomr(m,-)$ is a strong $\actle$-mo\-du\-le
functor, and by the same isomorphism we see that $\widetilde{S}(S(m)) \,{\cong}\, m$.
Analogously it follows that $S(\widetilde{S}(m))\,{\cong}\, m$. Thus $S$ 
and $\widetilde{S}$ define an equi\-va\-lence of categories. 
\end{proof}

\begin{Proposition} \label{proposition:Serre-twisted-ig}
Let $\calc$ be a GV-category and $\calm$ a left GV-module category over $\calc$.
  \Enumerate
  \item 
The relative Serre functor of $\calc$ is canonically a monoidal equivalence
  \be
  S_{\calc} \Colon \widehat{\calc}^{\tpr} \rarr~ \widehat{\calc}^{\tpl}
  \ee
with inverse monoidal functor $\widetilde{S}_{\calc}$.
  \item 
The relative Serre functor $S_{\calm}$ of $\calm$ is a twisted module functor,
in the sense that there are natural isomorphisms
  \be
  S_{\calm}(c \actrE m) \tocong S_{\calc}(c) \actlE S_{\calm}(m)
  \ee
for $c \iN \widehat{\calc}^{\tpr}$ and $m \iN \widehat{\calm}^{\tpr}$.
 \\
Analogously we have $\widetilde{S}_{\calm}(d \actlE n) \Cong 
\widetilde{S}_{\calm}(d) \actlE \widetilde{S}_{\calm}(n)$ for 
$d \iN \widehat{\calc}^{\tpl}$ and $n \iN \widehat{\calm}^{\tpl}$.
  \end{enumerate}
\end{Proposition}

\begin{proof}
To obtain the first statement, first note that the isomorphisms
\eqref{eq:iH1=id=coiHK} imply that $S(\TI) \Cong K$. The
isomorphism $S(x \tpR y) \Cong S(x) \tpL S(y)$ for $x,y \iN \widehat{C}^{\tpr}$ 
follows by applying the Yoneda Lemma to the composite of the isomorphisms
  \be
  \begin{aligned}
  \iHomr(x \tpR y,-) \cong \iHomr(x,\iHomr(y,-)) & \cong \icoHomr(Sx,\icoHomr(Sy,-))
  \nxl1
  & \cong \icoHomr(Sx \tpL Sy,-)
  \label{eq:Serre-mon}
  \end{aligned}
  \ee
of module functors from Proposition \ref{Proposition:weak-str-internal-hom} together
with the isomorphism $\iHomr(x \tpR y, -) \Cong \iHomr
             $\linebreak[0]$
(S(x \tpR y),-)$.
All isomorphisms in the sequence \eqref{eq:Serre-mon} are coherent, hence $S$ 
is a monoidal equivalence. 
 \\
The natural isomorphisms follow from the analogous computation for
$\iHomr(x \actrE m,-)$ with $x \iN \widehat{C}^{\tpr}$ and 
$m \iN \widehat{\calm}^{\tpr}$.
The statements for $\widetilde{S}$ are shown analogously. 
\end{proof}

For the case that the module category is algebraic we will obtain a stronger 
statement in Proposition \ref{proposition:Serre-twisted}.


\subsection{The relative Serre functor of $\calc$}

We now restrict our attention to the important regular case,
already considered in Proposition \ref{proposition:Serre-twisted-ig}, of $\calc$ as
a left GV-module category over itself. Since the internal Homs of $\calc$ can now
be expressed in terms of the monoidal structures, we can compute $S$ as follows. 

\begin{Proposition}\label{proposition:Serre-on-C-1}
Let $\calc$ be a GV-category. For any $c \iN \calc$ the following statements are 
equivalent:
  \Enumerate
  \item
$c \iN \calc$ is $\tpr$-admissible.
  \item 
For all $x,y \iN \calc$ the distributors $c \tpR (x \tpL y) \Rarr~ (c \tpR x) \tpL y$
are isomorphisms.
  \item 
$c$ has a right $\tpr$-rigid dual object.
  \end{enumerate}
\end{Proposition}

\begin{proof}
By applying $G$ to the distributors in statement 2 and invoking Proposition 
\ref{prop:G.dist} we see that these distributors are all isomorphisms if and only if
the oplax $\actle$-module functor $Y_{c}$ is strong. Statement 3 means that there
exists an object $c^\vee$ together with morphisms $c^\vee \tpR c \Rarr~ \TI$ and
$\TI \Rarr~ c \tpR c^\vee$ that obey snake identities. The equivalence of this 
statement with statement 1
is shown in \cite[Prop.\,5.2]{fssw}. Thus all three statements are equivalent. 
\end{proof}

By analogously working with the right module category $\calc_{\calc}$ and/or the 
left exact tensor product $\tpl$ one finds that the statement in Proposition 
\ref{proposition:Serre-on-C-1} has in total four incarnations:
For instance, for $\calc_{\calc}$ with the monoidal structure $\tpr$ we obtain the
statement that the objects $c$ that have a left $\tpr$-rigid dual are those whose
distributors $(x \tpL y) \tpR c \Rarr~ x \tpL (y \tpR c)$ are isomorphisms.
In the case of ${}_{\calc}\calc$ and $\tpl$, we see that the following statements 
are equivalent for every $c \iN \calc$:
   \Enumerate
  \item
$c \iN \calc$ is $\tpl$-admissible.
   \item 
For all $x,y \iN \calc$ the distributors $(c \tpL x) \tpR y \Rarr~ c \tpL (x \tpR y)$
are isomorphisms. 
  \item 
$c$ has a left $\tpl$-rigid dual object. 
  \end{enumerate}
In view of the third characterization of admissibility in Proposition
\ref{proposition:Serre-on-C-1}, by applying $G$ or $G\inv$ to the evaluation and
coevaluation one sees that for any $c\iN \widehat{\calc}^\tpr$ the object $G(c^\vee)$
is a right $\tpl$-rigid dual of $G(c)$. Analogous statements hold for the other 
three incarnations of Proposition \ref{proposition:Serre-on-C-1}.

\begin{Cor}
The functor $G^{2}$ preserves the $\tpr$- as well as the $\tpl$-admissible objects;
that is, the functors
  \be
  G^{2} \Colon \widehat{\calc}^{\tpr} \rarr~ \widehat{\calc}^{\tpr}
  \qquad \text{and} \qquad
  G^{2} \Colon \widehat{\calc}^{\tpl} \rarr~ \widehat{\calc}^{\tpl}
  \ee
are monoidal equivalences. 
\end{Cor}

\begin{proof}
Since $G^{2}$ is a monoidal equivalence for $\tpr$ as well as of $\tpl$, it preserves
the class of objects that have a right  $\tpr$-rigid dual, as well as the class of 
objects that have a left $\tpl$-rigid dual. 
\end{proof}

\begin{Lemma}\label{lemma:Serre-on-C-2}
Let $\calc$ be a GV-category for which the three equivalent conditions in
Proposition \ref{proposition:Serre-on-C-1} are satisfied. Then  $Y_{c}$ has 
  \be
  Y_{c}\ra(z)= z \tpL (G^{2}(x) \tpr K) 
  \label{eq:Yra}
  \ee
as a right adjoint functor. Moreover, the right $\tpr$-rigid dual
of any $c \iN \calc$ is given by 
  \be
  c^\vee \cong \TI \tpL Gc 
  \label{eq:cvee=1..Gc}
  \ee
and we have 
  \be
  y \tpL Gc  \cong y \tpr c^\vee
  \label{eq:y..Gc=y.cvee}
  \ee
as well as
  \be
  c \tpL y \cong c \tpR (\TI \tpL y)
  \label{eq:c..x=x.(1..x)}
  \ee
for every $c \iN \widehat{\calc}^\tpr$ and $y \iN \calc$.
\end{Lemma}

\begin{proof}
The isomorphism
$Y_{c}(y) \eq y \tpL Gc \Cong (y \tpR \TI) \tpL Gc \Cong y  \tpr (\TI\tpL Gc)$ 
implies that $Y_{c}$ has the functor \eqref{eq:Yra} as a right adjoint.
The expression \eqref{eq:cvee=1..Gc} for $c^\vee$ is shown in
\cite[Prop.\,5.2]{fssw}. The isomorphism \eqref{eq:y..Gc=y.cvee} follows,
for instance, from the fact that the two functors ${-} \tpL Gc$ and ${-} \tpr c^\vee$
are both right adjoint to ${-} \tpr c$. Finally, \eqref{eq:c..x=x.(1..x)} is the 
special case $x \eq 1$
of the isomorphism in the second statement of Proposition 
\ref{proposition:Serre-on-C-1}.
\end{proof}

Next we describe the relative Serre functor of $\calc$ explicitly:

\begin{Proposition} \label{proposition:coh-Serre}
The relative Serre functor $S$ of $\calc$ is given by
  \be
  S(c) \cong G^{2}(c) \tpr K
  \label{eq:partial-serre-expl}
  \ee
for $c \iN \widehat{\calc}^\tpr$.
This natural isomorphism is monoidal in the sense that the diagram 
  \be
  \begin{tikzcd}[row sep=0.5em]
  S(c \tpR d) \ar{r}{\cong} \ar{dddd}[swap]{\cong}
  & G^{2}(c \tpR d) \tpr K \ar{dr}{\cong}
  \\
  ~ & ~ & G^{2}(c) \tpR G^{2}(d) \tpR K \ar{dd}{\cong}
  \\ ~ \\
  ~ & ~ & G^{2}(c) \tpR (K \tpL (G^{2}(d) \tpR K)) \ar{dl}{\delta}
  \\
  S(c) \tpL S(d) \ar{r}{\cong} & (G^{2}(c) \tpR K) \tpL (G^{2}(d) \tpR K)
  \end{tikzcd}
  \label{eq:coherence-Serre}
  \ee
with $\delta \eq \deltal_{G^2c,K,G^2d\tpr K}$ (which is an isomorphism)
commutes for all $c,d \iN \widehat{\calc}^\tpr$. The inverse of $S$ is given by
  \be
  \widetilde{S}(c) \cong G^{2}(c) \tpL \TI
  \label{eq:widetilde-Serre}
  \ee
for $c \iN \widehat{\calc}^\tpl$.
\end{Proposition}

\begin{proof}
Note that because of $G^{2}(c) \iN \widehat{\calc}^\tpr$ the distributor
$\deltal_{G^2c,K,G^2d\tpr K}$ is an isomorphism. The expression
\eqref{eq:partial-serre-expl} for the relative Serre functor follows from 
$S(c) \eq Y_{c}\ra (K)$.
The monoidality follows from \eqref{eq:Serre-mon}, which expresses the monoidal 
structure of $S$; this is seen as follows. For $c \iN \widehat{\calc}^{\tpr}$
we have by Equation \eqref{eq:tildephi-Serre} the isomorphism
  \be
  G(\iHomr(c,x)) \cong \iHomr(x, Sc) 
  \label{eq:Serre-reformulated}
  \ee
for every $x \iN \calc$. If we insert $Sc \Cong G^{2}c \tpR K$, the isomorphism
\eqref{eq:Serre-reformulated} computes as
  \be
  G (x \tpl c) \cong Gc \tpR Gx \cong  (G^{2}c \tpr K) \tpL Gx \cong Sc \tpL Gx \,, 
  \label{eq:compute-Serre-ref}
  \ee
where the second isomorphism is a distributor. For the isomophisms 
\eqref{eq:compute-Serre-ref} the coherence diagram \eqref{eq:Serre-mon} translates 
to the diagram \eqref{eq:coherence-Serre}. 
The expression \eqref{eq:widetilde-Serre} for $\widetilde{S}$ follows analogously. 
\end{proof} 

\begin{Remark}
Also, upon a change of dualizing object from $K$ to
$\widetilde{K} \eq g \tpR K$ with invertible $g \iN \calc$, combining the formula 
\eqref{eq:partial-serre-expl} for the relative Serre functor with the resulting
change \eqref{eq:widetildeG**2} in $G^2$ shows that $S$ changes to
  \be
  S'(c) = \widetilde{G}^{2}(c) \tpR \widetilde{K}
  = g \tpR G^{2}c \tpR g^{-1} \tpR g \tpR K \cong g \tpR G^{2}c \tpR K
  \cong g \tpR S(c) \,. 
  \ee
It follows e.g.\ that if $S(c) \Cong  c$, then there exists an isomorphism 
$S'(c) \Cong c$ if and only if $g \tpR c \Cong c$.
\end{Remark}

Since  $G^{2} \colon \calc \Rarr~ \calc$ is a monoidal equivalence both for $\tpr$
and for $\tpl$, by composing $S$ with its inverse we obtain directly

\begin{Cor}
Let $\calc$ be a GV-category.
The functor
  \be
   {-} \tpR K \Colon \widehat{\calc}^\tpr \rarr~ \widehat{\calc}^\tpl
  \ee
is a monoidal equivalence with inverse functor $ {-} \tpL \TI$. 
\end{Cor}

We can now conclude that for algebraic module categories the relative Serre functor 
is twisted by $G^{2}$, a situation familiar from the rigid case 
\cite[Lemma\,4.23]{fScS2}:

\begin{Proposition} \label{proposition:Serre-twisted}
Let $\calm$ be an algebraic left $\calc$-module category. Then 
$S_{\calm} \colon \widehat{\calm}^\tpr \Rarr~ \calm$ is a twisted module functor
in the sense that there are coherent natural isomorphisms
  \be
S_{\calm}(c \actre m) \cong G^{2}c \actre S_{\calm}(m)
  \ee
for all $c \iN \widehat{\calc}^\tpr$ and $m \iN \widehat{\calm}^\tpr$. 
\end{Proposition}

\begin{proof}
By  Proposition \ref{Proposition:IsGVMod} the module distributors of $\calm$ can be
obtained from the distributors of $\calc$, which together with Proposition
\ref{proposition:Serre-on-C-1} implies that for all $c \iN \widehat{\calc}^\tpr$, all
$d \iN \calc$ and all $m \iN \calm$ the distributor
  \be
  c \actrE (x \actlE m) \rarr~ (c \tpR x) \actlE m 
  \ee
is an isomorphism. But then we obtain with Proposition \ref{proposition:Serre-twisted-ig}
for all $c \iN \widehat{\calc}^\tpr$ and $m \iN \widehat{\calm}^\tpr$ an isomorphism
  \be
  \begin{aligned}
  S_{\calm}(c \actrE m) & \cong S_{\calc}(c) \actlE S_{\calm}(m)
  \cong (G^{2}c \tpR K) \actlE S_{\calm}(m)
  \Nxl1
  & \cong G^{2}c \actrE (K \actlE S_{\calm}(m)) \cong G^{2}c \actrE S_{\calm}(m) \,.
  \end{aligned}
  \ee
It follows from Proposition \ref{proposition:coh-Serre} that the composite 
isomorphism is also coherent. 
\end{proof}

The following example shows that the category of admissible objects is in general 
not abelian, even if $\calc$ is:

\begin{Example} \Cite{Lemma\,5.5}{fssw}\\
Let $A$ be a finite-dimensional \ko-algebra and $\calc \eq A\bimod$ the corresponding 
GV-category. 
The following statements are equivalent for ${}_{A}M_{A} \iN \calc$:
  \begin{itemize}
  \item[(i)] 
$M$ is  admissible, i.e.\ $M \in \widehat{\calc}^{\tpr}$.
  \item[(ii)] 
$M$ has a $\ota$-right dual.
  \item[(iii)] 
$M_{A}$ is  projective as a right $A$-module.
  \end{itemize}
\end{Example}


\subsection{Relative Serre functors and Frobenius algebras} \label{sec:Serre+Frob}

Finally, generalizing an argument from \cite{shimi20}, we show that a trivialization
of the relative Serre functor equips the internal End algebra with a Frobenius structure. 
Recall from Equation \eqref{eq:tildephi-Serre} the isomorphism 
$\widetilde{\phi}_{m}(n) \colon \iHomr(n,S(m)) \Tocong G(\iHomr(m,n)).$

\begin{Lemma}
Let $m \iN \widehat{\calm}^{\tpr}$. For all objects $c\iN\calc$ and $m,n\iN\calm$
the diagram
  \be
  \begin{tikzcd}
  \Hom_{\calc}(c, \iHomr(n,S(m))) \ar{r}{} \ar{dd}[swap]{\widetilde{\phi}_{m}(n)}
  & \Hom_{\calm}(c \actrE n, S(m)) \ar{d}{\rho_{m}(c \actrE n)}
  \\
  ~ & \Hom(\iHomr(m, c \actrE n),K) \ar{d}{\cong}
  \\
  \Hom(c, G(\iHomr(m,n))) \ar{r}{} & \Hom(c \tpr \iHomr(m,n),K)
  \end{tikzcd}
  \label{eq:Serre-main-diag}
  \ee
commutes.
\end{Lemma}

\begin{proof}
{}From the proof of Lemma \ref{Lemma-represent-module-fun} we obtain a commuting 
diagram
  \be
  \begin{tikzcd}[row sep=2.3em]
  \Hom(\iHomr(m,n),c) \ar{d}[swap]{(\phi_{m}(n))^{*}} \ar{r}{\cong}
  & \Hom(n, Y_{m}\ra (c))\ar{d}{\cong} 
  \\
  \Hom(\icoHomr( Y_{m}\ra (K),n),c) \ar{r}{\cong} & \Hom(n, c \actlE Y_{m}\ra (K))
  \end{tikzcd}
  \label{eq:Serre-diag}
  \ee
Inserting the definition of the strong $\actle$-module functor structure of 
$Y_{m}\ra$, it follows that the outer hexagon in the diagram
  \be
  \begin{tikzcd}[row sep=2.3em,column sep=2.1em]
  \Hom(\icoHom(Sm,n),c)\ar{rr}{\cong} \ar{dr}{G^{-1}} \ar{dd}[swap]{\phi_{m}(n)^*}
  & ~ & \Hom(n, c \actlE S(m)) \ar{dd}{\cong} 
  \\
  ~ & \Hom(Gc, \iHomr(n,Sm)) \ar{d}[swap]{\widetilde{\phi}} \ar{dr}{\cong} & ~
  \\
  \Hom(\iHomr(m,n),c) \ar{d}[swap]{\cong} \ar{r}{G^{-1}} & \Hom(Gc, G(\iHomr(m,n)))
  & \Hom(Gc \actrE n, Sm)\ar{d}{\cong} 
  \\
  \Hom(Gc, \tpr \iHomr(m,n),K) \ar{rr}{\cong} \ar{ur}[swap]{\cong} & ~
  & \Hom(\iHomr(m, Gc \actrE n),K)
  \end{tikzcd}
  \label{eq:27}
  \ee
commutes. It is also readily seen that the inner triangle and the two inner
quadrangles commute. Thus the pentagon in the lower right corner commutes as well;
this yields the diagram \eqref{eq:Serre-main-diag}. 
\end{proof}

Following \Cite{Def.\,3.7}{shimi20} we give

\begin{Definition}
Let $m \iN \widehat{\calm}^{\tpr}$. The \emph{trace of $m$} is the composite morphism
  \be
  \tr_{m} \Colon \iHomr(m,S(m)) \rarr{\widetilde{\phi}_{m}(m)}
  G(\iHomr(m,m)) \rarr{G(u_{m})} G(1) = K 
  \label{eq:trace-m}
  \ee
with the unit $u_{m}\colon \TI \Rarr~ \iHomr(m,m)$.
\end{Definition}

We then obtain the following analogue of \Cite{Lemma\,3.8}{shimi20}:

\begin{Lemma} \label{Lemma:trace-comp-ev}
Let $m \iN \widehat{\calm}^{\tpr}$. With respect to the multiplication
$\imu{}_{m,n,k} \colon \iHomr(n,k) \tpR \iHomr(m,n) 
        $\linebreak[0]$
	\rarr{} \iHomr(m,k)$
we have a commuting diagram
  \be
  \begin{tikzcd}[row sep=2.9em,column sep=3.3em]
  \iHomr(n,Sm) \tpr \iHomr(m,n) \ar{r}{\imu_{m,n,Sm}}
  \ar{d}[swap]{\widetilde{\phi}_{m}(m)} & \iHomr(m,Sm) \ar{r}{\tr_{m}} & K
  \\
  G(\iHomr(m,n))\tpR \iHomr(m,n)
  \ar[yshift=-3pt,xshift=4pt]{urr}[swap]{\ev_{\iHomr(m,n)}} & ~ & ~
  \end{tikzcd}
  \ee
\end{Lemma}

\begin{proof}
This follows, as in \cite{shimi20}, by chasing the identity morphism in the 
upper left corner in \eqref{eq:Serre-main-diag} for $c \eq \iHomr(n,Sm)$ through 
the diagram. The downward-right-path yields the element 
$\ev_{\iHomr(m,n)} \cir \widetilde\phi_m(m) \iN \Hom(\iHomr(n,Sm) \tpR \iHomr(m,n),K)$.
Using Lemma \ref{Lemma:trace-comp-ev}, the other path in \eqref{eq:Serre-main-diag}
gives the morphism $\imu{}_{m,n,Sm}\iN \Hom(\iHomr(n,Sm) \tpR \iHomr(m,n),K)$.
\end{proof}

This leads us to the following analogue of Theorem 3.14 of \cite{shimi20}:

\begin{Theorem} \label{Theorem:Frob-Serre-trivial}
Let $m \iN \widehat{\calm}^{\tpr}$ be such that the object $S(m)$ is isomorphic
to $m$. Then for every choice of isomorphism $p \colon m \Rarr~ Sm$ in $\calm$,
$(\iHomr(m,m), \mu \eq\imu{}_{m,m,m})$
is a GV-Frobenius algebra in $\calc$ with Frobenius form
\be
\label{eq:lambda-from-rho}
  \lambda \Colon \iHomr(m,m) \rarr{\iHomr(m,p)} \iHomr(m,Sm) \rarr{\tr_{m}} K \,.
  \ee
\end{Theorem}

\begin{proof}
The proof given in \cite{shimi20} generalizes: By Lemma \ref{Lemma:trace-comp-ev}
for $A \eq \iHomr(m,m)$, the composite morphism
$\lambda \cir \mu \colon A \tpR A \Rarr~ K$ is equal to the composite
  \be
  \begin{aligned}
  \iHomr(m,m) \tpR A & \rarr{\iHomr(m,p)\tpr \id} \iHomr(m,Sm) \tpR A
  \nxl1 &
  \rarr{\widetilde{\phi}_{m}(m) \tpr \id} G(\iHomr(m,m)) \tpR \iHomr(m,m)
  \rarr{\ev_{\iHomr(m,m)}} K \,.
  \end{aligned}
  \ee
The first two morphisms in this composite are isomorphisms, and the last is the 
duality pairing, which is non-degenerate. Thus in total $\lambda$ is a Frobenius form
on $\iHomr(m,m)$.
\end{proof}

\begin{Remark}
Suppose that there is an isomorphism $p\colon m\Rarr~ S(m)$, i.e.\ that $m$ is a
fixed point under the action of the relative Serre functor. Then $m$ is 
an object in both $\widehat\calm^\tpr$ and $\widehat\calm^\tpl$. (But we do better
than just intersecting the classes of objects of these subcategories, because we 
also trivialize the relative Serre functor.) 
By Theorem \ref{Theorem:Frob-Serre-trivial}, for such a Serre fixed point $m$
the algebra $\iHomr(m,m)$ is Frobenius. In general it is not symmetric,
though -- there need not even exist a pivotal structure on $\calc$.
\end{Remark}

In the proof of Theorem \ref{Theorem:Frob-Serre-trivial}
no coherence requirement on the isomorphism $p$ needs
to be imposed. From Lemma \ref{Lemma:trace-comp-ev} and the naturality of the 
multiplication $\imu{}_{m,n,l}$ we can directly deduce a stronger non-degeneracy
of the Frobenius pairing:

\begin{Cor}
Let $m \iN \widehat{\calm}^{\tpr}$ be an object with an isomorphism $S(m) \Cong m$ and
a corresponding Frobenius form $\lambda \colon \iHomr(m,m) \Rarr~ K$. Then 
for every $n \iN \calm$ the pairing
  \be
  \iHomr(n,m) \tpr \iHomr(m,n) \rarr{\imu{}_{m,n,m}} \iHomr(m,m) \rarr{\lambda} K
  \label{eq:gen-non-deg}
  \ee
is non-degenerate.
\end{Cor}

This, in turn, characterizes the Frobenius structures on the algebra $\iHomr(m,m)$
that we obtain this way:

\begin{Proposition} \label{prop:serre-iso}
Let $m \iN \widehat{\calm}^{\tpr}$ and $\lambda \colon \iHomr(m,m) \Rarr~ K$ 
a morphism, such that the pairing \eqref{eq:gen-non-deg} is non-degenerate for every
$n \iN \calm$. Then there exists an isomorphism $p\colon m \Rarr\cong Sm$ such that
$\lambda$ is constructed from $p$ as the composite \eqref{eq:lambda-from-rho}.
\end{Proposition}

\begin{proof}
The statement follows directly from the Yoneda lemma applied to the composite
isomorphism $\iHomr(n,m) \Rarr\cong G(\iHomr(m,n)) \Rarr\cong \iHomr(n,Sm)$, in which
the first isomorphism comes from the non-degeneracy of \eqref{eq:gen-non-deg} and
the second from the definition of the relative Serre functor. 
\end{proof}

\begin{Cor} \label{cor:isodist-Frob}
If $x \iN \calc$ has a rigid right dual 
and there is an isomorphism $x \,{\cong}\, G^{2}(x) \tpR K$, then the internal End
$\iHomr(x,x)\,{\cong}\, x \tpL Gx$ has the structure of a GV-Frobenius algebra. 
\end{Cor}

Let us finally turn our attention to the situation that the GV-category $\calc$ 
admits a pivotal structure. Assume that a specific pivotal structure has been fixed.

\begin{Cor} \label{cor:pivotaldual}
Let $\calc$ be a GV-category with a pivotal structure.
Then for any $c \iN \widehat{\calc}^{\tpr}$ the rigid right dual obeys
  \be
  c^{\vee} \cong \TI \tpL Gc \cong G(c \tpR K) \,.
  \label{cvee=1..Gc}
  \ee
It follows that $c$ has $G(c)$ as a rigid right dual if  and only if there is an 
isomorphism $S(c) \Cong c$. In this case the internal End $\iHomr(c,c)$ has
the structure of a Frobenius algebra in $\calc$.
\end{Cor}

\begin{proof}
The presence of a pivotal structure implies that the relative Serre functor
\eqref{eq:partial-serre-expl} on $\widehat{\calc}^{\tpr}$
is given by $S(c) \Cong c \tpR K$. By Lemma \ref{lemma:Serre-on-C-2}
this provides the isomorphism \eqref{cvee=1..Gc}. It follows
that there is an isomorphism $S(c) \Cong c$ if and only if $c \tpR K \Cong c$ which,
in turn, is the case if and only if $c^{\vee} \Cong G(c)$. 
\end{proof}

Analogously, by \eqref{eq:widetilde-Serre} in a pivotal GV-category we have
$\widetilde{S}(c) \Cong c \tpL \TI$.

\begin{Example}
These observations connect nicely with the conditions given in \Cite{Def.\,3.8}{garW}
for the subcategory relevant to the classification of boundary conditions in the
so-called logarithmic triplet model of conformal field theory \Cite{Sect.\,3}{garW}.
There the case of a pivotal GV-category with monoidal isomorphism
\(\pi\colon \id_\calc \Rarr~ G^2\), as in Section \ref{sec:symmFrob}, is considered.
Denote the image of \(\pi\) under the GV-adjunction \(\pik\) of \eqref{eq:GV1} by
  \be
  f_c=\pik_{c,Gc}(\pi_c) \Colon c\tpr Gc \rarr~ K \,. 
  \ee
The subcategory classifying boundary conditions is constructed in two
stages. First there is the full subcategory \(\calc^r\) consisting of all
objects admitting both a left and a right dual (this is a two-handed version
of our subcategories of admissible objects). Within \(\calc^r\) there is
then a further restriction to the full subcategory \(\calc^b\) of all objects \(c\iN
\calc^r\) for which \(Gc\) is both a left and a right dual for \(c\) (an additional
non-degeneracy condition, which is not relevant here, is also imposed). It was
then noted that for \(c\iN\calc^b\) the internal End \(\iH(c,c)\) is isomorphic to
\(c\tpR Gc\) and shown in \Cite{Thm.\,3.10}{garW} that the composition of such an
isomorphism with \(f_c\) is non-degenerate (the choice of isomorphism was left 
implicit in \cite{garW}, but from Corollary \ref{cor:pivotaldual} we know that
such a choice is equivalent to a choice of isomorphism \(S(c) \Cong c\)). This
matches Corollary \ref{cor:isodist-Frob} asserting that such internal Ends
admit the structure of a Frobenius algebra.
\end{Example}

\begin{Example}
Note that the subcategory $\widehat\calc^\tpr$ of admissible objects of $\calc$ does
not depend on the choice of GV-structure, since it only involves the right exact
tensor product $\tpr$. In the pivotal case the subcategory of objects 
$c\iN \widehat\calc^\tpr$ for which there exists an isomorphism $c\Cong c\tpR K$ is 
of interest as well: the objects in this subcategory are precisely those for which
$c\Cong \TI\tpL c$; they also satisfy $c\Cong S(c)$ and their internal End is a
Frobenius algebra. This subcategory is obviously
very sensitive to the choice of dualizing object \(K\). For an abelian
group \(A\) and choice of normalized abelian \(3\)-cocycle \((F,\Omega)\),
consider the braided monoidal category $\vect_{\ko,A}^{(F,\Omega)}$: objects of
$\vect_{\ko,A}^{(F,\Omega)}$ are finite dimensional \(A\)-graded \(\ko\)-vector
spaces (so the isomorphism classes of simple objects are in bijection with the 
elements of \(A\)), morphisms are grade-preserving linear maps, the
tensor functor is the standard tensor product of graded vector spaces, and
the braiding and associativity on simple objects are given by \(\Omega\)
and \(F\), respectively. As is well known, this category is rigid (hence
\(\vect_{\ko,A}^{(F,\Omega)} \eq \widehat{\vect}{}_{\ko,A}^{(F,\Omega),\tpr}\)) and
every simple object $\ko_h$ is invertible. Thus any \(\ko_h\), \(h\iN A\), can
be taken to be a dualizing object. If we choose \(\ko_0\) to be the
dualizing object, then for \(m\iN \vect_{\ko,A}^{(F,\Omega)}\) the condition
\(m\Cong S_0(m)\Cong m\oti \ko_0\) is empty. On the other hand, if we choose
\(\ko_h\) with \(h\,{\ne}\,0\) to be the dualizing object, then an isomorphism
\(m\Cong S_h(m)\Cong m\oti \ko_h\) can only exist if \(m\) is a direct sum of the
form \(\bigoplus_{\tilde h\in \langle h\rangle\!} \ko_{g+\tilde h}\), where
\(\langle h\rangle\) is the subgroup of \(A\) generated by \(h\). Note that
this sum is contained in \(\vect_{\ko,A}^{(F,\Omega)}\) only if \(h\) has
finite order, otherwise completions are required. So objects \(m\) admitting
isomorphisms \(m\Cong S_h(m)\) are sums of ``cosets of \(\langle h\rangle\)
in \(A\)''.
\end{Example}


\vfill

\noindent
{\sc Acknowledgments:}\\[.3em]
We thank Max Demirdilek and Lukas Woike for valuable comments.
J.F.\ is supported by the Swedish Research Council VR under project no.\ 2022-02931.
S.W.\ is supported by the Engineering and Physical Sciences Research Council (EPSRC)
EP/V053787/1 and by the Alexander von Humboldt Foundation.    
C.S.\ is supported by the Deutsche Forschungsgemeinschaft (DFG, German Research 
Foundation) by SFB 1624 and under Germany's Excellence Strategy - EXC 2121 
``Quantum Universe'' - 390833306.

\newpage

\newcommand\wb{\,\linebreak[0]} \def\wB {$\,$\wb}
\newcommand\Arep[2]  {{\em #2}, available at {\tt #1}}
\newcommand\Bach[2]  {{\em #2}, Bachelor thesis (#1)}
\newcommand\Bi[2]    {\bibitem[#2]{#1}}
\newcommand\inBO[9]  {{\em #9}, in:\ {\em #1}, {#2}\ ({#3}, {#4} {#5}), p.\ {#6--#7} {\tt [#8]}}
\newcommand\J[7]     {{\em #7}, {#1} {#2} ({#3}) {#4--#5} {{\tt [#6]}}}
\newcommand\JO[6]    {{\em #6}, {#1} {#2} ({#3}) {#4--#5} }
\newcommand\JP[7]    {{\em #7}, {#1} ({#3}) {{\tt [#6]}}}
\newcommand\Jpress[7]{{\em #7}, {#1} {} (in press) {} {{\tt [#6]}}}
\newcommand\Jtoa[7]  {{\em #7}, {#1} {} (to appear) {} {{\tt [#6]}}}
\newcommand\BOOK[4]  {{\em #1\/} ({#2}, {#3} {#4})}
\newcommand\Mast[2]  {{\em #2}, Master's thesis #1}
\newcommand\PhD[2]   {{\em #2}, Ph.D.\ thesis #1}
\newcommand\Prep[2]  {{\em #2}, preprint {\tt #1}}
\newcommand\Prev[2]  {{\em #2}, preprint version {\tt #1}}
\newcommand\uPrep[2] {{\em #2}, unpublished preprint {\tt #1}}
\def\adma  {Adv.\wb Math.}
\def\ajse  {Arabian Journal for Science and Engineering}
\def\coma  {Con\-temp.\wb Math.}
\def\comp  {Com\-mun.\wb Math.\wb Phys.}
\def\coma  {Con\-temp.\wb Math.}
\def\joal  {J.\wB Al\-ge\-bra}
\def\jopa  {J.\wb Phys.\ A}
\def\jpaa  {J.\wB Pure\wB Appl.\wb Alg.}
\def\kyjm  {Ky\-o\-to J.\ Math.}
\def\quto  {Quantum Topology}
\def\sema  {Selecta\wB Mathematica}
\def\taac  {Theo\-ry\wB and\wB Appl.\wb Cat.}
\def\tams  {Trans.\wb Amer.\wb Math.\wb Soc.}

\end{document}